\documentclass [a4paper, 10pt] {amsart} 
\usepackage {amssymb,amsmath,bbm,latexsym,setspace,parskip,txfonts,psfrag}
\usepackage {graphicx, epstopdf}
\input diagrams.sty

\newtheorem{thm}{Theorem}
\newtheorem{prop}[thm]{Proposition}
\newtheorem{lemma}[thm]{Lemma}
\newtheorem{cor}[thm]{Corollary}
\newtheorem{dfn}[thm]{Definition}

\newtheorem{conj}[thm]{Conjecture} 

\numberwithin{thm}{section}

\newcommand {\essup} {\mathop \mathrm{ess \; sup}}
\newcommand {\Teich} {\mathrm{Teich}}
\newcommand {\mcg} {\mathrm{Mod}}
\newcommand {\psl} {\mathrm{PSL}(2,\bb C)}
\newcommand {\Belt} {\mathrm{Belt}}
\newcommand {\out} [1] {{}}
\newcommand {\ssm} {\scriptscriptstyle}
\newcommand {\shb} {\stackrel{!}{=}}

\newcommand {\oh} {\mathcal{O}}

\newcommand {\bb} [1] {\mathbbm{#1}}
\newcommand {\cal} [1] {\mathcal{#1}}
\newcommand {\ol} [1] {\overline{#1}}
\newcommand {\ul} [1] {\underline{#1}}

\newcommand {\cinf} {\hat{\bb C}}

\newcommand {\MPIs} {Max Planck Institute for Mathematics in the Sciences}
\newcommand {\be}{\begin{equation}\begin{aligned}}
\newcommand {\ee}{\end{aligned}\end{equation}}
\newcommand {\benn}{\begin{equation*}\begin{aligned}}
\newcommand {\eenn}{\end{aligned}\end{equation*}}
\newcommand {\QC} {\cal {QC}}

\psfrag {inf}{$\infty$}
\psfrag {zero}{$0$}
\psfrag {TM-a} {$S^1 = \partial \bb D$}
\psfrag {TM-b} {$w^\mu(\partial \bb D)$}
\psfrag {TM-c} {$w^\mu$}
\psfrag {P_1} {$P_1$}
\psfrag {P_2} {$P_2$}
\psfrag {P_3} {$P_1$}
\psfrag {P_4} {$P_1$}
\psfrag {P} {$P$}
\psfrag {Lambda}{$\Lambda$}

\title[Families of Line Bundles over Riemann Surfaces - A Constructive Approach via Automorphic Forms]{Families of Line Bundles over Riemann Surfaces, their Sections,  and their Degenerations  -- A Constructive Approach using Automorphic Forms}
\author{Guy Buss} 
\address{Mathematisches Institut der Universit\"at Bonn\\ Endenicher Allee 60, 53115 Bonn, Germany\\guybuss@math.uni-bonn.de}
\date{\today}

\begin{document}
\begin{abstract}
		In this paper we study the deformation problem of pairs consisting of a Riemann surface and a holomorphic line 				bundle over that surface, and also sections thereof. We emphasize a constructive approach throughout and work and 		use covering space techniques. In particular, we also describe the limits of such degenerations as the 		boundary of Teichm\"uller space is approached, and review the construction of augmented Teichm\"uller space in 				great detail.		
\end{abstract}

\setcounter{tocdepth}{2}
\maketitle
\tableofcontents
\section*{Introduction}
This paper is an expository article on holomorphic line bundles over Riemann surfaces and sections of these, their deformations as the complex structure varies and their limits over degenerated Riemann surfaces. Our aim is not only to treat this as an abstract moduli problem, but emphasize explicit constructions whenever possible. Therefore we have chosen an approach using covering space techniques. On the other hand, our point of view ultimately is geometric, and because of this we pause and translate our results into the differential-geometric language as often as possible. Before we go into more detail we want to share our initial motivation for the present work so that the content can be appreciated by the reader with this perspective in mind.\\

\textbf{Motivation.} Motivated by the classical theory of harmonic mappings of Riemann surfaces into other manifolds and their behaviour under changes of the conformal structure of the domain surface, in particular under its degeneration, we started to investigate similar questions for other two-dimensional geometric PDEs arising from mappings of Riemann surfaces. Because symmetries of the model under consideration have a deep impact on the resulting theory, we were particularily interested in supersymmetric models, as well as a model very much related to supersymmetry recently introduced by J\"urgen Jost and his coauthors in \cite{D-HM}, the solutions of which they dubbed \emph{Dirac-harmonic} maps.\\

All these models have a new ingredient in common: They are models for \emph{two or more fields}, one of which is again the map $\phi:\Sigma \rightarrow M$, whereas the others are typically sections $\psi_i$ of vector bundles over $\Sigma$ of the form $\cal L_i \otimes \phi^{-1} T^*M$ arising from holomorphic line bundles $\cal L_i\rightarrow \Sigma$. Because the $\psi_i$ are sections of bundles which are twisted by the pull-back of the cotangent bundle of the target via the map $\phi$, the equations are coupled in a geometrically very natural way. \\

In order to study the dependence on the conformal structure of the domain surface in such a model, one needs to understand in what sense pairs $(\Sigma, \cal L \rightarrow \Sigma)$ vary when the conformal structure of $\Sigma$ is changed. More formally, one needs to determine the universal object in the category of families of marked Riemann surfaces equipped with holomorphic line bundles. The next natural step is to study these objects under degenerations of the underlying Riemann surfaces, i.e., the question is: What happens with the line bundles? And even more importantly, what happens to sections of the corresponding bundles?\\

\textbf{Present Work.} In this paper we present a unified framework for dealing with the questions above, i.e., what precisely are deformations of pairs ($\Sigma, \cal L \rightarrow \Sigma)$, how do we construct them, is there a universal object and can one construct (holomorphic) sections of these families. The paper is divided into three parts, which we describe in detail below. 

But before, let us say that many of these questions have been addressed in the literature before. In fact, the portion of results in this paper that are new\footnote{Actually, 'new' is a difficult concept, as there are many small results we were not able to find anywhere in the literature although we are certain they are known to the experts. That being said, the results on extensions of sections to the boundary of Teichm\"uller space do qualify as new in our opinion, in particular in the precise formulation given here.} is rather small. However, the answers needed for analytical considerations involving PDEs must be formulated in explicit terms in order to be of use. And while an inflation of mathematical knowledge over the last half a century has taken place, many of the already existing results are distributed over numerous research papers, as well as formulated in different terms and with very different foci. That made the need for a coherent presentation in uniform notation and language, and this paper is the result of that effort.\\

Let us now give a more detailed breakdown of the content of the paper.\\

{\bf The Absoulte}, the first part of the paper, deals with prerequisites: here we review and compile different ingredients that will be put together later. We start by reviewing the uniformization of Riemann surfaces, and also introduce some more general vocabulary of Kleinian groups, since we will need more than standard Fuchsian uniformization in the third part, when we address the question of degenerations and boundary behaviour of the constructions. \\

After that, in Section \ref{LBvsFoA}, we review holomorphic line bundles on surfaces and show that for fixed degree their moduli space is a complex torus of dimension equal to the genus of the surface. We pull back the line bundles to the uniformizing surface and describe them using factors of automorphy for the covering group. Sections of the line bundles lift to functions on the universal cover with symmetry properties under the covering group, called automorphic forms. For so-called $s$-factors of automorphy there are very powerful Banach spaces techniques available, which we present in Section \ref{BSoAF}, and to be able to utilize them we have to show that such an $s$-automorphic uniformization of any line bundle exists. This is the content of Theorem \ref{uf} and its corollaries.\\

We end the first part by some analytic considerations concerning the Poincar\'e density of the unit disc which are needed for the boundary behaviour in part three, and finally, in Section \ref{RSaAF} we comment more precisely on the connection of automorphic forms and sections of line bundles for more general, i.e., non-compact, Riemann surfaces.\\

{\bf The Relative}. This is where we start to deform the objects under consideration. After giving appropriate general definitions of deformations in Section \ref{GenFam}, we proceed by introducing Teichm\"uller space in Section \ref{Teich}. This space is of utmost importance for the present work as it is the base space of the universal object in the category of families of marked Riemann surfaces. It allows several natural descriptions, and yet another one is presented in Section \ref{Bersemb}. In fact, this latter point of view allows Teichm\"uller space to be seen as a subset of a bigger space, namely the deformation space of a Fuchsian group. Yet again, this is needed for the description of boundary points of Teichm\"uller space in the third part of the paper.\\

Section \ref{FibreoverTeich} then introduces the Bers fibre space over Teichm\"uller space, which can be thought of as the universal family of uniformizations of the surface, and the universal curve, a fibre space where the fibre over a point in Teichm\"uller space is the corresponding Riemann surface. This now is the universal object for families of marked Riemann surfaces.\\

After having introduced all these objects we show in Section \ref{faforms} how a line bundle on the base fibre, described via factor of automorphy, extends to a line bundle over Teichm\"uller space, and how the factor of automorphy for the action of the covering group on the disc extends to a factor for the action of the covering group on Bers fibre space. In fact, we show in Theorem \ref{UniqueBundle} that the factor extends uniquely in such a way that it induces an $s$-factor on each fibre - a crucial fact for all analytic constructions of deformations.\\

In the remaining sections of this part we now go into detail how one can extend automorphic forms on the base surface to automorphic forms for factors of automorphy on the Bers fibre space, both in a smooth and in a holomorphic way. We prove that for any given section on the base surface, a holomorphic extenesion exists (Thm. \ref{ExtendSections}) and give a recipe. In fact, such a family of sections of line bundles can be constructed from a \emph{single polynomial on the unit disc}. Also due to the fact that these sections can be extended holomorphically, this leads to the fact that the vector bundle of global holomorphic sections of families of line bundles can, via this explicit procedure, be trivialised almost everywhere, i.e., outside of an analytic subset of complex codimension at least one (Thm. \ref{BasisAE}).\\

{\bf The Frontier.} This is the part where it becomes really interesting. In general, boundary phenomena of moduli problems are often the most challenging part of the story, both technically and conceptually. For degenerations of Reimann surfaces there is so much interesting mathematics going on that we could not resist in digressing from our main goal, namely the extension of automorphic forms to the boundary, and first give the reader a little feeling by what is going on. This we do from several points of view in Section \ref{DegenerofRS} thru Section \ref{WPComp}. \\

Since for our approach, which relies heavily on the varying coverings of the surfaces in the family, we need a description of the boundary objects \emph{together with their uniformizations}. This is the reason why we work with so-called \emph{augmented Teichm\"uller space} - by construction the limit points are given by quotients of domains on the Riemann sphere by the action of Kleinian groups. More precisely, this space is obtained from the Bers image of Teichm\"uller space in the space of bounded quadratic differentials by adding all boundary points corresponding to regular $b$-groups. It is a beautiful idea, but requires some work to understand the induced geometry and in which sense these are the limit points of degenerating Riemann surfaces. All this we present in quite some detail starting in Section \ref{BersBound}, and we also give a precise relation of this space to other 'enlarged' Teichm\"uller spaces.

In Section \ref{ExtoATS}, now having the augmented Teichm\"uller space, i.e., an extension of the base space of the universal family, at our disposal, we start thinking about the spaces over Teichm\"uller space that were considered earlier. We propose a definition of the augmented Bers fibre space, which to the best of our knowledge has not appeared in the literature before, and prove that it is a connected bounded set in $\bb C^{3g-2+n}$ (where $(g,n)$ is the genus and the number of punctures of the base surface), and that the fibres are either simply-connected or consist of infinitely many simply-connected components. We further show how factors of automorphy naturally extend to this space, how the automorphic forms over limit points look like and prove that the families of automorphic forms constructed in Theorem \ref{ExtendSections} extend continuously to automorphic forms over augmented Teichm\"uller space (Thm. \ref{ExttoBound}).\\

Finally, in Section \ref{AsofSections} we briefly discuss the asymptotic growth with respect to the Weil-Petersson product of the extended sections and show how this is related to the norm of the inclusion operator between integrable and bounded automorphic forms for function groups.\\

\textbf{Acknowledgenments.}
This work is part of my Ph.D. thesis written under the supervision of J\"urgen Jost at the \MPIs{} in Leipzig. I thank him for his constant support and guiding. I thank Will Kirwin for proofreading and conceptual help in structuring the final version of this paper, as well as Brian Clarke and Christoph Sachse for many inspiring discussions. The work was funded both by the Research Training Group \emph{Analysis, Geometry and their Interaction with the Natural Sciences} at the University of Leipzig and the \emph{International Max Planck Research School} at the MPI Leipzig.
%%%%%%%%%%%%%%%%%%%%%%%%%%%%%%%%%%%%%%%%%%%%%%%%%%%%%%%%%%%%%%%%%%%%%%%%%%%%%%%%%%%555
\addtocontents{toc}{\vspace{.5cm}}
\section*{\bf The Absolute}
\addtocontents{toc}{\vspace{.3cm}}
\setcounter{section}{1}
This part will deal with the situation of a single Riemann surface. Since the concept of \emph{uniformization} plays a key role in the whole exposition, we start by a brief discussion of Kleinian groups in Section \ref{KGandU}. This is also necessary, since we later, in the third part of the paper, will discuss limits of deformations of Riemann surfaces from the point of view of uniformization, and there lesser known Kleinian groups, so-called regular $b$-groups become important.\\

In the subsection thereafter we discuss the moduli space of line bundles on a surface, and then show how line bundles are described using covering space techniques by factors of automorphy. The link to the next section is Theorem \ref{AllBundles}, showing that any line bundle can be represented by an $s$-factor.\\

We will then proceed to introduce powerful Banach space techniques for $s$-factors of automorphy in Section \ref{BSoAF}, which are absolutely crucial for the constructions to come in the relative setting.\\

After covering some more necessary analytic background in Section \ref{IotPD}, Section \ref{RSaAF} then contains a precise discussion of the geometric meaning of these Banach spaces.
%%%%%%%%%%%%%%%%%%%%%%%%%%%%%%%%%%%%%%%%%%%%%%%%%%%%%%%%%%%%%%%%%%%%%%%%%%%%%%%%%%%55
\subsection{Kleinian Groups and Uniformization}\label{KGandU}
Uniformization roughly means the same as an explicit realization of the universal covering in the same category, i.e., we want to consider Riemann surfaces $\Sigma$ as given by $\Sigma \cong D/G$, where $D$ is a simply-connected Riemann surface and $G\subset \mathrm{Aut}(D)$ is a subgroup of the biholomorphisms of $D$. By Riemanns mapping theorem, the only simply-connected Riemann surfaces (up to biholomorphism) are $\bb D, \bb C$ and $\cinf$. Their automorphism groups are given by
\begin{eqnarray*} &&\mathrm{Aut}(\cinf) = \mathrm{PSL}(2,\bb C) \qquad \qquad
   \mathrm{Aut}(\bb C) = \{(a,b) \in \bb C^*\times \bb C, z \mapsto az + b\}\\
      &&\mathrm{Aut}(\bb D) = \{(a,\xi) \in \bb D \times \partial \bb D, z \mapsto \xi (z-a)(\bar az-1)^{-1}\}\;.
\end{eqnarray*}
Now, the quotient of a Hausdorff topological space $X$ by a continuous group action is again Hausdorff, iff the action is proper discontinuous, and if in addition $X$ is a manifold, the quotient is again a manifold, iff the action is freely discontinuous\footnote{$G$ acts properly discontinuously on $A \subset X$, iff any $x$ has a neighborhood $U_x$ such that the number $N_x:= \{g \in G| g(U_x) \cap U_x \neq \emptyset\}$ is finite. It acts freely discontinuously on $A \subset X$, iff  any $x$ has a neighborhood $U_x$ such that $N_x = 1$.}.\\

It turns out that there exists no non-trivial group $G\subset \mathrm{PSL}(2, \bb C)$ that acts freely discontinuously on all of $\cinf$, i.e., the only Riemann surface covered by $\cinf$ is $\cinf$ itself. The Riemann surfaces covered by $\bb C$ are $\bb C, \bb C \backslash 0$ and tori. These are called \emph{exceptional Riemann surfaces}.\\

From now on we will only consider non-exceptional Riemann surfaces, which we also call hyperbolic Riemann surfaces since due to being uniformized by the unit disc and the fact that automorphisms of the disc are hyperbolic isometries they can be equipped with a hyperbolic metric. It is clear by the previous discussion, that as long as a single Riemann surface $\Sigma$ is studied, it is sufficient to uniformize $\Sigma \cong \bb D/G$, with $G \subset \mathrm{Aut}(\bb D)$. However, as we deform the surface, a holomorphic theory of families of Riemann surfaces requires more technical flexiblity as we will see in Section \ref{Teich}.  
%%%%%%%%%%%%%%%%%%%%%%%%%%%%%%%%%%%%%%%%%%%%%%%%%%%%%%%%%%%%%%%%%%%%%%%%%%%%%%%%%%%%%%%%%%%%%%%%%%%%%%555
\subsubsection{From General Kleinian Groups to b-groups}
If no other reference is mentioned, the monograph \cite{KG} by Bernard Maskit is used as a basic reference for this Section. Let $G\subset$ PSL$(2, \bb C)$ be any subgroup. We define
the following subsets of the Riemann sphere,
\begin{eqnarray*}\label{domofKGDef}
^\circ\Omega(G)&:=&\{ z \in \cinf: G \textrm{ acts freely discontinuously at $z$}\}\\
      \Omega(G)&:=&\{ z \in \cinf: G \textrm{ acts properly discontinuously at      $z$}\}\\
      \Lambda(G)&:=&\{ z \in \cinf: \exists \; \textrm{a sequence of distinct } g_m \in G \textrm{ and a point } w\in \cinf \textrm{ s.t. }
      g_m(w) \rightarrow z\}\;.
\end{eqnarray*}
$\Omega(G)$ is called the \emph{region of discontinuity} of $G$ and
$\Lambda(G)$ the \emph{limit set} of $G$. It follows (see \cite{KG},
Ch.~ II), that $\Omega(G)$ is either empty or open and dense in
$\cinf$, $\Omega(G) \backslash ^\circ\Omega(G)$ is precisely the set of fixed points of elliptic elements, and it is discrete. $\Lambda(G)$ either consists of $0,
1, 2$ or uncountably many points. In the latter case,
$\Lambda(G)$ is a perfect set with empty interior, i.e., any point of
$\Lambda(G)$ is an accumulation point of $\Lambda(G)$. Moreover
$\Lambda(G)\cap\Omega(G) = \emptyset$ and $\Lambda(G)\cup\Omega(G) =
\cinf$.
\begin{dfn}
   A subgroup $G \subset \mathrm{PSL}(2,\bb C)$ is called \underline{Kleinian} iff $\Omega(G) \neq \emptyset$. A Kleinian group is called
   \underline{elementary} iff $\sharp \Lambda(G) \in \{0,1,2\}$, else it is called \underline{non-elementary}.
\end{dfn}
A Kleinian group is discrete in the natural quotient topology which $\mathrm{PSL}(2,\bb C)$ inherits 
from Mat$(2\times 2, \bb C) \cong \bb C^4$. This follows
from $\Omega(G) \neq \emptyset$. Denote the connected components of
$\Omega(G)$ by $D_i$, where $i \in I$, some index set. It is necessarily countable since $G$ is countable by discreteness, and will often implicity be assumed to be a subset of the natural numbers. Let $G$ act
on each $D_i$. Then any $\gamma \in G$ maps a connected component to
another one, $\gamma(D_i) = D_j$. If such a $\gamma$ exists for
given $i,j$, we say $D_i$ is \emph{conjugate} to $D_j$. The set of
elements $\gamma$ such that $\gamma(D_i) = D_i$ form a subgroup
$G_{D_i} \subset G$, the \emph{stabilizer} of $D_i$. Let
$$ I_i:= \left \{ j \mid \: \bigcap G(D_i) \cap D_j \neq \emptyset \right \}\;, \qquad m_i:= \mathrm{min} \{ j \mid j \in I_i\}\;. $$
Clearly, as sets $I_i = I_j$ if $D_i$ and $D_j$ are conjugate, and else
$I_i \cap I_j = \emptyset$. The component $D_{m_i}$ is a choice of a representative
for the class of conjugate domains and the new index set consisting
of all these representatives will be denoted by $K$. Now let
$$ \Delta_{k} := \bigcup_{i \in I_{k}} D_i\;.$$
The $\Delta_k$, called the \emph{components} of $G$\footnote{The
$\Delta_k$ are called components \emph{of G} in order to avoid
confusing them with the connected components $D_i$ of $\Omega(G)$.},
are clearly $G-$invariant and partition $\Omega(G)$. If a connected
component of $\Omega(G)$ already happens to be a component of $G$,
i.e., $D_i = \Delta_k$
(and then of course $i=k$), it is called an \emph{invariant component}. \\

In terms of the orbits, $\Delta_k / G \cong D_k/G_{D_k}$, so when we are 
interested in the quotient, we can always assume the components of a
group to be connected. Altogether
\begin{eqnarray}\label{decompo}
    \Omega(G)/G \cong \coprod_{k \in K} D_k / G_{D_k} \cong \coprod_{k \in K} \Delta_k / G\;,
\end{eqnarray}
and we denote by $X_k := D_k / G_{D_k}$ the connected components of
the quotient, each of which is a (possibly ramified) Riemann surface,
depending on whether the action of $G_k$ on $D_k$ has fixed points.
\begin{dfn}
   Let G be a non-elementary Kleinian group. G is called of \underline{finite analytic type} iff
   the Poincar\'e area of $\Omega(G)/G$ is finite.
\end{dfn}
So $G$ being of finite analytic type means that all $X_k$ occuring
in the disjoint decomposition of $\Omega(G)/G$ are Riemann
surfaces of finite analytic type. Moreover, the Poincar\'e area of a
Riemann surface of finite analytic type is given by
\begin{eqnarray} \label{VolRS}
   \mathrm{Area}(D_k/G_k) \geq 4\pi(g-1 + r)\geq 4\pi\;,
\end{eqnarray}
where $r\geq 0$ is a term coming from the ramification points, but
since we will only consider unramified Riemann surfaces in the
sequel, we skip the detailed expression. Therefore the
decomposition (\ref{decompo}) \emph{can only be finite}. With this
in mind, the following famous theorem by Ahlfors \cite{FGKG} can be
fully appreciated.
\begin{thm}[Ahlfors Finiteness Theorem] \label{AFT}
    Each finitely generated non-elementary Kleinian group is of finite analytic type.
\end{thm}
We will have to further distinguish Kleinian groups.
\begin{dfn}
    Let G be a non-elementary Kleinian group.
    \begin{itemize}
    \item{G is called a \underline{function group} iff it has an invariant
    component $\Delta$ such that $\Delta/G$ is of finite analytic
    type.}
    \item{A function group is called a \underline{b-group} iff the invariant component $\Delta$ is
    simply connected.}
    \item{A b-group  is called \underline{degenerate} iff the invariant component coincides with the region of discontinuity, i.e. $\Delta =
    \Omega(G)$.}
    \end{itemize}
\end{dfn}
By the Ahlfors Finiteness Theorem, the condition that $\Delta/G$ be of
finite type is automatically satisfied for finitely generated $G$.
Conversely, a function group is finitely generated (see \cite{KG},
Ch. X). This is a very remarkable fact, since in general a group
which is of finite analytic type need not be finitely generated,
i.e., the converse to Ahlfors theorem does not hold.\\
\begin{dfn} \label{Fuchs}
   A non-elementary Kleinian group $G$ is called \underline{Fuchsian} iff it leaves a round disc on $\cinf$ invariant.
   A Fuchsian group is called of \underline{first kind} iff the limit set is the whole boundary of the disc.
   Else $G$ is called of \underline{second kind}.
\end{dfn}
Obviously, if a group leaves a disc $D$ invariant, it also leaves the disc $D^c$ and the common boundary $\partial D$ invariant. It is easy to see that if there is a closed set invariant under the group action, then it contains the limit set. Hence the definition of first and second kind above makes sense because $\Lambda(G) \subset \partial D$ by the argument just given.\\

By a quasiconformal deformation of a group $G$ we mean a group resulting from conjugating $G$ by a quasiconformal homeomorphism of $\cinf$.
\begin{dfn}
    A group $G\subset \mathrm{PSL}(2,\bb C)$ is called \underline{quasi-Fuchsian} iff it is a quasiconformal deformation of a Fuchsian group. 
\end{dfn}
In the same way as above, a quasi-Fuchsian group is called of \emph{first} resp.~ \emph{second} kind iff the group out of which it was deformed was of first resp. second kind. \\

The region of discontinuity of a quasi-Fuchsian group has two
components if it is of first kind and one component if it is of
second kind. Observe that by Theorem \ref{AFT}, finitely generated quasi-Fuchsian groups
are function groups, and if they are of first kind they are even b-groups\footnote{Some authors exclude quasi-Fuchsian groups when they talk about b-groups, which is probably more in the spirit of thir name: The $b$ stands for \emph{boundary} because many b-groups lie on the boundary of Teichm\"uller spaces (Bers even conjectured, that all b-groups are boundary groups, (see Sect. \ref{BersBound}) while quasi-Fuchsian groups are interior points of Teichm\"uller spaces. Other authors distinguish between $b$- and $B$-groups. We will stick to the definition given above which agrees with the one given in \cite{KG}.}. \\

The key concept for studying the geometry of b-groups is the notion
of the Fuchsian equivalent.
\begin{dfn}
Let $(G, \Delta)$ be a b-group and $\psi: \bb D \rightarrow D$ a
Riemann mapping. The group $\Gamma:= \psi G \psi^{-1}$ is called a
\underline{Fuchsian equivalent} of $G$.
\end{dfn}
Obviously, a Fuchsian equivalent $\Gamma$ of $(G,\Delta)$ is a
finitely generated Fuchsian group of first kind and comes equipped
with a canonical isomorphism
$$\chi:G \rightarrow \Gamma\;, \qquad \chi(g):= \psi \circ g \circ
\psi^{-1}.$$ Of course $\chi$ and its inverse preserve elliptic
fixed points, since they lie in the interior of $\Delta$ resp. $\bb D$.
The corresponding statement for other elements is the following (see
\cite{OBoTSaoKG:1}, Sect.~ 5, Prop.~ I-III).
\begin{prop}
   The image under $\chi$ of a loxodromic element is always loxodromic while the image of a
   parabolic element of $G$ may be parabolic or loxodromic in
   $\Gamma$.
\end{prop}
\emph{Idea of proof.} Bers proves that the map $\psi:\bb D
\rightarrow \Delta$ can be extended to a map $\tilde \psi: \bb D^*
\rightarrow \Delta^*$, where $\bb D^*$ is the union of $\bb D$ and
the fixed points of elements of $\Gamma$ (and similarily for $\Delta^*$)
and that this map is a homeomorphism when restricted to the ends of
the axes of the transformations. But it needn't be globally
injective. However, he shows that the only way non-injectivity comes
about is the following: If $\tilde \psi(x_1)=\tilde \psi(x_2)$, then
the $x_i$ are fixed points of the same element $\gamma \in \Gamma$,
neccesarily loxodromic, and $\chi^{-1}(\gamma)$ is parabolic. $\Box$
\begin{dfn}
    Let $G$ be a b-group and $\Gamma$ its Fuchsian equivalent. An
    element $g \in G$ is called an \underline{accidental parabolic
    transformation} (abbreviated APT) iff $g$ is parabolic
    and $\chi (g)$ is loxodromic.
\end{dfn}
We have so far introduced some basic vocabulary and facts on Kleinian groups. We will investigate the geometry of b-groups more deeply in the third part of this paper. For this we still need some more general facts about the domain $\Omega(G)$ and the limit set $\Lambda(G)$, which we compile in the next section.
%----------------------------------------------------------------------
\subsubsection{The Geometry of $\Omega(G)$ and $\Lambda(G)$}
The two sets $\Omega(G)$ and $\Lambda(G)$ are geometrically very interesting sets and reflect the geometric complexity of Kleinian groups mentioned at the beginning of the section. We start by focussing on $\Omega(G)$ at first, though of course any statement about this set yields some information about $\Lambda(G)$ because $\Lambda(G) = \cinf \backslash \Omega(G)$. We list a few propositions, all contained in the monograph \cite{KG}, without further comment, and refer to the numbering in \cite{KG}.
% Every non-elementary Kleinian group contains a loxodromic element [V.E.1].
\begin{prop}[V.E.7] \label{2invcomp}
   Let $G$ be a non-elementary Kleinian group with $n\geq 2$
   invariant components. Then $n=2$.
\end{prop}
\begin{prop}[IX.D.5] Let $G$ be an analytically finite Kleinian
group with two invariant components $\Delta_1, \Delta_2$. Then
$\Omega(G) = \Delta_1 \cup \Delta_2$.
\end{prop}
\begin{prop}[V.E.8]
  Let $G$ be Kleinian with component $D$. Then either
  \begin{itemize}
     \item{$D$ is an invariant component,}
     \item{there is another component $D'$ such that $D\cup D'$ is
     invariant or}
     \item{there are infinitely many distinct translates of $D$.}
  \end{itemize}
\end{prop}
\begin{prop}[V.E.9]\label{finimp2comp}
    Let $G$ be a Kleinian group with finitely many components. Then
    $G$ has at most two components.
\end{prop}
These will be used in combination to understand the domains of discontinuity of b-groups on the limit set of Teichm\"uller space in Section \ref{BersBound}.
%%%%%%%%%%%%%%%%%%%%%%%%%%%%%%%%%%%%%%%%%%%%%%%%%%%%%%%%%%%%%%%%%%%%%%%%%%%%%%%%%%%%%%%
\subsection{Line Bundles vs. Factors of Automorphy}\label{LBvsFoA}
Sections of holomorphic vector bundles are interesting objects in complex analysis. For a fixed vector bundle $\cal V \rightarrow X$ they form a vector space themselves, denoted by $H^0(X,\cal V)$, whose dimension is usually denoted by $h^0(X,\cal V)$. In case $X$ is a compact manifold, these spaces are \emph{finite dimensional} due to the elliptic nature of the $\bar \partial$-operator\footnote{In fact, this is only a special instance of a more general result, namely that the space of global sections of a coherent sheaf on a compact complex space is finite dimensional.}.\\

For a line bundle over a Riemann surface there is a beautiful formula that relates the dimensions of certain spaces of holomorphic sections, the Riemann-Roch formula (see, e.g., \cite{ACaRS} or \cite{GoAC}).
\begin{thm} [Riemann-Roch] \label{R-R}
    Let $\cal L\rightarrow \Sigma$ be a line bundle over a compact Riemann surface $\Sigma$ of genus $g\geq 2$ and let $\cal K$ denote the canonical bundle of $\Sigma$. Then 
     \begin{eqnarray}h^0(\Sigma,\cal L) - h^0(\Sigma,\cal K \otimes {\cal L^{-1}}) = \mathrm{deg}(\cal L) + 1 - g\;.\end{eqnarray}
\end{thm}
This formula has several consequences. Recall that it is an elementary consideration to find out that bundles of negative degree have no holomorphic sections. 
\begin{itemize}
   \item{By applying Riemann-Roch to $\cal L \cong \cal O(\Sigma)$ one gets $h^0(\Sigma,\cal K) = g$.}
   \item{Likewise, by applying Riemann-Roch to $\cal L\cong \cal K$ one gets $\mathrm{deg}(\cal K)=2g-2$.}
   \item{For line bundles $\cal L$ with $\mathrm{deg}(\cal L) \geq 2g-1$, the Riemann-Roch formula gives an explicit formula for the dimension of the holomorphic sections,
   $$  h^0(\Sigma,\cal L) = \mathrm{deg}(\cal L) + 1 - g\;.$$}
\end{itemize}   
Bundles of degree greater than or equal to $2g-1$ are said to be in the \emph{topological range}, because the dimension of the space of global holomorphic sections only depends on the genus of the surface. If $\cal L$ is not in the topological range, the Riemann-Roch formula does not give the dimension explicitly, and moreover, the dimension is \emph{not topological}, i.e., the dimension possibly jumps when the complex structure is varied. This leads to interesting complex subvarieties of Teichm\"uller and moduli space (see, e.g., \cite{GoAC}).\\
\subsubsection{Moduli space of Line Bundles on a Riemann Surface}
Let us briefly describe the space of all holomorphic line bundles on a given Riemann surface. Recall that the group of isomorphism classes of complex line bundles on a compact Riemann surface is given by $\bb Z$, the integer being the degree of the bundle. On the other hand, this does not classify holomorphic line bundles - for each fixed degree we in fact get a $2g$-dimensional torus as we will see in a moment.\\

For any complex manifold, we have the following exact sequence of
sheaves\footnote{Here we use the standard notation $\oh_X$, resp.~ $\oh^*_X$, for the sheaf of holomorphic, resp.~ holomorphic nowhere vanishing functions. Further, a underlined ring means that we consider the locally constant sheaf, the stalks of which are precisely this ring.}
$$ 0 \rightarrow \underline{\bb Z} \stackrel{\iota}{\longrightarrow} \oh_X \stackrel{\ssm \textrm{exp}(2\pi i \cdot)}{\longrightarrow} \oh_X^* \rightarrow 0\;,$$
which, as usual, induces a long exact sequence in Cech cohomology. A part of it goes as follows:
$$\ldots \rightarrow \check H^1(X,\underline{\bb Z}) \stackrel{\iota_*}{\longrightarrow} \check H^1(X, \oh_X) \rightarrow \check H^1(X,\oh_X^*) \stackrel{c}{\longrightarrow} \check H^2(X,\underline{\bb Z})\rightarrow \ldots\;.$$
The connecting homomorphism $c$ is precisely the map sending a line bundle to its Chern class via the identification
$\check H^2(X,\underline{\bb Z}) \cong \bb Z$ for $X$ connected, which is the same as the degree of the associated
divisor (class). So looking at bundles of degree zero means looking at the kernel of $c$, and because of exactness,
$$ \textrm{ker} (c) \cong \textrm{im} (\textrm{exp}_*), \qquad  \textrm{ker} (\textrm{exp}_*) \cong \textrm{im} (\iota_*)\;,$$
we see that the line bundles of degree zero, denoted by Pic$^0(X)$, are given by
$$ \textrm{Pic}^0(X) \cong \check H^1(X, \oh_X) / \iota^* ( \check H^1(X, \underline{\bb Z})) \;.$$
Now, assume further that $X$ is compact. Since the sheaf $\oh_X$ is a sheaf of $\underline{\bb C}$-modules, $\check
H^1(X, \oh_X)$ is a complex vector space, and since $\oh_X$ is coherent, it is finite dimensional. Therefore Pic$^0(X)$
is a complex torus. Moreover, by Dolbeault's Theorem,
$$ \check H^q(X, \oh_X) \cong H^{0,q}_{\bar \partial}(X)\;.$$
Now, if $X$ is a compact Riemann Surface of genus $g$, $X$ is in particular a K\"ahler manifold, so the Dolbeault groups with
interchanged indices are isomorphic. Locally $dz$ spans the $(1,0)$-forms as a $C^\infty$-module; such a form is
$\bar \partial$-closed iff it is holomorphic and exact iff it is identically zero. By Riemann-Roch, the complex dimension
of this space is $g$, so
$$ \bb C^g \cong \Omega^{(1,0)}_{\textrm{hol}}(X) \cong   H^{0,1}_{\bar \partial}(X) \cong \check H^1(X, \oh_X)\;.$$
Hence, Pic$^0(X)$ in this case is a $g$-dimensional torus. \\

The full Picard group consists of bundles of arbitrary degree. They are now easily described, since $c:
\textrm{Pic}(X)\rightarrow \bb Z$ is a group homomorphism with kernel Pic$^0(x)$. The set of bundles of fixed
degree $q$, denoted (of course) by Pic$^q(X) \in \backslash \textrm{Pic}^0(X)$, is just a coset
space and so we have determined $\textrm{Pic}(X)$ fully as $\bb Z$ copies of $\textrm{Pic}^0(X)$.
%%%%%%%%%%%%%%%%%%%%%%%%%%%%%%%%%%%%%%%%%%%%%%%%%%%%%%%%%%%%%%%55
\subsubsection{Factors of Automorphy and Automorphic Forms}
Our next task is to connect line bundles with covering space techniques. Sections of line bundles are locally given by functions which are related via transition functions. Looking at this from the point of view of a covering, the lifted local functions must be related to each other via a periodicity condition encoded in the transition functions, and this periodicity condition is contained in the so-called factors of automorphy. In fact, using the precise definitions given below, it is not hard to work out that
\benn
		\left \{ \textrm{Iso.-classes of line bundles on}\; \Sigma\cong D / G \right\} \quad
		\leftrightarrow \quad \left \{ \textrm{factors of automorphy for the action of} \; G \; \textrm{on} \; D\right \}/\sim\;.
\eenn
Given a factor $\rho$, let us write $\cal L(\rho)$ for the corresponding line bundle (and vice versa). Furthermore, one also has the correspondence
\benn
   \left\{\textrm{Automorphic forms for}\; \rho\right \} \quad
		\leftrightarrow \quad \left \{\textrm{Sections of}\; \cal L(\rho) \right \}
\eenn
For an explicit account on these correspondences, see \cite{RSaGTF}.\\

Let us now come to the definitions as promised earlier.
\begin{dfn}\label{FactorOfAutomorphyDef}
    Let $D$ be a domain in $\cinf$ and $G \subset \mathrm{Aut}(D)$ a subgroup. A family of holomorphic functions $\rho_g: D \rightarrow \bb
C^*$, one for each $g \in G$, satisfying
\begin{eqnarray} \label{cocyc}
    \rho_{g_1g_2}(z) = \rho_{g_1} (g_2 z)
    \rho_{g_2}(z)\;, \qquad \forall z\in D\;,
\end{eqnarray}
is called a \underline{factor of automorphy} for the action of $G$ on $D$.
\end{dfn}
\begin{dfn}
   Two factors of automorphy $\rho, \eta$ are called \underline{equivalent} iff there exists a holomorphic function $h:D\rightarrow \bb C^*$ such
   that
   $$ \rho_g(z) = h(g z) \eta_g(z)h^{-1}(z)\;.$$
\end{dfn}
There are special kinds of factors of automorphy. The following classification is very useful.
\begin{dfn} \label{s-fact}
A factor of automorphy is called \underline{flat} iff each $\rho_g$ is constant and it is called
\underline{unitary flat} iff the constants are of modulus 1. Moreover, factors of automorphy satisfying
$$ |\rho_{g}(z)| = |g'(z)|^{-s} \qquad \forall \;
g \in G$$ for a fixed $s \in \bb R$ are called \underline {s-factors}.
\end{dfn}
Obviously, unitary flat factors are 0-factors, and conversely if a factor is an s-factor and a flat factor at once it is
unitary flat. Observe, that flat factors are nothing else than group homomorphisms $G \rightarrow \bb C^*$ and
unitary flat factors are group homomorphisms $G \rightarrow U(1) \subset \bb C^*$. Such homomorphisms are
determined by the images of a set of generators ${\gamma_i}$. In case $G$ is a cocompact Fuchsian group with standard generators, the images can be chosen freely, since the only relation in the group,
$$G \cong \left< a_1, \ldots, a_g, b_1, \ldots, b_g, \mid \prod_{i=1}^g a_ib_i a_i^{-1}b_i^{-1}= \bb 1_G \right >\;,$$
disappears automatically in case the target group is abelian.
\begin{dfn}\label{Dfn:AutomorphicForm}
   Let $\rho$ be a factor of automorphy. A function $f: D \rightarrow \bb C$ satisfying
   $$ f(g z) = \rho_{g} (z) f(z)\;, \qquad \forall
\gamma \in \Gamma_D\;,$$ is called an \underline{automorphic form} for the pair $(G,\rho)$. An automorphic form for an
s-factor will be called \underline{s-automorphic}.
\end{dfn}
Our next goal is to find a convenient representative in a given equivalence class of a factor of automorphy. We do this by adapting a theorem given by Gunning (\cite{RSaGTF}, Ch.~ II, Thm.~ 4). However, we use a different representative: an s-factor. This is crucial for the analytic construction of automorphic forms later on. We do this first for flat factors; the general case follows easily afterwards.
\begin{thm} \label{uf}
   Let $\Sigma\cong D/G$ be a compact Riemann surface and let $\rho \in \mathrm{Hom}(G, \bb C^*)$ be a flat factor on 		 $D$ for the group $G$. Then there exists a unitary flat factor equivalent to $\rho$.
\end{thm}
\emph{Proof.} Let $\{\omega_i \mid \;i=1 \ldots g\}$ be the basis of the vector space of automorphic forms for the canonical factor which corresponds to the basis for the holomorphic differentials on the Riemann surface $\Sigma$ dual to the $a$-generators of the symplectic basis of $H_1(\Sigma, \bb Z) \cong \langle a_1, \ldots, a_g, b_1, \ldots, b_g\rangle =:\langle \gamma_J\rangle$. We choose any point $z_0$ on $D$ and define the integral of $\omega_i$ by
$$ w_i(z) := \int_{z_0}^z \omega_i(\tau) d \tau\;.$$
The $\omega_i$ induce a map on the symplectic generators of $G$ by
$$\omega_i(\gamma_J) := w_i(\gamma_J z_0) =\int_{z_0}^{\gamma_{\ssm J} z_0} \omega_i(\tau)
d\tau\;.$$ Of course, the complex numbers $\omega_i(\gamma_J)$ are precisely the Riemann period matrix. The choice of the dual basis implies
$$\omega_i(a_j) = \delta_{ij}\;, \qquad \omega_i(b_j) = \tau_{ij}\;.$$
Now let $\rho \in \textrm{Hom}(\Gamma, \bb C^*)$ be given. For convenience, rewrite the numbers $\rho(\gamma_I)$
$$ \rho_{a_i} =: e^{2 \pi i \sigma_i}\;, \qquad \rho_{b_j} =: e^{2 \pi i \sigma'_j}\;. $$
Define a function $f$ by
 $$f(z):= \sum_{ij} C^{ji}\sigma_i w_j(z)\;.$$
This function satisfies
$$ f(a_k z) = f(z) + \sum_i C^{ki} \sigma_i\;, \qquad f(b_j z) = f(z) + \sum_{ik} C^{ki}\sigma_i
\tau_{kj}\;,$$ and moreover $f$ is holomorphic in $z$, so it's exponential $ h(z):= \textrm{exp}(2 \pi i f(z))$ is a
nowhere vanishing holomorphic function. We use this function to define an equivalent factor $\tilde \rho_\gamma:= h(\gamma
z) \rho_\gamma h^{-1}(z)$. Unitary flatness of $\tilde \rho$ translates into the $2g$ real equations
\begin{eqnarray}\label{Thenewcoc} |\rho_{\gamma_I}|^{-1} \shb \left|\frac{h(\gamma_I z)}{h(z)}\right|\;,\end{eqnarray}
which can be reduced to the following:
\begin{eqnarray*}
0&=&\textrm{Im}\left( \sum C^{ki}\sigma_i + \sigma_k \right) \\
0&=&\textrm{Im}\left( \sum C^{ji}\sigma_i \tau_{jk}+ \sigma'_k\right) \;.
\end{eqnarray*}
For the rest of the proof, the sum symbol will be omitted, that is, summation is assumed if two \emph{latin} indices appear in a product. Of course, if the sum is restricted, as in the next equation, it is written explicitly. \\

But now let us solve the equations. If all $\sigma_i, \sigma_j'$ vanish, the system is trivially satisfied. Suppose first $\sigma_\alpha \neq 0$ for some
$\alpha \in \{1,\ldots g\}$. Then the first $g$ equations, rewritten as
\begin{eqnarray} \label{1stEq} \textrm{Im}( C^{k\alpha}\sigma_{\alpha}) = - \textrm{Im}\left(\sum_{i \neq \alpha} C^{ki}\sigma_i + \sigma_k \right)\;,\end{eqnarray}
determine the entries of the $\alpha$-th column up to real-valued rescaling, thereby solving the equations. Now
suppose further that $\sigma_\beta \neq 0$ for some $\beta \neq \alpha$. Then the same trick works for the second $g$
equations since $\tau$ is non-singular,
$$ \textrm{Im}( C^{j\beta}\tau_{jk}\sigma_{\beta}) = - \textrm{Im}\left(\sum_{i \neq \beta} C^{ji}\tau_{jk}\sigma_i + \sigma_k \right)\;,$$
and the solution fixes the $\beta$-th column entries of $\tilde C^{jk} = C^{ji}\tau_{ik}$. \\

When such a $\beta$ does not exist, i.e., $\sigma_i = 0$ for all $i\neq \alpha$, we have to go back to the solution of the first $g$ equations and be a little more careful. $\sigma_\alpha \neq 0$ in particular implies either $\mathrm{Im} \sigma_\alpha \neq 0 $ or $\mathrm{Re} \sigma_\alpha \neq 0$. Assume $\mathrm{Im} \sigma_\alpha \neq 0$. Then we can rewrite \eqref{1stEq},
\begin{equation}\begin{aligned} \label{1stVar}\textrm{Re} (C^{k\alpha}) &= \frac{-1}{\textrm{Im} (\sigma_{\alpha})}\textrm{Im} \left(\sum_{i \neq \alpha} C^{ki}\sigma_i + \sigma_k  +  C^{k\alpha}\textrm{Re} (\sigma_{\alpha})\right)\\ &= \frac{-1}{\textrm{Im} (\sigma_{\alpha})}\left( \mathrm{Im} (\sigma_k)  + \mathrm{Im}(C^{k\alpha})\textrm{Re} (\sigma_{\alpha})\right)\;,\end{aligned}\end{equation}
which fixes $\textrm{Re} (C^{k\alpha})$ completely in terms of $\textrm{Im} (C^{k\alpha})$. We rewrite the second set of equations,
\begin{equation} \nonumber
   -\mathrm{Im} (\sigma'_k) = \mathrm{Im} \left( C^{j\alpha} \sigma_\alpha \tau_{jk}\right) = \mathrm{Im}( C^{j\alpha}\sigma_\alpha)\mathrm{Re} (\tau_{jk}) + \mathrm{Re}(C^{j\alpha}\sigma_\alpha)\mathrm{Im}(\tau_{jk})\;,
\end{equation}
and now use the first equations at two places. First we use $\mathrm{Im}( C^{j\alpha}\sigma_\alpha) = - \mathrm{Im} (\sigma_j)$, and then we insert the left-hand side of \eqref{1stVar} for $\mathrm{Re}(C^{k\alpha})$ after expressing the product $\mathrm{Re}(C^{j\alpha}\sigma_\alpha)$ in terms of the real and imaginary parts of the factors,
\begin{equation}\begin{aligned}\nonumber
   -\mathrm{Im} (\sigma'_k) &+ \mathrm{Im} (\sigma_j) \mathrm{Re} (\tau_{jk}) = \mathrm{Im}(\tau_{jk})\left(\mathrm{Re}(C^{j\alpha})\mathrm{Re}(\sigma_\alpha) - \mathrm{Im}(C^{j\alpha})\mathrm{Im}(\sigma_\alpha)\right) \\
    &= \mathrm{Im}(\tau_{jk})\left(\frac{-1}{\textrm{Im} (\sigma_{\alpha})}\left( \mathrm{Im} (\sigma_k)  + \mathrm{Im}(C^{k\alpha})\textrm{Re} (\sigma_{\alpha})\right)  \mathrm{Re}(\sigma_\alpha)- \mathrm{Im}(C^{j\alpha})\mathrm{Im}(\sigma_\alpha)\right)\;. \\
\end{aligned}\end{equation}
We view the left hand side of the equation as the inhomogeneous term of a linear system of equations and just write $B_k$ instead of the full expression. Observe that the first term on the right-hand side is also independent of our variables $\mathrm{Im}C^{k\alpha}$, so it is absorbed in the inhomogenity. We arrive at the following equations,
$$ B_k = \mathrm{Im}(\tau_{jk})\left(\frac{\mathrm{Re}^2(\sigma_\alpha) + \mathrm{Im}^2(\sigma_\alpha)}{\mathrm{Im}(\sigma_\alpha)}\right)  \mathrm{Im}(C^{j\alpha})\;,$$
and this system finally is seen to be solvable, because $\mathrm{Im}(\tau)$ is positive definite by Riemann's bilinear relations and the term in parentheses is a non zero number. A similar argument is possible for the case $\mathrm{Re}(\sigma_\alpha)\neq 0$. And finally, if all the $\sigma_i$ vanish, we redefine the function $f$ resp. $h$ by switching the role of the $\sigma$'s and $\sigma'$'. \\

So now equation \eqref{Thenewcoc} holds for all generators of the group, hence for all elements of the group and we have established $\tilde \rho$ as a unitary flat factor equivalent to $\rho$. $\Box$
\begin{cor}
    All bundles of degree zero can be expressed in terms of unitary flat factors
\end{cor}
\emph{Proof.} This follows immediately from Theorem \ref{uf}, since as is shown in \cite{RSaGTF}, all bundles of degree
zero can be expressed in terms of flat factors. $\Box$
\begin{cor} \label{AllBundles}
    All line bundles can be expressed in terms of s-factors.
\end{cor}
\emph{Proof.} Any line bundle of degree $d$ can be written as a product of a fixed line bundle of degree $d$ and a line
bundle of degree zero, so in order to prove the statement we only have to represent, for each integer $d$, a single
bundle of degree $d$ by an s-factor. The canonical bundle (see next section) is such a gadget, given by $\rho_\gamma(z) = \gamma'(z)$. Since the domain of the functions $D$ is simply connected, there exists a logarithm and hence any rational power $\rho^q_\gamma(z) := (\rho_\gamma(z))^q, \; q \in \bb Q$. Taking $q = (2g-2)^{-1}$ we obtain a factor\footnote{Although we can take arbitrary rational powers of the factor $\rho$, only powers $q \in \bb Z[(2g-2)^{-1}]$ will yield a factor of automorphy. This follows from applying the cocycle condition to the relation among the generators in the fundamental group of a compact surface. A more geometric argument would be, that other powers cannot be a factor, since the corresponding line bundle would not have integer degree.} corresponding to a bundle whose $(2g-2)$-th power equals the canonical bundle, hence it has degree $1$, and so $\rho^{qd}_\gamma(z)$ is a desired s-factor for any $d \in \bb Z$. $\Box$
%%%%%%%%%%%%%%%%%%%%%%%%%%%%%%%%%%%%%%%%%%%%%%%%%%%%%%%%%%%%%%%%%%%%%%%%%%%%%%%%%%%%55
\subsubsection{Examples: Half-canonical Bundles}
In this section we want to work out the factors of automorphy for half-integer tensor powers of the canonical bundle as an important class of examples. As always, we uniformize the Riemann surface $\Sigma$ by $\Sigma = D/G$.\\
 
A section of the canonical bundle\footnote{The canonical bundle of a complex manifols is the top exterior power of the holomorphic cotangent bundle; in complex dimension one this is the holomorphic cotangent bundle itself.} can locally be written $f(z)dz$. Setting $z':= g(z)$, we see that $f(g(z)) = g'(z)^{-1}f(z)$, and hence the corresponding factor of automorphy is given by $\rho^{\cal K}_{g}(z) := g'(z)^{-1}$. Observe $\rho^{\cal K}$ satisfies the cocycle condition \eqref{cocyc} because of the chain rule of differentation. This can be written in terms of one of the corresponding matrices in SL$(2, \bb C),$
$$  \rho^{\cal K}_{g}(z) = (cz + d)^2 \;,$$
and because of the square it is independent of the choice of the lifting of $\mathrm{PSL}(2,\bb C)$ to SL$(2, \bb C)$. \\

Because tensor products of line bundles correspond to the pointwise products of the factors of automorphy, the factor associated to $\cal K^n$ is given by 
$$\rho^{\cal K^n}_g(z) = g'(z)^{-n} = (cz + d)^{2n}\;.$$
All these factors are s-factors. We will take this example yet one step further: Recall that on a Riemann surface, any spin bundle $\cal S$ satisfies $\cal S^2 \cong \cal K$. The factors of automorphy are hence square roots of the canonical factor, which of course exist because $\rho^{\cal K}_g$ never vanishes and $D$ is simply-connected. The square roots are not unique. We have a $\bb Z_2$-ambiguity for the factor of each generator of the group and any two choices are non-equivalent\footnote{This follows from a general fact given below: For each factor of automorphy, there exists a unique s-factor in the same equivalence class. Any two choices are different s-factors and hence cannot be equivalent.}. Hence there exist $2^{2g}$ spin structures on a Riemann surface. Bundles of the type $\cal S^m \cong \cal K^{m/2}$ are called \emph{half-canonical}.\\

When working with spin structure there is one point where care is required: The square-root of the canonical factor, formally given by
$$\rho^{\cal S}_g(z) =\sqrt{ g'(z)} = (cz + d)\;,$$
depends on the choice of matrix representing $g$ in $\mathrm{SL}(2,\bb C)$. This dependence reflects the freedom of choice of square root: while there is now a $\bb Z_2$-ambiguity of the matrix representation of $g$, there seems to be a natural choice for the sqare root. Altogether there remains only a $\bb Z_2$-ambiguity, because if we choose the other matrix and the other square root the two signs that arise cancel. \\

The main point, however, is the following: Let $g_i$ be the generators of $G$, $M_i$ a choice of a matrix representing $g_i$ and $\tilde G$ the subgroup of $\mathrm{SL}(2,\bb C)$ generated by the $M_i$. Since the generators $g_i$ are subject to the relation
$$\prod_{i=1}^g g_i g_{i+g} g_i^{-1} g_{i+g}^{-1} = \bb 1\;,$$ in $\mathrm{PSL}(2, \bb C)$, we know that
that the product above with all $g_i$ substituted by $M_i$ is \emph{either} $\bb 1$ \emph{or} $-\bb 1$ in $\mathrm{SL}(2, \bb C)$.\\

In the latter case, the natural map $\tilde G \rightarrow G$ would be surjective with kernel $\bb Z_2$, while $\tilde G\cong G$ if the product is $\bb 1$. If we had $\tilde G \not \cong G$, the consequences would be drastic, e.g., there would not exist any non trivial automorphic forms for any non-integer half-canonical factors, because
$$ f(z) = f \left(\prod_{i=1}^g g_i g_{i+g} g_i^{-1} g_{i+g}^{-1} (z)\right) = \rho^{\cal K^{n/2}}_{\prod_{i=1}^g g_i g_{i+g} g_i^{-1} g_{i+g}^{-1}}(z) f(z) = (-1)^n f(z)\;,$$
and this of course implies $f\equiv 0$ for $n$ odd. Now that we have pointed this out, let us assure ourselves by  showing that this does not occur in the cases we are considering. We cite a result obtained by Kra that is more general than the one we need for our considerations in this section; the generality will come in handy later, however, when we consider noded Riemann surfaces and regular b-groups. 
\begin{thm} [\cite{OLKGtSL}]
    Let $G$ be a function group. Then $\tilde G \cong G$ iff $G$ does not contain an elliptic element of order two.
\end{thm}
%%%%%%%%%%%%%%%%%%%%%%%%%%%%%%%%%%%%%%%%%%%%%%%%%%%%%%%%%%%%%%%%%%%%%%%%%%%%%%%%%%%%%%%%%%%%%
\subsection{Banach spaces of Automorphic Forms}\label{BSoAF}
At several points later we will need Banach spaces of automorphic forms. In this section we briefly give the definitions and state the theorems we need later on. All material can be found, e.g., in \cite{AFaKG}.\\

From now on, $D$ will denote any hyperbolic open set of $\cinf$, i.e., an open set with more than one point in its complement, and $\lambda_D$ will denote its Poincar\'e density.\\

$L^p(D)$ will denote the space of $p$-integrable
measurable functions on $D$, and $\| \cdot \|_p$ will denote the $p$-norm. For an $s$ factor $\rho_s$ of a Kleinian group $G$ acting properly discontinuously on $D$, $L_{\rho_s}(D,G)$ will
denote the space of measurable automorphic forms, i.e., functions that satisfy $(f\circ g)(z) = f(z) \rho_s(g,z)$. On $L_{\rho_s}(D,G)$ the following expressions are well-defined norms,
\begin{eqnarray*}
   \|f\|_{L^p_s(D,G)} &:=& \left(\int_{\cal F} \lambda_{D}^{2-ps}(z)|f(z)|^p
d^2z\right)^{1/p}, \qquad 1 \leq p < \infty\;, \\
   \|f\|_{L^\infty_s(D,G)} &:=& \essup_{z \in \cal F} \left \{
\lambda_{D}^{-s}(z) |f(z)| \right \}\;.
\end{eqnarray*}
and the set of functions for which the norm is finite is denoted by $L^p_{\rho_s}(D,G)$. One can easily check that these are Banach spaces. Note that the integrals are not performed over $D$ but only over a
fundamental domain $\cal F$ for the action of $G$ on $D$. The subspace of holomorphic automorphic forms is denoted by \begin{eqnarray*}
   A^p_{\rho_s}(D,G)&:=& \cal O(D) \cap L^p_{\rho_s}(D,G)\;.
\end{eqnarray*}
In case $\infty \in D$ we will have to add the extra assumption\footnote{Otherwise the constructions will not be independent of the domain, e.g., the operation of pull-back introduced below would map functions holomorphic around $\infty$ to \emph{meromorphic functions} at the origin (see \cite{AFaKG} for details).} that $f \in O(|z|^{-2s})$. Also, since $L^p$-convergence of holomorphic functions implies local uniform convergence, these subspaces are closed and hence Banach spaces themselves.\\

We follow the tradition of denoting the spaces $A^\infty_{\rho_s}(D,G)$ by $B_{\rho_s}(D,G)$ and the spaces $A^1_{\rho_s}(D,G)$ by $A_{\rho_s}(D,G)$. Also, if we are dealing with an integer power $n$ of the canonical factor of automorphy, i.e., the factor $\rho_n(g,z) = g'(z)^{-n}$, we will use the subscript $n$ instead of $\rho_n$. And finally, when $G = \{\bb 1\}$ is the trivial group, we
simplify the notation and write $L^p_{s}(D)$ instead of $L^p_{\rho_s}(D,\bb 1)$, and similarily for the holomorphic subspaces. \\

These complex Banach spaces are of course independent of the chosen uniformization or, so to speak, invariant under
conjugation. More precisely, let $D$ be a simply-connected hyperbolic domain and $\phi:D \rightarrow \phi(D)$ be a biholomorphism. Then $\phi$ induces norm preserving isomorphisms for $1\leq p \leq \infty$, called the pull-back,
$$ \phi_s^*: L^p_{\rho_s}(\phi(D),\phi G \phi^{-1})\rightarrow
L^p_{\rho_s}(D,G)\;,\qquad (\phi_s^*f) (z) = (f\circ \phi) (z) \phi' (z)^s\;, $$ which respect the subspaces
$A^p_{\rho_s}$ of holomorphic functions.\\ \\
For conjugate numbers, i.e., $1/p + 1/p'=1$, we can introduce a product, called the \emph{Weil-Petersson pairing},
$$ L^p_{\rho_s}(D,G) \times L^{p'}_{\rho_s}(D,G)
\rightarrow \bb C\;, \qquad \langle f,g \rangle^G_s := \int_{\cal F} f(z) \ol {g(z)} \lambda^{2-2s}_D(z) d^2z \;.$$
\sloppy The integral is seen to converge by rewriting the integrand as follows,
$$f(z) \ol{g(z)} \lambda^{2-2s}_D(z) = \left(f(z) \lambda^{\frac{-sp}{p}}_D(z)\right) \cdot \left(\ol{g(z)} \lambda^{\frac{-sp'}{p'}}_D(z)\right) \cdot \lambda_D^2(z)\;,$$
and then applying H\"older's inequality with respect to the measure $\lambda^2_D(z) d^2z$. This also establishes the fact that the Weil-Petersson pairing induces an isometric isomorphism 
\begin{eqnarray} \label{IsoIso}L_{\rho_s}^{p'} (D,G) \cong \left(L_{\rho_s}^{p} (D,G)\right)^*\end{eqnarray}
for any $1\leq p<\infty$ and any group $G$ acting on $D$. It is remarkable that the subspaces of holomorphic automorphic forms respect this duality, albeit not isometrically. 
\begin{thm}[\cite{AFaKG}] \label{NondegenerateOnAs}
  For $1\leq p<\infty$ the anti-linear map
  $$  A_{\rho_s}^{p'} (D,G) \rightarrow \left(A_{\rho_s}^{p} (D,G)\right)^*\;, \qquad \psi \mapsto l_\psi:=\langle \cdot, \psi\rangle^G_s\;,$$
  is an isomorphism which satisfies the norm inequality 
  $$ c_s^{-1} \| \psi \|_{L^{p'}_s(D,G)} \leq \| l_\psi \| \leq \|\psi\|_{L^{p'}_s(D,G)}\;.$$
\end{thm}
Let us observe some relations between the different spaces.
\begin{lemma} \label{lincl}
    For $1<p<\infty$ and $G$ a function group with invariant component $D$, we have the following inclusions:
    \begin{eqnarray*}  &&L^\infty_{\rho_s}(D,G) \subset
L^p_{\rho_s}(D,G)\subset L^1_{\rho_s}(D,G), \qquad
A^\infty_{\rho_s}(D,G) \subset
A^p_{\rho_s}(D,G)\subset A^1_{\rho_s}(D,G)\;.
\end{eqnarray*}
\end{lemma}
\emph{Proof.} By definition, $D/G$ has finite area since $G$ is a function group. Hence the hyperbolic area of any fundamental domain $\cal F$ is finite, and we compute for $f \in L^\infty_{\rho_s}(D,G)$,
$$ \|f\|_{L^p_s(D,G)} =  \int_{\cal F} \lambda_D^{2-s}(z)|f(z)| d^2z \leq \|f\|_{L^\infty_s(D,G)} \int_{\cal F} \lambda_D^2(z) d^2z< \infty\;.$$ For the second inclusion we apply
H\"older's inequality with respect to the Poincar\'e measure,
\begin{eqnarray*}
\|g\|_{L^1_s(D,G)} &=& \int_{\cal F}|g(z)| \lambda^{-s}_D(z)\lambda^2_D(z)d^2z\\&\leq& \left(\int_{\cal F}
\lambda^{-sp}_D(z)|g(z)|^p\lambda^2_D(z)d^2z \right)^{\frac1p}\cdot\left( \int_{\cal F} \lambda^2_D(z) d^2z\right)^{\frac1{p'}} \leq c \|g\|_{L^p_s(D,G)}\;.
\end{eqnarray*}
The assertions for the subspaces of holomorphic functions clearly
follow. $\Box$ \\ 

There are also reverse inclusions, namely for any finitely generated Kleinian group $G$ with invariant simply connected
component $D$ (i.e., $G$ is a b-group) as we see in the next theorem. Observe that the isomorphisms here are only claimed for the subspaces of holomorphic automorphic forms. They are not true for all measurable automorphic forms.
\begin{thm}\label{incl}
     Let $G$ be a b-group with simply-connected invariant component $D$.
     Then the spaces
     $$ A^p_{\rho_s}(D,G) \cong A^q_{\rho_s}(D,G)\;,\qquad 1\leq p \leq q \leq \infty\;.$$
     are continuously isomorphic. In other words, the norms $\|\cdot\|_{A^p_s(D,G)}$ are all equivalent for $1\leq p \leq \infty$.
\end{thm}
In the case $p=1$ and $q=\infty$ these reverse inclusions, together with their norms, are crucial for understanding and bounding the degeneration of sections as the boundary of Teichm\"uller space is approached, which we will sketch in Section \ref{AsofSections}. Therefore we will look at this particular case more carefully, but will postpone this for a couple of pages until we have introduced the \emph{Bergman projection} and the \emph{Poincar\'e series} operator. \\

The automorphic Bergman kernel $K_{D,s}$ is the kernel of the Bergman projection operator, an operator mapping measurable automorphic forms to holomorphic ones. To describe it, let
\begin{eqnarray}\label{kernelinv}
 k_{\bb D}: \bb D \times \bb D \rightarrow \bb C\;, \qquad k_{\bb D}(z,w) := \frac{1}{\pi(1-z\bar w)^2}\;.
\end{eqnarray}
be the well-known Bergman kernel on the unit disc. For any other simply-connected domain $D$, we choose a Riemann mapping, i.e., a biholomorphism $\psi: \bb D \rightarrow D$, and define the kernel as follows,
$$ k_{D} (\psi(z),\psi(w)) \psi'(z)\ol{\psi'(w)} = k_{\bb D} (z,w)\;.$$
In case $D$ is multiply connected or has several components, each is hyperbolic and therefore covered by the unit disc. We can hence define the kernel for each component in the same way as above by choosing a holomorphic covering map. Observe that $k_D$ is well-defined because $k_{\bb D}$ behaves invariantly under $\mathrm{Aut}(\bb D)$. Finally, we define a slightly modified kernel,
\begin{eqnarray*} 
	\label{BergmanDef} K_{D,s}(z,w):= (2s-1)\pi^{s-1}\left(k_D(z,w)\right)^s\;, \qquad c_s:=\frac{2s-1}{s-1}\;,
\end{eqnarray*}
which we call the $s$-Bergman kernel. Before we come to the theorem on the projection operator, let us prove an important lemma we will use several times is the following.
\begin{lemma} \label{prep}
		For any $s>1$, we have the following formula,
   $$ \int_D \lambda^{2-s}_D(z) |K_{D,s}(z,w)| d^2 z = c_s \lambda^s_D(w)\;.$$
\end{lemma}
\emph{Proof.} We can choose a special domain to check the formula,
since everything is invariantly defined. Of course we choose the
unit disc. Let's also check it for $w=0$ first. Then,
\begin{eqnarray*}
   \int_{\bb D} \lambda^{2-s}(z) |K_{\bb D,s}(z,0)| d^2 z &=& \frac{2s-1}{\pi} \int_{\bb D} (1-|z|^2)^{s-2}d^2 z \\ &=& 2(2s-1) \int_{0}^1 (1-r^2)^{s-2}r dr \;=\; \frac{2s-1}{s-1} \;= \; c_s \lambda(0)^s\;.
\end{eqnarray*}
Now choose an automorphism of the disc mapping $w \mapsto 0$.
Since both sides have the same transformation behaviour, the
formula is valid for all $w \in \bb D$. $\Box$
\begin{thm} \label{last}
    The operator defined by the expression
    $$ f\mapsto (\beta_{\rho_s}f)(z):= \int_{D} \lambda^{2-2s}_D(w)K_s(z,w)f(w) d^2w\;,$$
    is a well-defined projection operator $L^p_{\rho_s}(D,G)\rightarrow A^p_{\rho_s}(D,G)$ of norm at most $c_s$.           Moreover it is symmetric with respect to the Weil-Petersson pairing, i.e.,
     $$\langle \beta_{\rho_s} f,g\rangle_s^G = \langle f,\beta_{\rho_s} g\rangle^G_s\;.$$
\end{thm}
The final theorem we need later on concerns the normal convergence of the following series called the \emph{Poincar\'e series},
$$\Theta_{\rho_s} [f] (z):= \sum_{g \in G} f(g z) \rho_{g}(z)^{-1}\;,$$
for given $s$-factor of automorphy $\rho_s$. This operator produces an automorphic form by averaging the value of the function along the orbit in an automorphic way.
\begin{thm}\label{PS}
   If $\|f\|_{L^1_s(D)} < \infty$, the series $\Theta_{\rho_s}[f]$
converges normally. Moreover $\Theta_{\rho_s} [f] \in
L^1_{\rho_s}(D,G)$ and $$ \| \Theta_{\rho_s} [f] \|_{L^1_{\rho_s}(D,G)} \leq \|f\|_{L^1_s(D)}\;.$$
\end{thm}
This means that the Poincar\'e operator $\Theta_{\rho_s}$ is a bounded linear operator $L^1_s(D)\rightarrow L^1_{\rho_s}(D,G)$. By normal convergence its restriction to $A^1_s(D)$ maps into $A^1_{\rho_s}(D,G)$. The Poincar\'e operator is compatible with the Weil-Petersson product in the following way, which can be checked by straight-forward calculation.
\begin{lemma} \label{scalar}
   Let $f \in L_{\rho_s}^\infty(D,G)$, $g \in L^1_{\rho_s}(D,G)$ and $g=\Theta_{\rho_s}[h]$ for some function
$h \in L_s^1(D)$. Then the scalar product can be computed by
$$ \langle f,g\rangle^G_s = \int_D f(w) \ol{h(w)} \lambda_D^{2-2s}(w) d^2w = \langle f, h\rangle_s^{\bb 1}\;.$$
\end{lemma}
Later on, we will need a very remarkable theorem by Marvin Knopp (\cite{ACTfAFaRT}, Thm.~ 3), which says that we can approximate automorphic forms arbitrarily well by ones that have a polynomial preimage under $\Theta$. 
\begin{thm}\label{Knopp}
   Let $G$ be a b-group with bounded simply-connected invariant component $D$, and let $s\geq 2$. Then
   the image of $\bb C[z]$ under $\Theta$ is dense in $A^1_{\rho_s}(D,G)$.
\end{thm}
Indeed, in the next section we will see that for $s\geq 2$ and $D$ as above, $\bb C[z]_{|D} \subset A^1_s(D)$ so the statement makes sense.\\

Let us now return to the isomorphism stated in Thm. \ref{incl} for the case $p=1, q=\infty$. This question of
whether such an isomorphism exists was an actively discussed research problem and was solved in a series of papers
(\cite{OtBoAF}, \cite{OIaBAF}, \cite{OIaBAF2}, \cite{BaIAF}) in increasing degrees of generality. \\

First of all, in case $G$ is a b-group with simply-connected invariant component $D$, we can restrict our attention to $D= \bb D$ due to everything being well-behaved under pull-back by a Riemann mapping. For the disc, however, Metzger
and Rao \cite{OIaBAF} proved a statement equivalent to the existence of an inclusion. It is stated in terms of a function $\alpha_s(z,w)$ which is exactly the Poincar\'e series with respect to $z$ of the modified Bergman kernel,
$$ \alpha_s(z,w) := \Theta[K_{\bb D,s}(\cdot, w)](z) = \sum_{g \in G} \rho^{-1}_g(z)K_{\bb D,s}(gz, w)\;.$$
This is well defined, since for fixed $w \in \bb D$ the function $K_{\bb D,s}(z, w)$ is contained in $A_s(\bb D)$, which can easily be checked using the \ref{prep}.
\begin{thm}
Let G be a Fuchsian group acting on $\bb D$. Then
$$A_s(\bb D,G) \subset B_s(\bb D,G) \quad \Longleftrightarrow \quad \sup_{z \in \bb D} \lambda_{\bb D}^{-2s}(z) \alpha_s(z,z) <
\infty\;.$$
\end{thm}
\emph{Proof}. Let $f \in A_s(\bb D,G)$. The equality $\beta_s f = f$ written explicitly using Lemma \ref{scalar}
yields
\begin{eqnarray*}
 f(u)&=&\int_{\bb D}K_{\bb D,s}(u,w)f(w)\lambda_\bb D^{2-2s}(w) d^2w\\
       &=&\int_{\cal F}\alpha_s(u,w)f(w)\lambda_\bb D^{2-2s}(w) d^2w\;.
\end{eqnarray*}
Now we take the modulus of this, multiply through and estimate the integral,
\begin{eqnarray}\label{estimate}
 |f(u)|&\leq& \sup_{w \in \cal F}\left(\lambda_{\bb D}^{-s}(w)|\alpha_s(u,w)|\right) \|f\|_{A_s(\bb D,G)}\;.
\end{eqnarray}
The equality $\beta_s f = f$ is true in particular for $f= \alpha_s(\cdot, w)$, in which case it yields
$$ \alpha_s(z,w) = \int_{\cal F} \alpha_s(u,w) \alpha_s(z,u)\lambda_{\bb D}^{2-2s}d^2 u\;,$$
and as a special case
$$ \alpha_s(w,w)=\int_{\cal F} |\alpha_s(w,u)|^2\lambda_{\bb D}^{2-2s}d^2w\;.$$
We can now apply the Cauchy-Schwarz inequality to the integral representation of $\alpha_s(w,w)$ and obtain the
estimate
$$ |\alpha_s(z,w)| \leq \sqrt{\alpha_s(z,z) \alpha_s(w,w)}\;.$$
If we multiply equation (\ref{estimate}) through with $\lambda^{-s}_{\bb D}$ and take the supremum over $u \in \bb D$,
we finally get
\begin{eqnarray*}
    \|f\|_{B_s(\bb D,G)} &\leq& \sup_{z,u \in \bb D} \left(\lambda_{\bb D}^{-s}(w)\lambda_{\bb D}^{-s}(z)\sqrt{\alpha_s(w,w)\alpha_s(z,z)}\right) \|f\|_{A_s(\bb D,G)}\\
                   &\leq& \sup_{z \in \bb D} \left(\lambda_{\bb D}^{-2s}(z)\alpha_s(z,z)\right)\|f\|_{A_s(\bb D,G)}
\end{eqnarray*}
To prove the converse, suppose there is a constant $M$ such that $\|f\|_{B_s(\bb D,G)} \leq M \|f\|_{A_s(\bb D,G)}$ for all $f \in A(\bb
D)$. Then, in particular,
$$ \alpha_s(w,w) \lambda_{\bb D}^{-s}(w) \leq M \| \alpha_s(\cdot, w)\|_{A_s(\bb D,G)}\;,$$
which can be estimated because of Lemma \ref{prep}:
\begin{eqnarray*}
   \|\alpha_s(\cdot,w)\|_{A_s(\bb D,G)} &\leq& \int_{\bb D} |K_{\bb D,s}(u,w)|\lambda_{\bb D}^{2-s}(u) d^2u
        \leq c_s \lambda_{\bb D}^s(w)\;.
\end{eqnarray*}
Hence altogether we arrive at the desired result
$$ \alpha_s(w,w) \lambda_{\bb D}^{-2s}(w) \leq M c_s\;. \;\Box$$\\
Metzger and Rao further proved that finitely generated Fuchsian groups without parabolic elements satisfied the
preceding condition, and in the succeding paper \cite{OIaBAF2} they established this for any finitely generated
Fuchsian group. Observe that the bound so obtained for $M$,
\begin{eqnarray} \label{EmbeddingConstant} M \leq \sup_{z \in \bb D} \left(\lambda_{\bb D}^{-2s}(z)\alpha_s(z,z)\right)\;,\end{eqnarray}
depends on $s$ \emph{and} on the group $G$, the latter dependence being hidden in the $\Theta$-operator, where $G$
enters via $|\rho_g(z)| = |g'(z)|^s$.\\
%%%%%%%%%%%%%%%%%%%%%%%%%%%%%%%%%%%%%%%%%%%%%%%%%%%%%%%%%%%%%%%%%%%%%%%%%%%%%%%%%%%%%%%55
%%%%%%%%%%%%%%%%%%%%%%%%%%%%%%%%%%%%%%%%%%%%%%%%%%%%%%%%%%%%%%%%%%%%%%%%%%%%%%%%%%%%%%5555555555
\subsection{Integrability of the Poincar\'e Density}\label{IotPD}
We want to make a short comment concerning the integrability of powers of the Poincar\'e density because this issue will come up in our construction of sections later. Particularily, if one wants to conclude that bounded functions have finite $L^p_s$-norm or if one wants to apply a H\"older inequality to integrands of the type $|f|^p\lambda^{2-ps}$, and other similar situations.\\

The following lemma follows from the fact that the Poincar\'e density is bounded away from zero for bounded domains.
\begin{lemma}\label{TwoMinusS}
   For a bounded domain $D$ and $s\geq 2$, the Poincar\'e density is $(2-s)$-integrable, i.e.,
   $$ \int_D \lambda_D^{2-s}(z) d^2 z < \infty \;.$$
\end{lemma}
For an unbounded domain this is of course not true, e.g., for $s=2$ the integral is equal to the euclidean area of $D$, which can be infinite. In fact, for the integral to have a finite contribution from infinity, the integrand would have to be of order $O(|z|^{2+\epsilon})$. But for an unbounded domain, $\lambda_D^{-1} \in \cal O(|z|^{-1})$, which would imply $2-s>2$, i.e., $s<0$. But then again we also have $\lambda_D \in \cal O(\mathrm{dist}(z,\partial D))$, hence terms with $s<0$ are not integrable when we approach a finite boundary point. Thus
\begin{lemma}
   Let $D$ be an unbounded hyp domain. Then
   $$  \int_D \lambda_D^{2-s}(z) d^2 z = \infty \qquad \forall\: s\in \bb R\;.$$
\end{lemma}
Let's look at the bounded case a little more closely. If we define
$$ s_0(D):= \mathrm{inf}\{s \mid  \|\lambda^{2-s}_D\|_1 < \infty\}\;,$$
then $s_0 \in [1,2]$ depends crucially and maybe somewhat surprisingly on the shape of the domain and in particular on the regularity of the boundary $\partial D$. For example, $s_0(D) = 1$ for domains such that $\partial D$ is rectifiable. For more on these issues, see \cite{OPDiAqD}, \cite{PAitBS} and the references given therein. For the unit disc though, $s_0 = 1$.
\begin{lemma} \label{PoincUnitDisc} For the unit disc we have
    $$ \int_{\bb D} \lambda_{\bb D}^{2-s} (z) d^2z < \infty
    \qquad
    \Leftrightarrow \qquad s>1\;.$$
\end{lemma}
\emph{Proof}. We plug in the Poincar\'e density in polar coordinates and calculate,
$$ \int_{\bb D} \lambda_{\bb D}^{2-s}d^2 z = C \int_0^1 (1-r^2)^{s-2}r dr\;.$$
The last expression can easily be seen to be finite iff $s>1$. $\Box$ \\

A general statement concerning integrability is given by the following lemma, which is valid for
arbitrary simply-connected domains, and, of course, by successive application also for domains with multiple components. However, it is of little practical use since it depends on a Riemann mapping, which is very hard to come by.
\begin{lemma}
    Let $D\subset \cinf$ be a simply connected domain, $\psi:D
    \rightarrow \bb D$ a Riemann mapping such that $\psi^{-1}(0) \neq \infty$ and $s>1$. Then the derivative of
the Riemann mapping belongs to $A^1_{s}(D)$, i.e.,
    $$\int_D \lambda_D^{2-s}(z) |\psi'(z)| d^2 z <
    \infty\;.$$
\end{lemma}
\emph{Proof}. By Lemma \ref{prep}, we have
$$ \int_D \lambda_D^{2-s}(z) |K_{D,s}(z,w)| d^2z = c_s \lambda_D^s(w)\;,$$
which is of course finite for all $\infty \neq w \in D$. Now choose $w:=\psi^{-1}(0)$,
\begin{eqnarray*}
K_{D,s}(z,w) &=& K_{\psi(D),s}(\psi(z),\psi(w)) \psi'(z) \overline
{\psi'(w)} = C \: \left(k_{\bb D}(\psi(z),\psi(w))\right)^s\psi'(z)
\overline
{\psi'(w)}\\
&=& C \: \frac{\psi'(z) \overline {\psi'(w)}}{\pi^s(1 - \psi(z)
\overline{\psi(w)})^{2s}}
 = \tilde C \: \psi'(z)\;. \: \Box
\end{eqnarray*}
%%%%%%%%%%%%%%%%%%%%%%%%%%%%%%%%%%%%%%%%%%%%%%%%%%%%%%%%%%%%%%%%%%%%%%%%%%%%%%%%555
%%%%%%%%%%%%%%%%%%%%%%%%%%%%%%%%%%%%%%%%%%%%%%%%%%%%%%%%%%%%%%%%%%%%%%%%%%%%%%%%%%5
\subsection{The Relation between Bounded Automorphic Forms and Sections of
Line Bundles} \label{RSaAF} 
As we know from Section \ref{LBvsFoA}, automorphic forms for a factor $\rho$ can be identified with sections of a line bundle which corresponds to this $\rho$. The new ingredients we have at hand now are the Banach norms on the spaces of automorphic forms, and in this section we will discuss the \emph{boundedness} of the norm as a geometric feature.\\

First of all, we noticed in Theorem \ref{incl} that all norms $\|\cdot\|_{L^p_s(D,G)}$, $p \in [1, \infty]$, are equivalent if $G$ is a b-group with invariant component $D$, so it makes no difference which one we consider when uniformizing a line bundle over a surface of finite analytic type, which is what we generally want to do.\\

Now if $\Sigma = D/G$ is a compact Riemann surface, any automorphic form will be bounded, since any fundamental domain $\cal F\subset D$ has compact closure in $D$ and the $p=1$ norm is obtained by integration over $\cal F$. The boundedness of the norm can hence in general only induce some restriction for the automorphic form at a boundary point of $D$, which corresponds to a puncture of $D/G$. It is therefore clear that
$$ H^0(\Sigma, \cal L(\rho)) \cong A^p_{\rho}(D,G)\;,$$
for any $p\in [1, \infty]$. \\

The corresponding statement for a surface with punctures is a little more complicated. First of all, recall that a holomorphic line bundle on a non-compact Riemann surface is trivial, so the symbol $\cal L(\rho)$ has to be interpreted correctly. We do this with the help of the filled-in surface $\tilde \Sigma$ of $\Sigma$. By this we mean the surface obtained from $\Sigma$ by filling in the punctures, which is a canonical procedure (for more details, see \cite{ACaRS}, p.66). We denote by $N$ the injective holomorphic map $N: \Sigma \hookrightarrow \widetilde \Sigma$.\\

Take an automorphic form $\phi \in A^p_\rho(D,G)$ and consider its zero divisor, i.e., the formal sum of points $p_i$ where $\phi$ vanishes. This descends to the surface $\Sigma$, and via the map $N$ we can push forward this divisor to a divisor $D(\rho) \in \mathrm{Div}(\tilde \Sigma)$. The puncture divisor $P \in \mathrm{Div}(\tilde \Sigma)$ is defined as the divisor given by the sum over all the points $p \in \widetilde \Sigma$ which don't have a preimage under $N$, with coefficient one. And finally, for a real number $s$, let $[s]$ denote the biggest integer which is strictly smaller than $s$\footnote{So if $s$ is an integer, $[s] = s-1$ whereas if $s$ is any positive non-integer, $[s]$ is the number in front of the decimal point in the decimal representation of $s$.}. The following statement should be clear by the discussions above and the estimates in Section \ref{IotPD}.
\begin{prop}
   Let $G$ be a b-group with invariant simply connected component $D$. Let $\tilde \Sigma$ denote the compact surface obtained by filling the punctures of the surface $\Sigma=D/G$, $P \in \mathrm{Div}(\tilde \Sigma)$ the puncture divisor and $\rho_s$ an s-factor for the action of $G$ with $s>1$. Then
   $$ H^0(\widetilde \Sigma, \cal L(\rho_s)) \cong A_{\rho_s}^p(D,G)\qquad \mathrm{where} \qquad \cal L(\rho_s) = \cal L(D(\rho) + [s]P)\;.$$
\end{prop}
From this point of view, it is now easy to understand why we had to restrict to the case $s>1$. s-factors with $s>1$ correspond precisely to the topological range, i.e., to line bundles of degree $\geq 2g-1$. This follows from the canonical factor being a 1-factor and $\mathrm{deg}(\cal K) =2g-2$. But it is nevertheless very interesting to observe that from an analytic point of view the possible non-existence of sections of line bundles of degree $\leq 2g-2$ can be explained by the fact that the Poincar\'e series operator is not well-defined on the space $A_s^1(D)$ for $s<1$. \\

The Riemann-Roch formula \eqref{R-R} now translates into a dimension formula for the spaces $A^p_{\rho_s}(D,G)$.
\out{There is, however, one subtle point in the formula: If punctures are present, we explained in Section \ref{LBvsFoA} that for any $s \in \bb R$ there exists an $s$-factor for the group,} 
\begin{lemma}\label{DimensionFormula}
    Let $G$ be a b-group with simply-connected invariant component $D$, $\rho_s$ an s-factor for the action of $G$ on $\bb D$ with $s>1$, and let $n$ denote the number of punctures on the surface $\Sigma= D/G$. Then the following dimension formula holds:
    $$\mathrm{dim} A^p_{\rho_s}(D,G) = (2s - 1)(g-1) + n [s]\;.$$
\end{lemma}
In particular, $A^p_{\rho_s}(D,G)$ is finite dimensional iff $D/G$ is of finite type.
%%%%%%%%%%%%%%%%%%%%%%%%%%%%%%%%%%%%%%%%%%%%%%%%%%%%%%%%%%%%%%%%%5
%%%%%%%%%%%%%%%%%%%%%%%%%%%%%%%%%%%%%%%%%%%%%%%%%%%%%%%%%%%%%%%%%5
%%%%%%%%%%%%%%%%%%%%%%%%%%%%%%%%%%%%%%%%%%%%%%%%%%%%%%%%%%%%%%%%%5
\addtocontents{toc}{\vspace{.7cm}}
\section*{\bf The Relative} 
\setcounter{section}{2}
\setcounter{subsection}{0}
\addtocontents{toc}{\vspace{.3cm}}
In this section we start the study of families of objects rather than the objects themselves. Namely, in order to undersand the 'collection of all objects', i.e., their 'moduli space' so to speak, it is inevitable to know all possible families of objects in the category. This was Grothendiecks revolutionary insight, which in the end culminates at the notion of an 'universal object'. But we will not digress further into vague explanations and rather give a brief outline of this section.\\

The first subsection contains the general definitions in the context of families, and is followed by Section \ref{Teich}, which is devoted to introduce a model for the base space of the universal family of marked Riemann surfaces, namely Teichm\"uller space. \\

The universal object and some related ones is presented in Section \ref{FibreoverTeich}, while finally, in Section \ref{faforms} we come to the main topic of this section, namely we show in great detail how families of automorphic forms are constructed.
%%%%%%%%%%%%%%%%%%%%%%%%%%%%%%%%%%%%%%%%%%%%%%%%%%%%%%%%%%%%%%%%%%%%%%%%%%%%%%%%%%%%%%%%%%%5555555
\subsection{Families of Complex Manifolds and Line Bundles}\label{GenFam}
In this section we briefly gather the general definitions of families we will consider later on. The first one is due to Kodaira and Spencer \cite{CMaDoCS}.
\begin{dfn}
    A \underline{family of compact complex manifolds} consists of a
pair of connected complex manifolds $X$ and $B$ together with a
proper surjective holomorphic submersion $\pi: X \rightarrow B$.
Such a family is called a \underline{deformation of $\Sigma$} iff
$\Sigma \cong  \pi^{-1}(t)$ for some $t \in B$.
\end{dfn}
Let us remark a few things about such families.
\begin{itemize}
	\item{The data of a family will be abbreviated by either
$(X,B,\pi)$ or simply $\pi:X\rightarrow B$.} 
	\item{Since $\pi$ is proper, the preimages of points are compact. Moreover, $d\pi$ has
constant rank, hence the fibres are compact complex submanifolds of $X$,
which will be denoted by $X_t := \pi^{-1}(t)$, $t \in B$.}
    \item{By a theorem of Ehresmann, the differentiable structure of
the fibres cannot change, i.e., $X_{t_1}\cong_{C^\infty}X_{t_2}$
for any $t_1, t_2 \in B$, and hence neither the topology.}
   \end{itemize}
The notion of isomorphism in this category is the following.
\begin{dfn}
   Two families over the same base $(X_1, B, \pi)$ and $(X_2, B,
\psi)$ are \underline{equivalent} iff there is an equivariant
biholomorphism $\Phi: X_1\rightarrow X_2$, i.e. $\pi = \psi \circ
\Phi$.
\end{dfn}
Moreover, we have the fundamental construction of pull-back; in more fancy language, the fibered product is
representable in the category of complex manifolds over $B$. This means the following: Given a family $(X,B,\pi)$ and
a holomorphic map $f:A\rightarrow B$, we obtain a family $(f^*(X), A, \pi_1)$ by setting
$$ f^*(X) := A\times_f X = \{(a,x) \in A \times X, f(a) = \pi(x)
\}\;,\qquad \pi_1(a,x) = a \;. $$ One sees immediately that $\pi_1^{-1}(a) \cong X_{f(a)}$. \\

The definition of a family of complex manifolds equipped with line bundles is rather obvious.
\begin{dfn}
   A family of complex manifolds equipped with line bundles is a family of complex manifolds $(\pi:\cal X\rightarrow B)$ together with a line bundle $\cal L \rightarrow \cal X$ over the total space of the family.
\end{dfn}
Observe that in this way for each $p \in B$, $X_p:=\pi^{-1}(p)$ is a complex manifold equipped with the line bundle $\cal L_p := \cal L|_{X_p}$, and hence, the manifolds not only form a family, but the individual line bundles also line up to form a global bundle over the total space
%%%%%%%%%%%%%%%%%%%%%%%%%%%%%%%%%%%%%%%%%%%%%%%%%%%%%%%%%%%%%%%%%%%%%%%%%%%%%555
\subsection{Teichm\"uller Spaces and Beltrami Differentials}
\label{Teich}
%-------------------------------------------------------------------------------------
In this paragraph we present a brief construction of the base space of the universal object in the category of families of marked Riemann surfaces, known as Teichm\"uller space. We will only cover the material we need later on, and refer the interested reader to the textbooks on this interesting and broad subject (e.g., \cite{TCAToTS}, \cite{ItTS}, \cite{UFaTS}, \cite{TTvol1vol2} or the survey article \cite{FDTSaG}).\\

The set-theoretic definition of the Teichm\"uller space of a Riemann surface is the following.
\begin{dfn} \label{DefTeich}
    Let $\Sigma$ be a Riemann surface. The \ul{Teichm\"uller space} $\mathrm{Teich}(\Sigma)$ of $\Sigma$ is the set of
    equivalence classes of triples $(\Sigma, f, \Sigma')$, where $\Sigma'$ is another Riemann surface and $f:\Sigma \rightarrow \Sigma'$
    is a quasiconformal homeomorphism, under the following equivalence relation:
    $$f:\Sigma \rightarrow \Sigma' \sim g:\Sigma \rightarrow \Sigma''
    \quad \Longleftrightarrow \quad  \exists \textrm{ holomorphic
    }\sigma:\Sigma' \rightarrow \Sigma''\; \textrm{s.t.}
    \;\bar g^{-1}\circ \bar \sigma\circ \bar f  \simeq \bb 1_\Sigma \: \mathrm{rel} \: \partial \Sigma\;,$$ 
    where $\bar h$ denotes the extension of the quasiconformal $h$ to $\Sigma \cup \partial \Sigma$. The equivalence classes will be denoted $[\Sigma,f,\Sigma']$, and
    $\Sigma$, which will be identified with $[\Sigma, \bb 1_{\Sigma},
    \Sigma]$, is called the base point of $\mathrm{Teich}(\Sigma)$.
\end{dfn}
It turns out that Teichm\"uller space itself is a complex manifold in a natural way as we will see in the course of the next two sections sections. 
\subsubsection{The Beltrami Model of Teichm\"uller Space}
The complex-analytic theory of Teichm\"uller spaces is founded on the solution of the Beltrami equation on the plane as presented by Ahlfors and Bers\footnote{This theorem has a quite complicated history. The reason we quote Ahlfors and Bers paper here instead of earlier work is, that they furthermore proved a crucial reqularity statement for varying $\mu$'s. A particular case important for Teichm\"uller theory is the fact that if $\mu$ varies holomorphically in the Banach space, then $\mathrm{ev}_z: L^\infty(\bb C)_1 \rightarrow \cinf$ is a holomorphic map.} \cite{RMTfvM}.
\begin{thm} \label{BeltSol}
    Given any $\mu \in L^\infty(\bb C)$ of norm $\|\mu\|_\infty =
    k < 1$, there exists a $\frac{1+k}{1-k}$-quasiconformal homeomorphism
    $w[\mu]:\cinf \rightarrow \cinf$ solving the Beltrami equation 
    $$ \partial w[\mu]= \mu \; \bar \partial w[\mu]\;.$$
    Any other solution is obtained by post-composition with a suitable M\"obius transformation.
\end{thm}
Let us now describe how this theorem fits into the picture. To any $\mu \in L^\infty(\bb C)$ of norm less than one we can associate a quasiconformal homeomorphism $w[\mu]$ given by the theorem by choosing a normalization, i.e., a representative of the M\"obius orbit of solutions. If we are given a Kleinian group $G$, and an invariant component $D$ of $\Omega(G)$, we obtain a new domain $D^\mu := w[\mu](D)$ together with a group $G^{\mu} := w[\mu]G w[\mu]^{-1}$ which preserves $D^\mu$. All \emph{topological}
properties of the group action are preserved, since it is obtained by composition
with homeomorphisms. In particular, the domain of discontinuity of
$G^\mu$ is $w[\mu](\Omega(G))$, and $D\mu$ is an invariant component. By direct calculation one can check that $G^\mu$ acts as biholomorphisms of $D^\mu$, iff $\mu$ additionally satisfies 
\begin{eqnarray} \label{BeltFact} 
\mu(z) = \mu (g z) \frac{\overline{g'(z)}}{ g'(z)}\;, \qquad \mathrm{ for \;almost\; all} \: z \in \Omega(G), \; \forall \;\gamma \in
G\;.\end{eqnarray}
The space of all such $\mu$ for a given group $G$ will be denoted by $\mathrm{Belt}(D,G)$. Such $\mu$ are called, by a slight abuse of nomenclature,  \emph{Beltrami differentials}\footnote{Beltrami differentials should be thought of as being measurable sections sections of $T_{(1,0)}\Sigma \otimes T^*_{(0,1)}\Sigma$ with $L^\infty$-norm bounded by $1$. As for other sections of line bundles, we have a natural identification of these with $\mathrm{Belt}(D,G)$ so there should be no confusion.} for $G$. In particular, if \eqref{BeltFact} holds almost everywhere on $\cinf$, e.g., on all of $\Omega(G^\mu)$, then $G^\mu$ is again a Kleinian group.\\

If we started with a Fuchsian group of first kind, acting on $\bb D$, the domain $D^\mu$ we obtained is of course a quasidisc and $G^\mu$ is a quasi-Fuchsian group. Observe that the type of the Riemann surface $X^\mu:= D^\mu/G^\mu$ is the same as the type of $X:=\bb D/G$, since the type is determined by the nature
of the generators of $G$. This again is determined by the number of their fixed points and their location (i.e., in the regular or in the limit set); data which is obviously preserved by the deformation.\\

So starting from a Beltrami differential $\tilde \mu \in \Belt(\Sigma)$, we lift it to the covering $\Sigma = \bb D/G$. This function we call $\mu$. We have to extend $\mu$ to the whole complex plane in order to apply Thm. \ref{BeltSol}, and there are two natural choices, of which we choose to set $\mu (z) \equiv 0$ on $\bb D^c$. By construction\footnote{Due to $\mu$ restricted to $\bb D^c$ being in $\mathrm{Belt}(\bb D^c,G)$, we are actually constructing \emph{two} new Riemann surfaces, $D^\mu/G^\mu$ and $(D^\mu)^c/G^\mu$, a process called \emph{simultaneous uniformization}. Of course, $\mu$ can be extended in a way that it is \emph{not} in $\mathrm{Belt}(\bb D^c,G)$ when restricted to $\bb D^c$, but this possibility has, to the best of our knowledge, not led to any new insights or applications.}
, $\mu$ is now in 
$$\Belt(\Omega(G),G)  \cong \mathrm{Belt}(\bb D,G)\oplus \mathrm{Belt}(\bb D^c,G)\;.$$
We now pick the solution $w^\mu$ to the Beltrami equiation, which we normalize\footnote{We pick this normalization to have a uniform bound on the quasidiscs $D^\mu$: Viewed in the chart around $\infty$, $w^\mu|_{\bb D}$ is a normalized schlicht function, and hence due to Koebe's $\frac 14$ Theorem, $\bb D_{\frac 14} \subset w^\mu (\bb D) = (D^\mu)^c$. Hence, in the chart around the origin, this means the domains $D^\mu$ are all contained in $\bb D_4$.} by $w^\mu(\infty) = \infty$, and using a chart $\eta = z^{-1}$ at infinity, by $(w^\mu)'_{|\eta=0} = 1$ and $(w^\mu)''_{|\eta=0}  = 0$. Maps normalized in this way will be called \emph{1-point normalized} (at the origin)\footnote{In contrast to the normalization by prescribing the image of three points under the map, which we would call a 3-point normalization.}.\\

For $p \in \Sigma$, it's preimage in $\bb D$ is the orbit $Gp'$, where $p'$ is an arbitrary point in the preimage. $w^\mu$ maps this orbit again to an orbit: $$w^\mu(G p')=w^\mu (G(w^\mu)^{-1}w^\mu p') = G^{\mu} w^\mu p'\;.$$ Therefore, by setting $\Sigma^\mu = D^\mu/G^\mu$, $w^\mu$ descends to a well defined quasiconformal homeomorphism $w^\mu: \Sigma \rightarrow \Sigma^{\mu}$ and so represents a point $[\Sigma, w^\mu, \Sigma^\mu]$ in $\mathrm{Teich}(\Sigma)$. We call the map
$$ \pi_T: \mathrm{Belt}(\Sigma) \rightarrow \mathrm{Teich}(\Sigma)\;,  \qquad \tilde \mu \mapsto [\Sigma, w^\mu, \Sigma^\mu] \;,$$
the Teichm\"uller projection map. This map is by construction surjective: For any point
$[\Sigma,f,\Sigma']\in \mathrm{Teich}(\Sigma)$, $\pi_T(\mu_f) =
[\Sigma,f,\Sigma']$, where $\mu_f$ is the dilatation of $f$.\\ 

Since the map is surjective, the next thing is to get control over the preimages of points. Obviously this defines an equivalence relation, which is called \emph{Teichm\"uller equivalence}: Two Beltrami differentials are said to be equivalent iff they lie in the same fibre of $\pi_T$. Equivalence can conveniently be expressed in terms of
the solutions to the Beltrami equation as explained, e.g., in \cite{TCAToTS}, namely
\begin{eqnarray} \label{Belteq}\mu \sim \nu \qquad \Rightarrow \qquad w^\mu |_{\partial \bb D} = w^\nu|_{\partial \bb D}\;.\end{eqnarray}
Hence we can represent Teichm\"uller space as equivalence classes
of Beltrami differentials,
\begin{equation}\label{BeltramiModel}\mathrm{Teich}(\Sigma) \cong \cal T_B(\Sigma):=\mathrm{Belt}(\Sigma)/\sim \;.\end{equation}
The subscript $B$ stands for \emph{Beltrami model}. $\cal T_B(\Sigma)$ is a complex manifold with the complex structure induced from the complex structure on the Banach space $L^\infty_{(-1,1)}(\Sigma)$.\\

Observe that the equivalence condition \eqref{Belteq} implies that the
domains $w^\mu(\bb D)$ don't depend on the choice of representative
$\mu$ for the equivalence class. On the other hand they do depend on the way one extends $\mu$ to the plane, e.g., if $\mu$ is extended by 'reflection', $w^\mu(\bb D) = \bb D$ for all $\mu$. \\
%%%%%%%%%%%%%%%%%%%%%%%%%%%%%%%%%%%%%%%%%%%%%%%%%%%%%%%%%%%%%%%%%%%%%%%%%%%%%%%%%%%%%%%%%%%%%55
\subsubsection{The Moduli Space of Riemann Surfaces}
Although it will not occur often explicitly in this paper, the moduli space of Riemann surfaces always lurks in the background, and hence we will very briefly introduce it. It is the 'best approximation' to the base space of a universal object\footnote{As is well known, this category doesn't have a universal object. It is only a course moduli space, which is due to the fact that some surfaces have nontrivial automorphism groups.} in the category of families of Riemann surfaces, i.e., if we drop the condition that the families be marked. In fact, many mathematicians view Teichm\"uller space only as an intermediate object, the ultimate object of interest being the moduli space.\\

Set theoretically, the moduli space of a Riemann surface $\Sigma$ of finite analytic type, denoted by $\cal M(\Sigma)$, is given by a definition very similar to the one of $\mathrm{Teich}(\Sigma)$, namely as equivalence classes of quasiconformal mappings $(\Sigma, f, \Sigma')$ under moduli equivalence $\sim_m$:
$$ (\Sigma,f,\Sigma') \sim_m (\Sigma,g,\Sigma'')\quad \Leftrightarrow \quad \exists \;\mathrm{conformal} \;\sigma:\Sigma' \rightarrow \Sigma'' \quad \mathrm{s.t.} \quad \sigma \circ f = g\;.$$
In other words, moduli space is the space of complex structures on the surface $\Sigma$. Since all the maps are required to be quasiconformal homeomorphisms, they are in particular invertible, and hence the definition of Teichm\"uller space and the definition of moduli space can be rewritten as 
$$ \mathrm{Teich}(\Sigma) = \cal C(\Sigma) / \QC_0(\Sigma) \;, \qquad \cal M(\Sigma) = \cal C(\Sigma) / \QC(\Sigma)\;,$$
where $\QC(\Sigma)$ is the group of quasiconformal homeomorphisms of $\Sigma$, $\QC_0(\Sigma)$ the component of $\QC(\Sigma)$ containing the identity and $\cal C(\Sigma)$ the set of complex structures on the surface $\Sigma$. The point of this rewriting is that $\QC_0(\Sigma)$ is a normal subgroup of $\QC(\Sigma)$, so we can form the quotient $\mathrm{Mod}(\Sigma):= \QC(\Sigma) / \QC_0(\Sigma)$, called the \emph{Teichm\"uller modular group}\footnote{This group is closely related to the more well-known concept of a mapping class group: the \emph{mapping class group} of a topological space $X$, denoted by $\cal {MCG}(X)$, is given by the quotient of all orientation-preserving homeomorphisms by the subgroup of consisting of homeomorphisms isotopic to the identity. There is a natural map $\mathrm{Mod}(\Sigma) \rightarrow \cal {MCG}(\Sigma)$, which is an isomorphism for surfaces of finite analytic type. For compact surfaces this is easy to see: Any orientation-preserving homeomorphism is isotopic to a diffeomorphism, and any orientation-preserving diffeomorphism of a compact surface is of course automatically quasiconformal. For punctured surfaces one has to combine the above argument with removable singularity theorems for quasiconformal mappings, see \cite{TCAToTS}.}, and relate Teichm\"uller space to moduli space,
$$ \cal M(\Sigma) = \mathrm{Teich}(\Sigma) / \mathrm{Mod}(\Sigma)\;.$$
In this way, the moduli space becomes a topological space. But due to some surfaces having a non-trivial automorphism group, the action of $\mathrm{Mod}(\Sigma)$ on $\mathrm{Teich}(\Sigma)$ is \emph{not free}\footnote{In fact, the action is only free for $\Sigma = \bb C, \cinf$. For all other surfaces of finite type, there exists a quasiconformally eqiuivalent surface with automorphisms (see \cite{TCAToTS}, p.125).}. It is properly discontinuous, though, and hence moduli space carries the structure of a normal complex space\footnote{Of course one can consider quotients of Teichm\"uller space by \emph{normal subgroups} of the modular group, and in fact, the modular group always has a subgroup of finite index which acts freely.}. Since the automorphism group of stable surfaces is by definition \emph{finite}, it can also be thought of as a space with an orbifold complex structure. And finally, since Riemann surfaces are the same as projective algebraic curves over $\bb C$, there is also an algebraic construction of moduli space and in this framework it has been proven to be quasi-projective. \\

In fact, since Teichm\"uller space is topologically trivial, all the topology of moduli space is generated by the $\mathrm{Mod}(\Sigma)$-action on $\mathrm{Teich}(\Sigma)$. So even though we are only concerned with Teichm\"uller space, we will pay some attention to the action of the modular group when we talk about degenerations in part three.
%-------------------------------------------------------------------------------------
\subsection{The Bers Embedding of Teichm\"uller Space} \label{Bersemb}
%-------------------------------------------------------------------------------------
In this section we will see holomorphic realizations of Teichm\"uller spaces as bounded contractible domains in  complex vector spaces, via the so-called \emph{Bers embedding}. For this purpose, we will identify $\cal T_B(G)$ with a certain subset of $\cal O(\bb D)$ called the \emph{schlicht model} of Teichm\"uller space and denoted by $\cal T_S(G)$, which allows us to use methods of geometric function theory to study Teichm\"uller spaces.\\

The Bers embedding utilizes the \emph{Schwarzian derivative}, a classical non-linear differential operator, which we introduce in the first section.\\

The second section is then devoted to deformation spaces of Fuchsian groups, which in particular contain the sets $\cal T_S(G)$, which are the \emph{schlicht models} of Teichm\"uller spaces, and the third section brings the two concepts together and introduces the Bers embedding.\\

The reason for working with the schlicht model and the images under the Bers embedding instead of sticking to the Beltrami model will only become transparent in Section \ref{BersBound}. At the heart lies the fact that they are 
subsets of the deformation space of Fuchsian groups, so one can talk about boundary points and closures resp. completions, and also that the boundary points have a natural geometric interpretation again. In particular, a main object of interest later on will be \emph{augmented Teichm\"uller space}, which is a partial completion of the image of Teichm\"uller space under the Bers embedding and which is of fundamental importance, since it is a preimage of the Deligne-Mumford compactification of Moduli space for the quotient by the modular group.\\

As a technical remark, in this section we interchange the role of $\bb D$ and $\bb D^c$ in our considerations of Teichm\"uller spaces, i.e., we are working on the chart around infinity on the Riemann sphere. Hence the Beltrami differentials we discussed previously will now be supported on $\bb D^c$, and the associated quasiconformal solutions of the Beltrami equation will be a holomorphic functions on the disc.  
\subsubsection{Schlicht Functions and the Schwarzian Derivative} \label{SchDer}
The Schwarzian derivative is a remarkable non-linear differential operator which has been studied extensively. All of the unreferenced material in this section can be found, e.g., in \cite{UFaTS}.\\

In order to define the Schwarzian, let $f:D \rightarrow
\bb C$ be locally injective and thrice differentiable on a domain $D \subset \bb C$. The Schwarzian of $f$ is given by
$$S_f(z) : = \left(\frac{f''(z)}{f'(z)}\right)' -\frac{1}{2}\left(\frac{f''(z)}{f'(z)}\right)^2\;.$$
A simple computation shows that $S_f(z) \equiv 0$ on an open set iff $f \in \textrm{PSL}(2,\bb C)$ and that the chain rule for the Schwarzian of a composition is given by
\begin{eqnarray} \label{trafoSchw}
    S_{f\circ g} = (S_f \circ g)(g')^2 + S_g\;.
\end{eqnarray}
Hence $S_f = S_{1/f}$ since $1/f$ is of the form $M\circ f$ with $M\in \textrm{PSL}(2,\bb C)$, so we can
define the Schwarzian derivative on the set of locally injective meromorphic functions, which we denote by $\cal
M_{li}(D)$. To do so, for $f \in  \cal
M_{li}(D)$, set
$$S_f(p) = S_{\frac{1}{f}}(p)\;,$$ 
for all poles $p$ of $f$. Observe that local injectivity forces the pole
to be simple. The chain rule also allows us to define the Schwarzian acting on domains in $\cinf$ containing $\infty$ as follows: We set $\phi(z) = f(1/z)$ and define
$$ S_f(\infty) = \lim_{z\rightarrow 0} z^4 S_\phi (z)\;.$$
\begin{thm} \label{SchwSol}
    For any sc-domain $D \subset \cinf$, the Schwarzian derivative is a surjective operator \mbox{$S: \cal M_{li}(D) \rightarrow
    \cal O(D)$}. Moreover, the solution $f$ to the equation $S_f = \phi$ is unique up to post-composition by M\"obius transformations.
\end{thm}
\emph{Proof}. One easily sees that the previous definition implies that the Schwarzian of a meromorphic locally injective function is holomorphic. Now let
$y_1,y_2$ be two linearly independent solutions of the
differential equation
$$ y'' - \frac 12 \phi y = 0\;, \qquad \phi \in \cal O(D)\;.$$
Since the highest coefficient does not vanish and all coefficients are holomorphic, the above-mentioned solutions
necessarily exist. One checks by direct computation that
$$ S_{\frac{y_1}{y_2}} = \phi\;,$$
so we have shown surjectivity. Now any other pair of independent
solutions can be written as $\tilde y_1 := ay_1 + by_2$ and $\tilde
y_2:= cy_1 + dy_2$ with $a,b,c,d$ constituting a non-singular
$2\times2$ matrix. More precisely,
$$ \frac{\tilde y_1 }{\tilde y_2} = \frac{ay_1 + by_2}{cy_1 + dy_2}
=\frac{a\frac{y_1}{y_2} + b}{c\frac{y_1}{y_2} + d}\;,$$
which corresponds to a post-composition by a M\"obius transformation. $\Box$\\

Due to the local character of the Schwarzian derivative, this theorem has an interesting corollary about locally injective extensions of holomorphic functions 
\begin{cor}
   Let $D\subset \bb C$ be a domain and $f \in \cal M_{li}(D)$. If the Schwarzian derivative $S_f$ of $f$ has an extension $\hat F$ to a domain $D' \supset D$, then there is a locally injective extension of $f$, say $\hat f$, such that $S_{\hat f}= \hat F \in \cal O(D')$.
\end{cor}
\emph{Proof.} Let $\hat F$ be an extension of $S_f$ to $D'$. Then by the theorem above, there exists a \mbox{$\hat f\in \cal M_{li}(D')$} such that $S_{\hat f} = \hat F$. But since the Schwarzian is a local operator, the restriction of $\hat f$ to $D$ has the same Schwarzian as $f$, hence $f$ and $\hat f$ differ only by a M\"obius transformation on $D$. Applying this transformation to $\hat f$ we obtain a locally injective meromorphic extension of $f$, the Schwarzian of which is precisely $\hat F$. $\Box$\\

Let $\cal{S}(D) \subset \cal M_{li}(D)$ denote the set of \emph{injective holomorphic} functions. Elements of $\cal S(D)$ are called \emph{schlicht} functions, or by some authors also \emph{univalent} functions. It will become clear in the next section why we are especially interested in schlicht functions. We will start by quoting a basic but very useful theorem about schlicht functions.
\begin{thm}[Koebe's $\frac 14$-theorem]\label{K1QThm}
		Let $f$ be a schlicht function on the disc $\bb D$ that fixes the origin and with $f'(0)$ = 1. Then the image $f(\bb D)$ contains the disc $\bb D_{\frac 14}$.
\end{thm}
Recall that we introduced several norms on the space of automorphic forms\footnote{Here we interpret the holomorphic functions $\cal
O(D)$ as automorphic forms with respect to the trivial group so we can equip $\cal O(D)$ with any of the norms introduced in
Section \ref{BSoAF}.} on a domain in Section \ref{BSoAF}, one of these being the $B_2$-norm
$$ \| f\|_{B_2(D)}:= \textrm{sup}_{z\in D} \left\{|f(z)|\lambda^{-2}_D(z)\right\}\;.$$
A famous theorem by Kraus (rediscovered and often attributed to Nehari) says the following:
\begin{thm}[Kraus--Nehari Theorem]
Let f be a schlicht function on $\bb D$. Then
$$ |S_f(z)| (1-|z|^2)^2 \leq 6\;,\qquad \forall \; z \in \bb D\;,$$ and the bound is sharp. Since $\lambda_\bb D^{-2} = (1-|z|^2)^2$, this can be rewritten as
$$ \| S_f\|_{B_2(\bb D)} \leq 6 \qquad \forall\;f\in \cal S(\bb D)\;.$$
\end{thm} 
In fact, for an arbitrarily given sc-hyp domain $D$ there exists a constant $C(D) \subset [6,12]$ such that the same statement holds.\\

On the other hand, the Schwarzian provides a criterion for a function being schlicht. It is very remarkable since it shows the sensitivity of the Schwarzian to the phenomenon of quasiconformality, and on the other hand gives a further characterization of quasidiscs. The statement and a proof can be found in (\cite{UFaTS}, Section II.4.6).
\begin{thm}[Gehring] \label{SchlichtCrit} 
   Let $D$ be a sc-hyp domain. $D$ is a quasidisc iff there exists a constant $\delta>0$ such that for all $f \in \cal M_{\mathrm{li}}(D)$, 
   $\|S_f\|_{B_2(D)} \leq \delta$ implies that $f \in \cal S(D)$.
\end{thm} 
The norm of the Schwarzian induces a way to measure M\"obius invariant distances of domains, which is very much related to Teichm\"uller distance. For more on the geometry of the space of domains we refer the reader to \cite{UFaTS}.\\

Let us denote the set of image points of the schlicht functions under the Schwarzian by $\bb S$, i.e.,
$$ \bb S := S(\cal S(\bb D)) \subset B_2(\bb D)\;.$$
It is rather easy to establish that $\bb S$ is closed in $B_2(\bb D)$. Now let $\cal T_S(\bb 1) \subset \cal S(\bb D)$ denote\footnote{Actually the symbol $\cal T_S(\bb 1)$ also requires the functions to be 1-point normalized at the origin (see Def. \ref{TeichofGroup} below). This of course has no impact on the image under the Schwarzian.} the subset of schlicht functions which admit an extension to a quasiconformal homeomorphism of $\cinf$ and $T(\bb 1)$ its image under the Schwarzain. The following highly interesting theorem can be found in \cite{TCAToTS}.
\begin{thm} [Gehring] \label{IntClos}
    The interior of $\bb S$ is precisely $T(\bb 1)$ but the closure of $T(\bb 1)$ is a proper subset of
    $\bb S$.
\end{thm}
%========================================================================================================
\subsubsection{Deformation Spaces of Fuchsian Groups and Teichm\"uller Spaces}\label{DeformationSpaces}
In this section we approach Teichm\"uller spaces from a function-theoretic prespective, the prerequisites for which we have prepared in the previous section. The advantages are twofold: First of all, points of Teichm\"uller space will not be represented by equivalence classes of Beltrami differentials, as in the Beltrami model, but by a \emph{single schlicht function}, which is a much easier object to deal with and allows techniques from classical complex analysis and geometric function theory to be applied. The second main advantage is that it provides an extrinsic viewpoint of Teichm\"uller spaces as subsets of the so-called \emph{deformation spaces of Fuchsian groups}\\

The notion of a deformation space of a Fuchsian group was introduced by Kra in \cite{DoFG}, \cite{DoFG2}.
\begin{dfn}
   Let G be a Fuchsian group acting on $\bb D$. A \underline{deformation} of G is a pair $(\chi,f)$ where $\chi:G\rightarrow \mathrm{PSL}(2,\bb
   C)$ is a homomorphism and $f: \bb D \rightarrow \cinf$ a locally injective meromorphic function that satisfies the
   compatibility equation $$ f \circ g = \chi(g) \circ f\;, \qquad \forall g \in G\;.$$
   Two deformations $(\chi_1, f_1)$ and $(\chi_2, f_2)$ are called equivalent iff
   $$ \exists M \in \mathrm{PSL}(2, \bb C): \qquad f_2 = M \circ f_1\:, \qquad \chi_2(g) = M \circ \chi_1(g) \circ
   M^{-1} \;.$$
\end{dfn}
\begin{dfn}
	 The set of all equivalence classes of deformations of $G$ is denoted by $\mathrm{Def}(G)$ and is called the \underline{deformation space} of the Fuchsian group $G$.
\end{dfn}
The most general notion of a deformation of $G$ would be any homomorphism of $G$ into $\psl$. But
as we will see below in Proposition \ref{DefQDCorr}, defining a deformation as a pair as above yields a complex linear structure on the space of all deformations which will additionally tie together with geometric deformations of the Riemann surfaces.\\

Observe also that the data of the definition is somewhat redundant: The developing map $f$ determines $\chi$, since by local injectivity, for any point $z \in f(\bb D)$ there is a neighborhood $U_z$ where $f$ is invertible. Hence
$$ \chi(g) (w) = (f \circ g \circ f^{-1})(w) \qquad \forall \: w \in U_z\;,$$
and this determines $\chi(g)$ completely, since a M\"obius transformation is characterized by its value on three points. But also conversely, any $f \in \cal M_{\mathrm{li}}(D)$ \emph{is} the developing map of a deformation as we will see in the proof of \ref{DefQDCorr}. With this in mind, we will often identify functions in $\cal M_{\textrm{li}}(\bb D)$ with the corresponding deformations $(f, \chi)$ they
induce. For clarity in later use we will denote the forgetful maps by 
\begin{eqnarray} \label{DevHomDef}
\begin{array}{ll}
   \mathrm{dev}: \textrm{Def}(G) \rightarrow \cal M^0_{\textrm{li}}(\bb D)\;,& \qquad \mathrm{dev}([f,\chi]) = \tilde f\\
   \mathrm{hom}: \textrm{Def}(G) \rightarrow \mathrm{Hom}(G,\mathrm{PSL}(2,\bb C))\;,& \qquad \mathrm{hom}([f,\chi]) = \tilde \chi\;.
\end{array}
\end{eqnarray}
The upper index $0$ in $\cal M^0_{\textrm{li}}(\bb D)$ indicates that we consider the subset of functions which are 1-point normalized at the origin\footnote{Recall that this means $f(0) = 0 = f''(0)$ and $f'(0) = 1$. The upper $0$ will have the same meaning in the symbol $\cal S^0(\bb D$).}.
\begin{prop}\label{DefQDCorr}
   There following map is bijective and hence induces a $1-1$ correspondence
   between the equivalence classes of deformations of $G$ and the quadratic differentials of $G$:
   $$c_G: \mathrm{Def}(G)\rightarrow \cal Q(G):=\{f \in \cal O(\bb D): f = (f\circ g) (g')^2 \; \forall \: g \in G\}\;, \qquad c_G:= S \circ \mathrm{dev}\;.$$
\end{prop}
\emph{Proof.} In Theorem \ref{SchwSol} we established the correspondence of $\mathrm{Def}(\bb 1) \cong \cal M^0_{li}(\bb D)$ with $\cal Q(\bb 1) \cong \cal O(\bb D)$. As a next step, we show that $c_G([f,\chi]) = S_{\mathrm{dev}([f,\chi])}$
is a quadratic differential for $G$ whenever $[f,\chi] \in \mathrm{Def}(G)$. Recall the transformation behaviour
(\ref{trafoSchw}) for $g \in$ PSL$(2,\bb C)$,
$$ S_{f\circ g} = (S_f \circ g)(g')^2\;, \qquad S_{g\circ f} = S_f\;.$$
If $g \in G$ we have $f\circ g =  \chi(g) \circ f$ so altogether we
get $$ S_f = (S_f \circ g)(g')^2\;,$$ i.e. the Schwarzian of the
developing map behaves like a quadratic differential for $G$. Now let $f \in \cal M_{\mathrm{li}}(\bb D)$ be a funcion such that $S_f = \phi $ is a quadratic differential for the group $G$. Then for all $g \in G$,
$$ S_{f \circ g} = (S_f \circ g)( g')^2 = (\phi\circ g) ( g')^2 = \phi\;,$$ 
and hence by the uniqueness part of the solution theorem of the Schwarzian differential equation there exists a M\"obius transformation $\chi(g)$ such that $\chi(g) \circ f = f \circ g$. This association is a homomorphism,
$$ \chi(g_1g_2) \circ f = f \circ (g_1g_2) = (f \circ g_1)\circ g_2 = \chi(g_1) \circ f \circ g_2 = \chi(g_1)\chi(g_2) \circ f\;,$$
and the pair $(\chi,f)$ therefore a deformation of $G$. $\Box$ \\

Now we can introduce a subset of the deformation space called the schlicht model of Teichm\"uller space.
\begin{dfn} \label{TeichofGroup}
   The \underline{schlicht model of Teichm\"uller space} of a Fuchsian group G acting on $\bb D$, denoted by $\cal T_S(G)$, is given by
   $$ \cal T_S(G):=\left\{
        \begin{array}{l}
            f \in \cal S^0(\bb D):\; \exists \: \textrm{qc.-hom.}\; w: \cinf \rightarrow \cinf\; \textrm{s.t.} \; w|_{\bb D}=f \\
            \textrm{and $w$ is compatible with $G$, i.e. $wGw^{-1}$ is again Kleinian}
        \end{array}\right\} \;.$$
   The preimage on $\cal T_S(G)$ under $\mathrm{dev}$ is denoted by $\mathrm{Def}_{\mathrm{QF}}(G)$.
\end{dfn}
The following proposition justifies the name of the set $\cal T_S(G)$.
\begin{prop} \label{CorrOfTwoModels}
   The map $\pi_{BS} : \cal T_B(G) \rightarrow \cal T_S(G)$ defined by $\pi_{BS}([\mu]) = w^\mu_{|\bb D}$ is a bijection.
\end{prop}
\emph{Proof}. First of all, the map is well defined and injective, since Teichm\"uller equivalence was defined as
\begin{eqnarray}\label{BeltEquiv} \mu \sim \nu \qquad \Longleftrightarrow \qquad w^\mu|_{\bb D} = w^\nu|_{\bb D}\;,\end{eqnarray}
in the last section. The map is also surjective: Let $f \in \cal T_S(G)$ be given and $w$ be any quasiconformal extension of $f$. 
Recall that $wGw^{-1}$ is again Kleinian iff the complex dilatation of the map $w$ is a Beltrami differential with
respect to the group $G$ acting on $\bb D$. Hence $\pi_{BS}^{-1}(f) = \hat \mu (f)$. $\Box$ \\

\subsubsection{The Bers Embedding}
In the previous section we saw a 1-1 correspondence of the Teichm\"uller space of a group $G$ with a subset of the space of quadratic differentials for $G$ on $\bb D$. This complex vector space can be equipped with the $B_2$-norm, and we have already seen  that this nowm is intimately connected with the Schwarzian derivative and schlicht functions. In the present section we will show that we are dealing with much more than a correspondence of sets, namely with a biholomorphism. Let us start by giving the map a name.
\begin{dfn}\label{BersEmbeddingDef}
   Let $G$ be a Fuchsian group acting on $\bb D$. The \underline{Bers embedding} $\beta$ of Teichm\"uller space is given by
   $$ \beta: \cal T_B(G) \rightarrow B_2(\bb D, G)\;, \qquad \beta = S \circ \pi_{BS}\;.$$
   The image $\beta(\cal T_B(G))$ will be denoted by $T(G)$. The lift of this map,
   $$\tilde \beta : \bb D_1\left(L^\infty_{(-1,1)}(\bb D^c, G)\right) \rightarrow  B_2(\bb D, G)\;, \qquad \tilde \beta := \beta \circ \pi_T$$ is 				 called the \underline{Bers projection}.
\end{dfn}
The target space of the map as stated above is correct, since \mbox{$\mathrm{dev}^{-1}(\cal T_S(G)) = \mathrm{Def}_{\textrm{QF}}(G) \subset \mathrm{Def}(G)$}, which by \ref{DefQDCorr} is mapped to the quadratic differentials for $G$. Moreover, since the elements of $\cal T_S(G)$ are schlicht functions, we have
$$ T(G):=S(\cal T_S(G)) \subset \bb D_6\left(B_2(\bb D,G)\right)\;.$$ \\

Let us remark that The Beltrami model of Teichm\"uller space behaves contravariantly under the inclusion of Fuchsian groups $G' \subset G$, by what we mean that $G'\subset G$ implies $\cal T_B(G) \subset \cal T_B(G')$. This is obvious from the corresponding isometric inclusion of the vector spaces 
$$ L^\infty_{(-1,1)}(\bb D^c,G) \hookrightarrow L^\infty_{(-1,1)}(\bb D^c,G')\;.$$
The Bers embedding \emph{respects this structure}, because we also have isometric embeddings of the target spaces,
\begin{eqnarray}\label{BoundIsom} B_2(\bb D, G) \hookrightarrow B_2(\bb D, G')\;.\end{eqnarray}
Observe that the isometry in \eqref{BoundIsom} \emph{only} works for the $L^\infty_2$-norm, but not for any finite $p$ due to the fact that the group $G$ only enters the norms in terms of a fundamental domain for its action. Whereas the $p$-norms depend on the fundamental domain, because integration is performed over them, the $\infty$-norm does not see the fundamental domain. In particular, since the hyperbolic volume of the disc is infinite, there would be no \emph{universal} Teichm\"uller space in the $p$-norm setting.\\

The following theorem due to Bers justifies calling the map an embedding. For a detailed proof we refer to the textbook \cite{TCAToTS}.
\begin{thm} \label{BembThm}
    Let $G$ be a Fuchsian group acting on $\bb D$. The Bers projection $\tilde \beta$
    is a holomorphic submersion and factors precisely through $\pi_T$ yielding a holomorphic embedding of $\cal T_B(G)$ into $B_2(\bb D, G)$ as 			an open, bounded and contractible domain.
\end{thm}
The theorem has a wealth of consequences. Recall that by Theorem \ref{SchlichtCrit}, there exists a constant $\delta$ such that the $\delta$-ball in $B_2(\bb D)$ around the origin is contained in $\bb S$, hence the $\delta$-ball around the origin in $B_2(\bb D, G)$ is also contained in $\bb S \cap B_2(\bb D, G)$. Further, by Theorem \ref{IntClos}, 
$$ T(G) = T(\bb 1) \cap B_2(\bb D, G) = \textrm{int}(\bb S) \cap B_2(\bb D, G) \supset \bb D_\delta\;,$$
so $\mathrm{dim}_\bb C T(G) = \mathrm{dim}_\bb C B_2(\bb D, G)$. If $G$ is a Fuchsian group of first kind without elliptic elements it represents a surface of type $(g,n)$. By the Riemann--Roch Theorem \ref{R-R}, we get:
\begin{cor} \label{DimTeich} 
Let $G$ be a Fuchsian group of first kind uniformizing a surface of type $(g,n)$. Then  
    $$\mathrm{dim}_\bb C T(G) = 3g-3 +n  \;,$$
whereas if $G$ is Fuchsian of second kind, $T(G)$ is infinite-dimensional.
\end{cor}
On the other hand, the holomorphic shape of $T(G)$ is very intriguing and not yet fully understood, but it is known to be of fractal nature. As another example, which is interesting from the point of view of geometric function theory, let us look at the two quantities
$$ i(G) := \sup_{\delta \in \bb R} \{\bb D_\delta \subset T(G)\}\;, \qquad o(G):= \inf_{\delta \in \bb R} \{ \bb D_\delta \supset T(G)\}\;,$$
called the \emph{inradius} resp.~ the \emph{outradius} of the Teichm\"uller space $T(G)$. They obviously satisfy $i(G)\geq 2$ and $o(G)\leq 6$. The following facts concerning the in- resp.~ outradius can be found in \cite{OIRoTS} resp.~ \cite{OtOotTS}: $i(G)$ is strictly greater than two for any finitely generated Fuchsian group $G$ of first kind but there exists a sequence of groups $G_i \subset \mathrm{Def}_{\mathrm{QF}}(G)$ such that $i(G_i)\rightarrow 2$ for $i\rightarrow \infty$. $o(G)$ equals $6$ for Fuchsian groups of second kind and $o(G)$ is strictly less than $6$ for finitely generated Fuchsian groups of first kind. Yet given a finitely generated Fuchsian group of first kind $G$, there exists a sequence $G_i$ of groups in $\mathrm{Def}_{\mathrm{QF}}(G)$ such that $o(G_i)\rightarrow 6$ for $i\rightarrow \infty.$\\
%-------------------------------------------------------------------------------------
\subsection{Fibre Spaces over Teichm\"uller Space}
\label{FibreoverTeich}
%-------------------------------------------------------------------------------------
In this section we will introduce two fibre spaces over Teichm\"uller
space, $\cal F(G)$ and $\cal V(\Sigma)$, and explain their geometric and conceptual significance. While $\cal F(G)$ will be very useful in analytic constructions of objects in families, its conceptual importance is due to the fact that the related fibre space $\cal V(\Sigma)$, called the \emph{universal curve}, is the universal object in the category of marked families of Riemann surfaces of finite analytic type.\\

As before, let the given hyperbolic Riemann surface $\Sigma$ be uniformized as $\Sigma=\bb D/G$ so we can identify $\mathrm{Teich}(\Sigma)$ with $\cal T_B(G)$.\\ 

Let us start with the definition of the fibre space $\cal F(G)$, which, among others, was extensively studied by Bers in \cite{FSoTS}. 
To every point in $\cal T_B(G)$ we attach the
appropriate quasidisc as a fibre,
\begin{eqnarray*} \label{FibreSpaceDefforlater} \cal F(G) \subset \cal T_B(G) \times \bb C\;, \quad \cal F(G):=\left \{([\mu], z): z \in w^\mu(\bb D)\right\}\;.\end{eqnarray*}
The projection map $\tilde \pi_\Sigma: \cal F(G)\rightarrow \cal T_B(G)$ onto the first factor is obviously
holomorphic. The group $G$ viewed as an abstract group acts
on $\cal F(G)$ by
$$ G \times  \cal F(G) \rightarrow  \cal F(G)\;, \qquad
(g,([\mu], z)) \mapsto \left([\mu], g^{\mu} (z) := w^\mu \circ g \circ (w^{\mu})^{-1}(z) \right)\;.$$ The action
is holomorphic and preserves the fibres. Moreover, if $G$ does not contain elliptic elements, the quotient $\cal V_\bullet(\Sigma):=\cal F(G)/G$ again
is a complex manifold and the projection $\tilde \pi_\Sigma$ factors
through
$\pi_\Sigma: \cal V_\bullet(\Sigma) \rightarrow \cal T_B(G)$.\\

It is tautological in the sense that the preimage of a point
$t=[\Sigma,f,\Sigma'] \in \cal T_B(G)$ is precisely the Riemann surface
it represents, i.e., $\tilde \pi_\Sigma^{-1}(t) \cong \Sigma'$ which may be non-compact depending on the original Riemann surface $\Sigma$. Therefore $\tilde \pi_\Sigma$ is
not a proper map if $\Sigma$ is non-compact, and $\cal V_\bullet(\Sigma)$ is not a family of complex
manifolds in the sense of Kodaira and Spencer. \\

At least when $\Sigma$ is of finite analytic type there is a canonical way of obtaining a compact surface $\widetilde \Sigma$ by inserting the missing points. Hence we obtain a family $\pi_\Sigma:\cal V(\Sigma) \rightarrow \cal T_B(G)$, where the fibre of a point is the compact `filled-in' surface $\widetilde \Sigma'$ of the surface $\Sigma'$ it represents.\\
\begin{dfn}
The family $\pi_\Sigma: \cal V(\Sigma) \rightarrow \cal T_B(G)$ is called the
\underline{universal curve}.
\end{dfn}
The universal curve of a compact surface of genus $g\geq 3$ has \emph{no holomorphic sections}, i.e., there doesn't exist a holomorphic map $s:\cal T_B(G) \rightarrow \cal V(\Sigma)$ such that $\pi_\Sigma\circ s = \bb 1_{\cal T_B(G)}$. This was proved by Hubbard in \cite{SlsadlcudT}.\\

On the other hand, if $\Sigma$ has $n$ punctures, the projection $\pi_\Sigma$ has $n$ disjoint sections $s_i:\cal T_B(G) \rightarrow \cal V(\Sigma)$ that correspond to the position of the punctures. This is possible, in contrast to the compact case, because the universal curve $\cal V(\Sigma)$ is larger\footnote{The dimension of the base differs by $n$, see Corollary \ref{DimTeich}.} than the universal curve $\cal V(\widetilde \Sigma)$ although they have the same fibres topologically. \\

The following important theorem justifies calling $\cal V(\Sigma)$ the universal
curve. For compact surfaces it was first proven by  Grothendieck \cite{TdCeGA}, and the generalization to surfaces of finite type was proven by Engber in \cite{TSaRoF}. A very detailed and recommendable discussion of universality properties of fibre spaces over Teichm\"uller space can be found in (\cite{TCAToTS}, Ch.~ 5).
\begin{thm} \label{UniversalityOfUniversalCurve}
   Let $\pi: \cal X \rightarrow B$ be an arbitrary family of marked Riemann
surfaces of finite type and let $\tilde {\cal X} \rightarrow B$ be the family of compact marked Riemann surfaces with marked points obtained by filling in the punctures by marked points in all fibres. Then there exists a unique holomorphic map $f: B \rightarrow \cal T_B(G)$ such that $\tilde {\cal X} \cong f^*\cal V(\Sigma)$ in the sense of families of marked Riemann surfaces with marked points.
\end{thm}
We now understand how the families $\cal V(\Sigma)$ and $\cal V(\widetilde \Sigma)$ are related: There is a unique map $F: \cal T_B(G) \rightarrow \cal T_B(\widetilde G)$ such that $\cal V(\widetilde \Sigma) = F^*\cal V(\widetilde \Sigma)$. The fibre $F^{-1}(p)$ of a point in $\cal T_B(\widetilde \Sigma)$, which corresponds to a marked Riemann surface $\Sigma_p$, is given by
$$ F^{-1}(p) = \left\{ (\Sigma_p ,P_1, \ldots, P_n)\mid \{P_1,\ldots,P_n\} \in (\Sigma_p)^n \backslash \Delta \right \} \;,$$
with the `same' marking, i.e., choice of homotopy classes of loops generating $\pi_1(\Sigma_p, P)$. By $\Delta$ we simply mean the union of all diagonals, i.e., the points $\{P_1, \ldots, P_n\}$ in each $n$-tuple have to be distinct.
%-------------------------------------------------------------------------------------
%-------------------------------------------------------------------------------------
%-------------------------------------------------------------------------------------
%-------------------------------------------------------------------------------------
%-------------------------------------------------------------------------------------
\subsection{Families of Factors of Automorphy}
\label{faforms}
%-------------------------------------------------------------------------------------
In this short section we are going to bring together several of the previous ingredients to show that there is a unique way to extend line bundles over a surface to families of line bundles over deformations of the surface, namely using the $s$-factor representation. We assume $\Sigma$ to be of finite type, and as we commented earlier, we present the construction only for the universal family of marked Riemann surfaces, since it is functorial under pull-back.\\

Suppose we are given a family of line bundles over the universal curve, i.e., a holomorphic line bundle $\cal L
\rightarrow \cal V(\Sigma)$. Since $\cal V_\bullet(\Sigma) \cong \cal F(G) / G$ and line bundles over $\cal F(G)$ are trivial, we can
describe $\cal L$ in terms of a factor of automorphy for the action of $G$ on $\cal F(G)$, i.e., a holomorphic map $\cal R: G \times \cal F(G) \rightarrow
\bb C^*$ satisfying the cocycle condition. \\

For fixed point $t \in \cal T_B(G)$ in the base we get a factor of automorphy on the fibre $X_t$ representing $\cal L_t$,
which we denote by $\rho^t_g(z):=\cal R_g(t,z)$. A section $\Psi \in \Gamma(\cal L)$ corresponds to an automorphic form on $\cal F(G)$ with respect to $\cal R$,
$$ \Psi(g(t,z))= \Psi(t,z)\cal R_g(t,z)\;,$$
but since the action of $G$ on $\cal F(G)$ preserves the fibres, this equation reads
$$ \Psi(t, g^t(z)) = \Psi(t,z) \rho^t_g(z)\qquad \forall \; t\;,$$
i.e., $\psi^t(z):= \Psi(t,z)$ is an automorphic form with respect to $\rho^t$. Of course equivalence
classes of them are again in 1-1 correspondence to isomorphism classes of line bundles over $\cal V(\Sigma)$. \\

In Lemma \ref{uf}, we showed that a flat factor is always equivalent to a unitary flat one, i.e., that any line bundle
of degree zero can be represented by a 0-factor. The following can be considered as the family version that lemma. 
\begin{prop} \label{AboutUFFs}
    Let $\cal L_{t_0} \rightarrow X_{t_0}$ be a line bundle of degree zero over a Riemann surface $X_{t_0}$, thought of as the fibre of $\cal V(\Sigma)$
    over $t_0 \in \cal T_B(G)$. Then
    \begin{enumerate}
        \item{Any deformation $\cal L \rightarrow \cal V(\Sigma)$ of $(X_{t_0}, \cal L_{t_0})$ over $\cal T_B(G)$ can be
        represented by a factor $\cal R$ such that $\cal R(t_o, \cdot):=\rho^{t_0}$ is unitary flat.}
        \item{There is a unique deformation of $\cal L \rightarrow \cal V(\Sigma)$ of $(X_{t_0}, \cal L_{t_0})$ over $\cal T_B(G)$
         such that $\cal L$ can be represented by a factor $\cal R$ inducing unitary flat factors on all fibres.}
    \end{enumerate}
\end{prop}
\emph{Proof.} Any such deformation is given by a line bundle $\cal L \rightarrow \cal V(\Sigma)$, which again is given in
terms of a factor of automorphy $\cal R$ acting on $\cal F(G)$. Now this factor induces a factor on the fibre $D^{t_0}$
over $t_0$ by $\rho^{t_0}:=\cal R(t_0,\cdot)$ which is equivalent to any factor representing $\cal L_{t_0}$. This factor
can, by Lemma \ref{uf}, be represented as a 0-factor by conjugating with a function of the form
$$ h(z):= e^{2 \pi i \sum_{ij} C^{ji}\sigma_i w_j(z)}\;.$$
The function $h(z)$ can now be extended to a function $\tilde h(z,t)$ on $\cal F(G)$ by replacing the $w_j(z)$ by the
functions $\alpha_j(z,t)$ defined by Bers in \cite{HDaFotM}, i.e., automorphic
functions representing the canonical basis of holomorphic 1-forms. Therefore $\alpha_j(z,t_0)= w_j(z)$ and the factor
$\tilde{\cal R}$ obtained by conjugation with $\tilde h$ is unitary flat over $t_0$.\\

To prove the second statement, observe that the constant factor $\cal R(t, \cdot):= \rho^{t_0}$ proves the existence
and the uniqueness stems from the fact that unitary flatness says that the modulus satisfies $|\cal R_{g}(t)|=1$
for all $t$, and since $\cal R$ is holomorphic it is necessarily constant. $\Box$\\

Intuitively this is true in any degree, not only in degree zero. To establish this, recall the natural $s$-factors given by powers of the canonical factor,
$$ \kappa^s_g(z):= g'(z)^s\;,$$
for any $s \in \bb R$ such that $(2g-2)^{-1}s \in \bb Z$. Since over $\cal T_B(G)$ the group elements $g^t \in G^t$ depend holomorphically\footnote{Think of $t$ some complex coordinates on $\cal T_B(G)$.} on $t$, one easily checks that they fit together to yield a factor on $\cal F(G)$ that will be denoted by $\cal C^s$. For
non-integer $s$, these factors of course depend on the additional choice of the root of $g'$, but since they
are always discrete we get the following:
\begin{lemma}\label{AboutCan}
    Let $\cal K_{t_0}^s \rightarrow X_{t_0}$ be an $s$-th power of the canonical bundle over the Riemann surface $X_{t_0}$.
    Then there is a unique deformation $\cal K^s \rightarrow \cal T_B(G)$ of $(X_{t_0},\cal K_{t_0}^s)$ such that
    the line bundle induced on any fibre is isomorphic to the $s$-th power of the canonical bundle of the fibre.
    Moreover, $\cal K^s$ can be represented by a factor, denoted $\cal C^s$, which induces an $s$-factor on every fibre.
\end{lemma}
From this, together with proposition \ref{AboutUFFs}, we get the statement for arbitrary degree.
\begin{thm}\label{UniqueBundle}
    Let $\cal L_{t_0} \rightarrow X_{t_0}$ be a line bundle over a Riemann surface. Then there is a unique
    deformation $\cal L \rightarrow \cal V(\Sigma)$ of $(X_{t_0}, \cal L_{t_0})$ over $\cal T_B(G)$ such that $\cal L$
    can be represented by a factor $\cal R$ inducing s-factors on all fibres.
\end{thm}
\emph{Proof.} Let $\cal L_{t_0}$ be of degree $d$. The existence follows from noting that $\cal L_{t_0}$ can be
represented as a product of a flat bundle and any choice of an $s=(2g-2)^{-1}d$-th root of the canonical bundle of
$X_{t_0}$. Now both factors have natural deformations in terms of s-factors by Proposition \ref{AboutUFFs} and Lemma
\ref{AboutCan}, so their tensor product yields an s-deformation of $\cal L_{t_0} \rightarrow X_{t_0}$. Uniqueness
follows from the group structure of s-factors together with the previous uniqueness statements.
%%%%%%%%%%%%%%%%%%%%%%%%%%%%%%%%%%%%%%%%%%%%%%%%%%%%%%%%%%%%%%%%%%%%%%%%%%%%%%%%%%%%%%%%%%%%%%%%%555
%%%%%%%%%%%%%%%%%%%%%%%%%%%%%%%%%%%%%%%%%%%%%%%%%%%%%%%%%%%%%%%%%5
\addtocontents{toc}{\vspace{.2cm}}
\setcounter{section}{3} \setcounter{subsection}{0}
\subsection{Extensions of Automorphic Forms over Teichm\"uller Space}\label{aghs}
In this section, we finally address the issue of extending automorphic forms given on a surface to any marked deformation, i.e., any marked family of Riemann surfaces containing the surface as a fibre. There are several ways to do this, but as far as we know only one way to extend the sections holomorphically in both fibre and base direction. We describe this in Thm. \ref{ExtendSections}. Unfortunately this method only works for $s\geq2$.

We also deal with the more general problem of extending a given basis of the space of automorphic forms on a fibre. It will turn out that the extended elements will yield a basis on \emph{almost all} fibres, more precisely, it may fail to be a basis on an analytic subset of positive codimension (Thm. \ref{BasisAE}). And again, we only do this for the universal family because of aforementioned reasons.\\

For the rest of the paper $\Sigma\cong \bb D/G$ be a given surface of finite analytic type and $\rho_s$ an $s$-factor for the action of $G$ on $\bb D$. We identify $\Sigma$ with the 0-fibre of $\mathrm{Teich}(\Sigma)\cong T(G)$, which is the universal marked deformation space of $\Sigma$, and denote by $t$ some global analytic coordinate on $T(G)$.\\

Recall that for a given factor of automorphy $\rho_s$ there is a unique $s$-deformation with base $\cal T_B(G)$ by Theorem \ref{UniqueBundle}, i.e., a factor $\cal R$ acting on $\cal F(G)$ such that the restriction to the fibre over $t \in T(G) \cong \cal T_B(G)$ is an $s$-factor on the quasi-disc $D^t$ for the group $G^t$. \\

There are two natural questions we want to solve in this section, the second one of course depending on a positive solution of the first one. In fact, on an abstract level the answers are easily
obtained, as we explain below. We are looking for a concrete way of constructing the objects of concern.
\begin{enumerate}
   \item{Given an element $\psi \in A^1_{\rho_s}(\bb D, G)$, is there a holomorphic function $\Psi$ on $\cal F(G)$
   such that $\Psi|_{D^t} \in A^1_{\rho^t_s}(D_t, G^t)$ for all $t$ and $\Psi|_{\bb D} = \psi$? How do we construct it?}
   \item{Given a basis $\{\psi_1,\ldots,\psi_N\}$ of $A^1_{\rho_s}(\bb D, G)$, is there an extension $\Psi_i$ of each
   $\psi_i$ such that $\{\Psi_1|_{D^t}, \ldots, \Psi_N|_{D^t}\}$ is a basis of $A^1_{\rho^t_s}(D^t, G^t)$ for all $t \in T(G)$? }
\end{enumerate}
Both problems have been studied and implicitly solved by Earle and Kra in
\cite{PDaPSoVRS} and \cite{HCDoVRS}. We follow their line of argument, but also give an explicit recipe for the construction. Also, our notation differs heavily from theirs.\\

The following corollary of the famous Theorem B by Cartan motivates the attempt at finding such functions $\Psi$.
\begin{thm}
    Let $X$ be a Stein manifold, $A$ an analytic subspace and $f\in \cal O(A)$. Then there is an $F\in \oh(X)$ such that
    $F|_A = f$.
\end{thm}
Since Teichm\"uller space is known to be a contractible domain of
holomorphicity, in particular Stein, the first problem is reduced to
the question of whether there is an automorphic extension among all
extensions. \\

Our strategy goes as follows:
\begin{itemize}
   \item{Let $\psi \in A^1_{\rho_s}(\bb D,G)$ be given. Since $\Theta_{\rho_s}: A^1_{s}(\bb D)\rightarrow A^1_{\rho_s}(\bb D, G)$ is surjective, there exist preimages $h_\psi$, i.e., functions such that $\Theta_{\rho_s} [h_\psi] = \psi$. In fact, since the source is infinite-dimensional and the target finite-dimensional, there are infinitely many such preimages.}
   \item{In the next step, we try to find such an $h$ among all the preimages, with the additioinal property  that it has an extension $H: \cal F(G) \rightarrow \bb C$ and such that the restriction of $H$ to all fibres lies in $A^1_s(D_t)$.}
   \item{Now we can apply $\Theta_{\rho_s^t}$fibrewise. We obtain functions in $A^1_{s,\rho^t}(D_t, G^t)$ for all          $t$.}
\end{itemize}
The function we obtain in the last step, which we can write more suggestively as 
\begin{eqnarray} \label{ThetaFamily}
 \Psi(z,t) := \Theta_{\cal R}[H](z,t) = \Theta_{\rho_s^t}[H|_{D^t}](z)\;,
\end{eqnarray}
will be a holomorphic automorphic extension, as we show in Proposition \ref{ReqFamTheta}. The trick is finding an appropriate function $H$ in the second step.\\

For this we have to come up with a method producing a
function on $D^t$ for varying $t$ from a given one on $\bb D$, say
$h:D\rightarrow \bb C$. The only obvious functorial and elegant way
to do this is by pull-back, and the only functions $D^t
\rightarrow \bb D$ we have at our disposal are the inverse functions of the quasiconformal maps
$w^t: \bb D \rightarrow D^t$ and eventually a non-canonical family
of Riemann mappings. The latter ones are no good here, since they
definitely don't depend holomorphically on $t$. The quasiconformal
homeomorphisms, however, do yield holomorphic functions in $t$ for
fixed $z$, yet of course the $z$-direction will only be
quasiconformal.\\

To remedy the non-holomorphicity one could apply the projection operator $\beta_s$ from
Section \ref{BSoAF}, $\beta^{D^t}_s: L^1_s(D^t) \rightarrow A^1_s(D^t)$, fibrewise, to get holomorphic functions
$$\beta^{D^t}_s (h \circ (w^t)^{-1}): D^t \rightarrow \bb C$$
on all fibres. But this procedure destroys the holomorphic dependence on $t$, since the family of operators $\beta^{D^t}_s$, being integral operators with kernel function $K_{D^t,s}$, are not holomorphic in $t$, because the  kernel function does not depend holomorphically\footnote{The kernel function is proportional to a power of $\lambda_{D^t}$ depending on the value of $s$. The Poincar\'e density, however, depends only real-analytically on $t$, since it transforms with (powers of) the \emph{modulus} of the derivative of a Riemann mapping.} on $t$! Achieving holomorphicity in both Teichm\"uller and fibre directions, i.e., on all of $\cal F(G)$, is the real difficulty. \\

The solution to the problem is so trivial that it is easy to miss it at first glance: We simply use the \emph{same} function $h$ for all fibres. This trick, however, only works for $s \geq 2$.
\begin{prop} \label{ReqFamTheta}
    Let $\cal R$ be an s-factor with $s \geq 2$, let $\epsilon > 0$ and $h \in \oh(\bb D_{4 + \epsilon})$.
    Then $\Psi(\cdot,t) := \Theta_{\rho_s^t}[h_{|D^t}] \in A^1_{\rho_s^t}(D^t, G^t)$ for all $t \in T(G)$.
\end{prop}
\emph{Proof}. By the Koebe-1/4 Theorem, the quasidiscs $D^t$ are all contained in the disc of radius 4, and if $h$ is holomorphic on any larger disc $\bb D_{4+\epsilon}$, it is of course \emph{bounded} on the disc of radius 4. Since for $s\geq 2$, boundedness implies
integrability (see Lemma \ref{TwoMinusS}), $h|_{D_t} \in A^1_s(D^t)$ for all $t \in T(G)$. If we now define
$$ \Psi(t,z) := \Theta_{\cal R}[h](z) := \sum_{g \in G}(\rho^t_{g})^{-1}(z) (h \circ g) (z)\;,$$
then $\Psi(t,z)$ is a holomorphic  automorphic extension of the function $\psi = \Theta_{\rho_s}[h]$, since $\rho^t$ is holomorphic in $t$ and the sum converges normally for each $t$ by Theorem \ref{PS}, hence we may differentiate term by term. $\Box$ \\

Now this, together with the Theorem of Knopp \ref{Knopp}, gives us an affirmative answer to the first problem, since any $D^t$ is a Jordan region, $A^1_{\rho^t_s}(D^t,G^t)$ is finite-dimensional, and since polynomials are entire functions.
\begin{thm} \label{ExtendSections}
    For any s-factor $\rho_s$ for a Fuchsian group $G$ of the first kind acting on $\bb D$ with $s\geq 2$,
    any $\psi \in A^1_{\rho_s}(\bb D, G)$ is the restriction of an element $\Psi \in A(\cal R)$. Moreover, this $\Psi$ can be written as the Poincar\'e series of a fixed polynomial.
\end{thm}
This is in fact the \emph{only place we need} $s\geq 2$! Recall from Section \ref{IotPD} that the integrability of $\lambda^{2-s}$ depended heavily on the shape of the domain $D$, and in general the smallest power for which $\lambda^{2-s}$ is integrable is a function $s_0(D) \in [1,2]$. For the unit disc, or any other rectifiable domain, one has $s_0 = 1$, but except for the 0-fibre, \emph{all other fibres} are quasi-discs with Hausdorff dimension of their boundary greater than one \cite{UFaTS}.\\

%=============================================================================================================
By our general considerations, we have attached a vector space to any point $t\in T(G)$ for any factor $\cal R$, namely
$A^1_{s,\rho^t}(D_t, G^t)$,
$$ \pi_{\cal R}: \cal V(\cal R):=\bigcup_{t \in T(G)} A^1_{\rho^t_s}(D^t, G^t) \rightarrow
T(G)\;.$$
\begin{thm}\label{HermHolBundle}
    The space $\pi_{\cal R}: \cal V(\cal R)\rightarrow T(G)$ is a
    holomorphic vector bundle. It carries a natural smooth hermitean structure given by the Weil-Petersson pairing $\langle \cdot,\cdot\rangle_s^{G^t}$     in each fibre, which we denote by $\mathrm{WP}_{\cal R}$. 
\end{thm}
\emph{Proof.} We need to show local triviality. To this end, let $t
\in T(G)$ and $\{\psi_1,\ldots,\psi_N\}$ be a basis of the fibre.
By theorem \ref{ExtendSections}, the $\psi_i$ can locally be extended,
and since they are linerly independent on the fibre over $t$, they
are also linearly independent on the fibres over a neighborhood $U_t$ of $t$. This yields a local trivialization
$$ \cal V(\cal R)|_{\pi_{\cal R}^{-1}(U_t)} \cong U_t \times \bb
C^N\;.$$ Now if two such trivializations overlap on $V\in T(G)$,
i.e., if we have two $N$-tuples of linearly independent functions
$\{\Psi_i, \Phi_j\}$, they can without loss of generality be seen as
extensions of the same basis of the fibre over a point $t_0 \in V$
by taking appropriate linear combinations. One then easily sees
that the transition matrix is holomorphic. \\

The Weil-Petersson product induces a hermitean structure, since for a quasi-Fuchsian group of first kind, the spaces $A_{\rho_s^t}^p(D^t,G^t)$ are isomorphic for any $1\leq p\leq \infty$ by Theorem \ref{incl}. For $p=2$ the Weil-Petersson pairing is a product and therefore the fibres are Hilbert spaces. Let $\Psi_1,\Psi_2$ be local holomorphic sections given by $\Psi_i = \Theta_{\cal R}[h_i]$ for some functions $h_i$. Then their product is given by
$$\mathrm{WP}_{\cal R} (\Psi_1,\Psi_2)(t) := \langle \Psi_1, \Psi_2 \rangle_s^{G^t} = \int_{\cal F_t} \Psi_1(t,z) \overline {\Psi_2} (t,z) \lambda_{D^t}^{2-2s}(z) d^2z\;.$$
The Poincar\'e density depends smoothly\footnote{In fact, even real analytically. This follows from the fact that the solutions $w^t$ of the Beltrami equation may be chosen to depend real-analytically on $t$, which follows from the Ahlfors-Weil section (see, e.g., \cite{TCAToTS}).} on $t$. Also, the fundamental domain depends on $t$. This we can get rid of by using Lemma \ref{scalar}, with the help of which we rewrite
$$ \mathrm{WP}_{\cal R} (\Psi_1,\Psi_2)(t) = \int_{D^t} \Psi_1(t,z) \overline{h_2}(z) \lambda_{D^t}^{2-2s}(z) d^2z\;.$$
And finally, we use the solution $w^t$ to turn the integral into an integral over $\bb D$. All ingredients are now smooth in $t$. $\Box$ \\
%%%%%%%%%%%%%%%%%%%%%%%%%%%%%%%%%%%%%%%%%%%%%%%%%%%%%%%%%%%%%%55
\subsection{Families of Bases of Automorphic Forms}\label{PolAB}
A well-known theorem of Grauert says that any holomorphic vector bundle over a contractible
Stein manifold is trivial, which is equivalent to saying that a global basis exists. Since
the given basis on the fibre over $t_0 \in T(G)$ is a linear
combination of this global basis, we get the affirmative answer to
the second question, again as an existence result.
\begin{thm}
   Let $\{\psi_i\}$ be any basis of a fibre of $\cal V(\cal R)$.
   Then there are global sections $\Psi_i$ extending the $\psi_i$
   and inducing a basis on every fibre.
\end{thm}
We are interested in extending the sections of a basis on a fibre by the previous method and see how close this construction comes to a global basis. So given a basis $\{\psi_1,\ldots,\psi_N\}$, we know that any
$\psi_i$ can be extended to a global section of $\cal V(\cal R)$ by
Theorem \ref{ExtendSections}, and that this extension yields a local
basis in some neighborhood $U \ni t_0$ in $T(G)$. But `how big' is the set to which the basis extends?
\begin{thm}\label{BasisAE}
    Let $\{\psi_i\}$ be a basis of the fibre over $t_0 \in T(G)$, and $\{\Psi_i\}$ be any global
    extensions. Then they induce a basis on almost every fibre. More
    precisely, the set
    $$ N(\Psi):=\left\{t\in T(G): \{\Psi_i|_{\pi^{-1}_{\cal R}(t)}\}
    \:\textrm{is not a basis}\right\}\;,$$
    is a proper analytic submanifold of $T(G)$.
\end{thm} 
\emph{Proof.} The proof of this fact follows from a more general consideration: If we have N functions defined on an open set $\Omega \subset \bb R^N$, each of class at least $C^N(\Omega)$, we determine thier linear independence by looking at the zeroes of the Wronskian,
\begin{eqnarray*}
  w(f_1, \ldots, f_N)(x):= \mathrm{det} \left(\begin{array}{cccc} f_1(x) & f_2(x) & \ldots & f_N(x) \\ f'_1(x) & f'_2(x) & \ldots & f'_N(x) \\ \vdots & & \ddots &\\
                f^{(N)}_1(x) & f^{(N)}_2(x) & \ldots & f^{(N)}_N(x) \end{array} \right)\;.
\end{eqnarray*}
It is well-known that the functions are linearly independent in a neighborhood of any point where the Wronskian doesn't vanish. The same holds true for complex functions and linear independence over $\bb C$. \\

Now consider the functions $\Psi_1, \ldots, \Psi_N$. The induced functions $\psi^t_i(z):=\Psi_i(t,z)$ on any fibre $D^t$ are linearly independent iff they are linearly independent on some small neighborhood $U^t$ in the fibre. On the other hand, we know already that the $\Psi_i(t,0)$ are linearly independent on some open neighborhood $V\subset T(G)$ containing $t$\footnote{Recall that $0 \in D^t$ for all $t \in T(G)$ by the equi-H\"older continuity of families of $K$-quasiconformal mappings for any $K<\infty$.}. Because the Wronskian $w(\Psi_1(t,0), \ldots,\Psi_N(t,0))$ is a holomorphic function in $t$ and non-zero on $V$, the vanishing locus $N(w)$ of $w(\Psi_1(t,0), \ldots,\Psi_N(t,0))$ is a proper analytic subset of $T(G)$. This again implies that the functions, considered as functions of $z$, are linearly independent on some neighborhood $U^t\subset D^t$ for all $t \in T(G) \backslash N(w)$. $\Box$\\

It would be very interesting to find out in which sense the sets $N(\Psi)$ depend on the chosen preimages $h_i$ of $\psi_i$ on the 0-fibre. Recall that although all choices produce the same basis on the 0-fibre, this will not be the case for other fibres, i.e., the choice of preimages does affect the outcome. We will explain the occurence of the $N(\Psi)$ at the end of Section \ref{DontExtSect}.\\

Since we obtained a \emph{proper} analytic submanifold $N(\Psi)$, we can repeat the process by choosing a different extension $N(\Psi')$ of the basis. Now the set on which both extensions fail to be a basis is then given by $N(\Psi) \cap N(\Psi')$, which we can assume to be a proper analytic subset of both $N(\Psi)$ and $N(\Psi')$. By repeating this process finitely often, we have covered all of Teichm\"uller space with constructed bases.
%%%%%%%%%%%%%%%%%%%%%%%%%%%%%%%%%%%%%%%%%%%%%%%%%%%%%%%%%%%%%%%%%%%%%%%%%%%%%%%%%%%%%%%%%%%%5
\subsection{Other Methods of Extension}\label{OtherMOExtension}
We want to make a brief comment on two other ways of extending sections. The first method is very beautiful and relies on an integral operator $\cal L_C^s$ defined by Bers in \cite{An-sIEwAtQM}, while the second method follows from the existence of a canonical connection, the Chern connection, on holomorphic vector bundles over K\"ahler manifolds.\\

To describe the first method, let $C$ be a Jordan curve, $D_1$ and $D_2$ the two components of $\cinf \backslash C$
each equipped with the hyperbolic metric
$\lambda^2_{D_j}(z)|dz|^2$. Given a function $\psi :D_1 \rightarrow
\bb C$ and $s \in \bb R$, the operator $\cal L_C^s$ applied to $\psi$ formally produces a function $\phi: D_2
\rightarrow \bb C$ by 
$$ \phi (z)  = \left(\cal L^{s}_C \psi\right)(z) :=
\int_{D_1}\frac{\lambda_{D_1}^{2-2s}(w)}{(w-z)^{2s}}\ol{\psi(w)}
d^2w\;.$$ 
The following theorem is a synthesis of parts of the
Theorems 1, 2 and 3 of \cite{An-sIEwAtQM}.
\begin{thm}\label{trewq}
   Let $G$ be a quasi-Fuchsian group, $C$ its fixed curve and $2\leq s \in \bb N$. Then
   $$\cal L^{s}_C: B_s(D_1,G) \rightarrow B_s(D_2,G)\;$$
   is a continuous anti-linear isomorphism.
\end{thm}
With the help of this operator one can extend a section $\psi \in A_{\rho_s}^1(\bb D, G)$ to a section of $\cal V(\cal R)$ as follows:
\begin{itemize}
   \item{Since $A_{\rho_s}^1(\bb D, G) \cong B_{\rho_s}(\bb D,G)$ for $G$ of Fuchsian of first kind, we apply $\cal              L^s_{\partial \bb D}$ to $\psi$ and denote the result my $\psi^c \in B_{\rho_s}(\bb D^c,G)$.}
   \item{Because the restriction of the quasiconformal homeomorphism $w^\mu$ to $\bb D^c$ is holomorphic, we can transport this $\psi^c$ to any fibre via pull-back:
   $$  \Psi^c(z,t):= (w^{\mu(t)}|_{\bb D^c}^{-1})^*_s \psi(z)\;.$$
   This $\Psi^c$ is a function on the `complement' to Bers' fiber space, i.e., the space 
   $$\cal F^c(G):= \{(t,z) \in T(G) \times \bb C \mid z \in (D^t)^c\}\;,$$
   hence $\Psi^c$ is a holomorphic section of the following bundle,
      $$ \pi^c: \cal V^c(\cal R):=\bigcup_{t\in T(G)} B_{\rho_s^t}((D^t)^c, G^t)\;.$$}
\item{We obtain a section of the bundle $\cal V(\cal R)$ by applying $\cal L^s_{\partial D^t}$ fibrewise, i.e.,
      $$ \Psi(z,t) = \left(\cal L^s_{\partial D^t}\circ (w^{\mu(t)}|_{\bb D^c}^{-1})^*_s\right) \psi(z)\;.$$}
\end{itemize}
The section so obtained, however, is only real-analytic because the kernel of the integral operator $\cal L^s_{C}$ depends on the Poincar\'e density, which only varies real-analytically. Moreover, this construction only works for \emph{integer} values of $s$. The proof of Theorem \ref{trewq} given by Bers can be generalized to \emph{half-integers}, but not to any $s \in \bb Z[(2g-2)^{-1}]$ or yet more generally, to any real $s>1$.\\

Let us now describe the other method of extending sections we mentioned above. The good news is that this method will yield \emph{holomorphic sections} of the bundle $\cal V(\cal R)$; the drawback is that they are not easy to compute. The following general theorem is needed for their construction and can be found, e.g., in \cite{DGoCVB}.
\begin{thm} \label{ChernConnection}
   Let $\pi: \cal V \rightarrow B$ be a holomorphic vector bundle of finite rank over a complex manifold $B$ equipped with a hermitean form inducing a Hilbert space structure on all fibres. Then there exists a unique connection $\nabla^{\mathrm{Ch}}$ on $\cal V$, called the Chern connection, compatible with the complex and the hermitean structure.
\end{thm}
This connection is analogous to the Levi-Civita connection in Riemannian geometry. 
Let us write the compatibility conditions more explicitly. First of all, a connection $\nabla$ on a complex vector bundle $\cal V\rightarrow M$ is a map
$$\nabla: \Gamma(\cal V) \rightarrow \Gamma(\cal V \otimes \Omega^1_{\bb C}(M))\;,$$
where $\Omega^1_{\bb C}(M)$ is the complexified cotangent bundle. This bundle allows a splitting $\Omega^1_{\bb C}(M) = \Omega^{(1,0)}(M) \oplus \Omega^{(0,1)}(M)$ into holomorphic and anti-holomorphic forms, and so the connection splits into two operators $\nabla=\nabla_1 + \nabla_2$. Moreover, any holomorphic bundle carries a natural connection $\bar \partial$ induced from the complex structure. A connection $\nabla$ is said to be \emph{compatible with the complex structure} iff $\nabla_2 = \bar \partial$.\\

Let $h$ now be a hermitean structure on $\cal V \rightarrow M$. A connection $\nabla$ on $\cal V$ is called \emph{compatible with $h$} iff it satisfies
$$ d h(\Psi_1,\Psi_2) = h(\nabla\Psi_1,\Psi_2) + h(\Psi_1, \nabla\Psi_2)\;.$$ \\
Observe that Theorem \ref{ChernConnection} applies to the bundle $\pi: \cal V(\cal R) \rightarrow \cal T_B(G)$, i.e., we know $\cal V(\cal R)$ has the Chern connection. Now given a point $\psi \in A^2_{\rho_s}(\bb D,G)$, we define its extension to any other fibre by \emph{parallel transport} with respect to the Chern connection. The extension does depend on the curves used to connect the other fibres with the 0-fibre in general, but A consistent choice of such curves nevertheless exists because of the contractibility of $T(G)$.\\
% and because of the compatibility conditions satisfied by the connection, the section we obtain will be automatically %\emph{holomorphic}.\\

In the particular case of the torus, this point of view has very interesting and deep connections to Theta functions and modular forms (see, e.g., \cite{GQoCST}). In fact, the generalization of these matters to surfaces of higher genus is an ongoing project of ours.
\addtocontents{toc}{\vspace{.5cm}}
\section*{\bf The Frontier}\label{Frontier}
\addtocontents{toc}{\vspace{.3cm}}
\section{Riemann Surfaces and their Limits}
We start to explore boundary phenomena by geometric means in this section. At first we describe what happens geometrically at the boundary from different points of view and then study the boundary points as geometric objects, namely \emph{noded Riemann surfaces}. In the last subsection we then extend the formal definition of Teichm\"uller spaces to the so-called \emph{deformation spaces} of Riemann surfaces.
%%%%%%%%%%%%%%%%%%%%%%%%%%%%%%%%%%%%%%%%%%%%%%%%%%%%%%%%%%%%%%%%%%%%%
\subsection{Degenerations of Riemann Surfaces}\label{DegenerofRS}
\begin{figure}[h]
    \includegraphics [scale=.7]{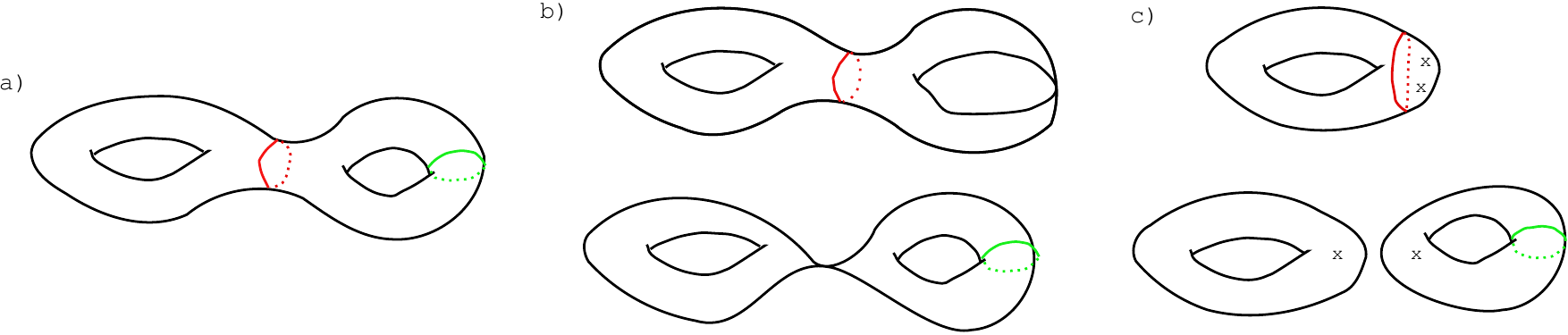}
    \caption{\label{DegRS} a) A genus two surface equipped with a homologically trivial (red) and a homologically  non-trivial (green) curve. b) The topological picture of the surfaces obtained by pinching the green resp.~ the red curve. c) The complex analytic viewpoint of the limit surfaces.}
\end{figure}

If one tries think in terms of hyperbolic objects\footnote{By which we meant that we cannot do wild topological deformations, e.g., shrink arbitrary closed sets etc., but that the surface should be viewed as carrying a hyperbolic metric throughout the deformation.} the only degenerations one comes up with after some throught are the ones displayed in Figure \ref{DegRS}.\\ 

These degenerations\footnote{Recall that the length of a simple closed geodesic in the free homotopy class of a loop is one of the Fenchel-Nielsen coordinates, so changing its length is a displacement in Teichm\"uller space.} are called \emph{pinching simple closed geodesics}. The figure shows the two principally different cases, namely pinching a homologically trivial, and a homologically non-trivial, geodesic. That this intuition is indeed correct is the content of a corollary to a rather general result by Mumford \cite{ARoMaCT}. The corollary is known as Mumford's Compactness Theorem\footnote{Mumford stated the result for $\Sigma$ being of finite type without punctures. The more general result given here is due to Bers \cite{ARoMuCT}.}.
\begin{thm}[Mumford's Compactness Theorem]
    Let $\Sigma$ be a Riemann surface of finite type. Then the subset of $\cal M(\Sigma)$, consisting of surfaces such that the hyperbolic length of all closed geodesics on them is bounded from below by a constant $C>0$, is compact.
\end{thm}
Of course, this result deals with the moduli space instead of Teichm\"uller space, but from many points of view, Teichm\"uller space is only an intermediate object and so the guiding ideas should be taken from the point of view of the moduli space. Now, the quotient map from dividing out the mapping class group is not a proper map, since the orbit of a point in Teichm\"uller space is infinite. Hence, if we take the viewpoint from moduli space, we are not necessarily looking for a \emph{compactification} but rather for a \emph{partial completion} of Teichm\"uller space, regarding non-compactness of the orbits of $\mcg(\Sigma)$ as irrelevant. The degenerations we have to cope with at the level of Teichm\"uller spaces hence are precisely the pinching of geodesics.\\

Let us, just for convenience, rephrase Mumford's Compactness Theorem at the level of Teichm\"uller spaces.
\begin{cor}
   Let $[\Sigma, f_i, \Sigma_i]$ be a sequence in $\mathrm{Teich}(\Sigma)$ such that the hyperbolic length of the closed geodesics is bounded from below by $C>0$. Then there exists a subsequence $[\Sigma, f_{i_j}, \Sigma_{i_j}]$ and a sequence $\gamma_j \in \mcg(\Sigma)$ of elements of the mapping class group such that $ [\Sigma, \gamma_j \circ f_{i_j}, \gamma_j^{-1}(\Sigma_{i_j})]$ converges.
\end{cor}
Let us think of the Beltrami model of Teichm\"uller space $\cal T_B(\Sigma)$ for a moment. We would intuitively think that a sequence of Beltrami differentials $\mu_i$ such that $\pi_T(\mu_i) \in \cal T_B(\Sigma)$ corresponds to a degenerating sequence of Riemann surfaces has divergent maximal dilatation, i.e., $K(\mu_i) \rightarrow \infty$ or equivalently $\|\mu_i\|_\infty \rightarrow 1$. This is indeed true and follows from a lemma by Wolpert \cite{TLSaMfCRS}, the proof of which goes by a basic estimate in terms of moduli of rectangles.
\begin{lemma}
   Let $\gamma$ be a closed curve on $\Sigma$ and $[\Sigma,f,\Sigma']$ a point in $\mathrm{Teich}(\Sigma)$. Let $l_h(\gamma)$ denote the hyperbolic          length\footnote{In case there is none, e.g., when $\gamma$ wraps around a puncture, set $l_h(\gamma) = 0$.} of the unique closed geodesic freely homotopic to $\gamma$. Then
   $$ K(f)^{-1} l_h(\gamma) \leq l_h(f(\gamma)) \leq K(f) l_h(\gamma)\;,$$
\end{lemma}
A conclusion of this lemma is that a $K$-quasiconformal map between Riemann surfaces is a $K$-quasi-isometry with respect to the hyperbolic metrics.\\

This lemma also tells us something about the mapping class group
action: The orbit of any point in $\mathrm{Teich}(\Sigma)$ under the
$\mcg(\Sigma)$-action also has accumulation points of infinite
maximal dilatation. This is so because the mapping class group is
generated by Dehn twists about the simple closed geodesics, e.g., let
$\gamma$ be such a geodesic and denote by $\Theta_\gamma \in \mcg(\Sigma)$ the
twisting around $\gamma$. Denote by $´\gamma'$ another simple closed
geodesic which intersects $\gamma$. Then
$$l(\Theta^n_\gamma(\gamma')) \sim  l(\gamma') + n l(\gamma) \rightarrow \infty\;,$$
and so necessarily, by abuse of notation, $K(\Theta_\gamma)\rightarrow \infty$.\\

In algebraic geometry there is a standard compactification of the moduli space of curves known as the Deligne-Mumford compactification \cite{DMComp}, which we will denote by $\overline{\cal M}_{\mathrm{DM}}(\Sigma)$. It is an object of intense study. The beauty of its construction is that it works in far more generality than for Riemann surfaces, i.e., for curves over any (possibly non-algebraically closed) field of arbitrary characteristic. The price one has to pay is that one has to work in a very technically demanding category, that of (Deligne-Mumford) stacks, and this can obscure a lot of the geometry. We will not go into more detail on these matters here, but just mention that the description of the compactification is similar to the one we give later by adding noded Riemann surfaces as boundary points, which are known as \emph{stable algebraic curves}.
%%%%%%%%%%%%%%%%%%%%%%%%%%%%%%%%%%%%%%%%%%%%%%%%%%%%%%%%%%%%%%%%%%%%%%%%%%%%%%%%%%%%%%%%%%%%%%%%%%%%%%%%%%%%%%%%%%%%%%%%%%%%5
\subsection{Noded Riemann Surfaces}\label{nodedRS}
Noded Riemann surfaces are the geometric objects one obtains as
limits of degenerating Riemann surfaces. Let's give a more formal
definition which uses the language of \emph{complex spaces} (see,
e.g., \cite{CAS}).
\begin{dfn}\label{DefNodedRS}
    A \underline{noded Riemann surface} $\Sigma'$ is a connected complex space where any point $x$ has a neighborhood $U_x$ such that
    $$ (\Sigma', \cal O(\Sigma'))_{|U_x} \cong (\bb D, \cal O(\bb D))  \qquad \mathrm{or} \qquad (\Sigma', \cal O(\Sigma'))_{|U_x} \cong (\bb D^2, \cal O(\bb D^2))/\langle zw \rangle\;,$$ where $z,w$ are the coordinates on $\bb D^2$ and $(0,0)$ corresponds to $x$. Moreover, the surface obtained by removing the \underline{nodes} $N(\Sigma')$, i.e., the points having the latter type of neighborhood, must be a Riemann surface of finite analytic type.
\end{dfn}
By definition, the set of nodes $N(\Sigma')$ are the only singularities of the complex space $\Sigma'$. The smooth part therefore is given by
$$\Sigma' \backslash N(\Sigma') = \Sigma_1 \cup \ldots \cup \Sigma_p \;.$$
The surfaces $\Sigma_i$ are called the \emph{parts} of $\Sigma'$; denote their type by $(g_i, n_i)$. A noded Riemann surface is called \emph{stable} iff all parts are stable (which means that $\chi(\Sigma_i) = 2-2g_i - n_i < 0$, or equivalently the automorphism group of the surface is finite) and it is called \emph{maximally degenerate} iff all parts are thrice punctured spheres.\\

\begin{figure}[t]
    \label{NodedRS}
    \includegraphics [scale=.7]{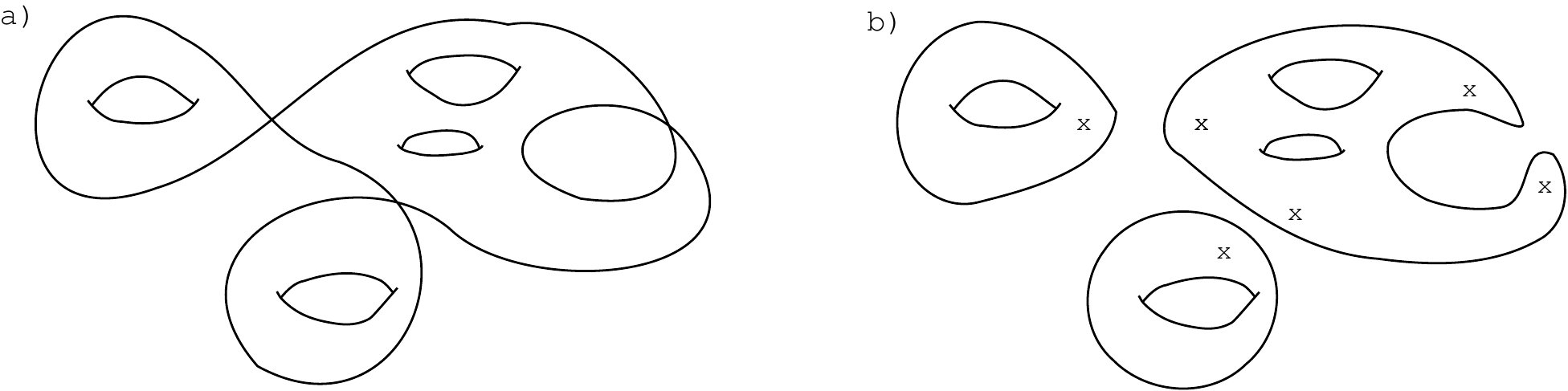}
    \caption{a) displays a noded Riemann surface of genus 5 while b) shows its regular part which is the disjoint union of two once-punctured tori and a four times punctured Riemann surface of genus two.}
\end{figure}
From the point of view of complex structures, the decomposition into parts is the right way to view these objects. In particular, the Teichm\"uller space of a noded Riemann surface $\Sigma'$ is identified with the product of the Teichm\"uller spaces of its parts,
$$ \mathrm{Teich}(\Sigma') = \mathrm{Teich}(\Sigma_1) \times \ldots \times \mathrm{Teich}(\Sigma_p)\;.$$
Given any simple closed curve $\gamma$ on a Riemann surface $\Sigma$ there is a natural way to associate a noded Riemann surface $\Sigma_\gamma$ to it by \emph{pinching} it, i.e., $\Sigma_\gamma = \Sigma / \gamma$ in the sense of a topological quotient. The hyperbolic, resp.~ complex, construction of this procedure is explained in Section \ref{Strat}. \\

There is just one little danger in the pinching picture: After the process of pinching, the surface seems to touch itself `tangentially' at the node. But if we look at the description of the model neighborhood of a node in Definition \ref{DefNodedRS}, the two parts of the surface actually intersect \emph{transversely} at the node.
%=========================================================================================================
\subsection{The Strong Deformation Space of a Riemann Surface}
Recall the Definition of Teichm\"uller space as equivalence classes of maps between Riemann surfaces (Def. \ref{DefTeich}). Bers enlarged this definition to also include deformations of the surface pinching some curves or opening some nodes; see \cite{SoDRS} or the review article \cite{FDTSaG}.
\begin{dfn}\label{DefSigmaofSurface}
   Let $\Sigma$ and $\Sigma'$ be noded Riemann surfaces.
   \begin{enumerate}
   \item{A \underline{strong deformation} of $\Sigma$ onto $\Sigma'$ is a continuous surjective map $f:\Sigma \rightarrow \Sigma'$, also denoted $(\Sigma,f,\Sigma')$, such that
   \begin{itemize}
      \item{the image of a node on $\Sigma$ is a node on $\Sigma'$,}
      \item{the preimage of a node on $\Sigma'$ is either a node or a closed Jordan curve on $\Sigma$ omitting all nodes and}
      \item{$f$ restricted to $\Sigma \backslash f^{-1}(N(\Sigma'))$ is an orientation-preserving homeomorphism.}
      \end{itemize}}
   \item{The \underline{strong deformation space} of $\Sigma'$, denoted $\mathrm{Def}(\Sigma')$, is the set equivalence classes of all deformations onto $\Sigma'$, under the equivalence relation $ (\Sigma_1,f_1,\Sigma') \sim (\Sigma_2,g,\Sigma')$ iff there exist homeomorphisms $\psi: \Sigma_1 \rightarrow \Sigma_2$ and $\phi: \Sigma' \rightarrow \Sigma'$, where $\psi$ is homotopic to a biholomorphism and $\phi$ is homotopic to the identity, on $\Sigma'$
   such that the following natural diagram commutes,
   \begin{diagram}
         \Sigma_1 & \rTo^{f_1} & \Sigma'\\
         \dTo_\psi & &\dTo_\phi\\
         \Sigma_2 & \rTo_{f_2} & \Sigma'\;.
   \end{diagram}}
   \end{enumerate}
\end{dfn}
Observe here that in contrast to the definition of $\mathrm{Teich}(\Sigma)$ we have switched the source and target: An element of $\mathrm{Teich}(\Sigma)$ is an equivalence class of a map $f$ \emph{from} $\Sigma$ to some other Riemann surface while an element of $\mathrm{Def}(\Sigma)$ is an equivalence class of a map from some other surface \emph{onto} $\Sigma$. The reason for this is simple: The definition of $\mathrm{Teich}(\Sigma)$ doesn't really care about the order, since all maps are homeomorphisms and therefore invertible\footnote{Which is also the reason for the shorter formulation of the equivalence relation in Def. \ref{DefTeich} though it is the same idea.}, so one could give an equivalent definition of $\mathrm{Teich}(\Sigma)$ in terms of maps onto $\Sigma$. This is \emph{not} the case for $\mathrm{Def}(\Sigma)$: The definition given here allows nodes of $\Sigma$ to be opened, and hence for a surface without nodes, $\mathrm{Def}(\Sigma) \cong \mathrm{Teich}(\Sigma)$, whereas if we had defined it the other way around, $\mathrm{Def}(\Sigma)$ would not allow the opening of nodes but rather the pinching of curves. In that case $\mathrm{Def}(\Sigma)$ would correspond to $\mathrm{Teich}(\Sigma)$ plus the Teichm\"uller spaces of all possible pinchings of curves on $\Sigma$.\\

At first it could be confusing to call this space the deformation space of $\Sigma$, since we introduced the deformation space $\mathrm{Def}(G)$ of a Fuchsian group in Section \ref{Bersemb} but actually the analogy to the definitions of $\cal T_S(G)$ and $\mathrm{Def}(G)$ and their mutual relationship are quite deep.\\

The deformation space is not only a set, but a topological space. A set $A \subset \mathrm{Def}(\Sigma)$ is open in this topology iff for any point $p = [\Sigma',f_p,\Sigma] \in A$ there are finitely many curves $C_1,\ldots , C_L$ on the parts of $\Sigma'$ such that for all deformations $ [\Sigma'',f_h,\Sigma']$ satisfying
$$|l_h(C_i) - l_h(f^{-1}_h (C_i))| < \epsilon \qquad \mathrm{and} \qquad |l_h(f_h^{-1}(n)| < \epsilon \qquad \forall 1\leq  i \leq L\;, n \in N(\Sigma')\;,$$
the deformation $[\Sigma'', f_p \circ f_h, \Sigma] \in A$. One easily sees that this topology agrees with the topology of Teichm\"uller space, e.g., by considering geodesic length functions as coordinates, on the subset $\mathrm{Def}(\Sigma_0) \subset \mathrm{Def} (\Sigma)$, where $\Sigma_0$ is the surface obtained from \emph{opening the nodes} (the procedure of opening nodes will be explained in Section \ref{Strat}) of $\Sigma$. 
%%%%%%%%%%%%%%%%%%%%%%%%%%%%%%%%%%%%%%%%%%%%%%%%%%%%%%%%%%%%%%%%%%%%%
\addtocontents{toc}{\vspace{.3cm}}
\section{Completions of Teichm\"uller Spaces}
We start this section by building a topological model space for a completion of Teichm\"uller space, taking into account the pinching of simple closed curves. We then proceed by describing briefly the so-called Weil-Petersson metric and the completion of $\mathrm{Teich}(\Sigma)$ with respect to that metric using Masur's estimates, from which follows that the two described spaces are homeomorphic.\\

After that we turn our attention to the Bers completion of Teichm\"uller space, which is a compactification, and study the relation between boundary points of this space and the associated geometric objects. This is done in order to single out a subset of the Bers compactification, called the \emph{Bers augmented Teichm\"uller space} $\hat T^B(G)$ where $\Sigma = \bb D/G$. This subset of the Bers compactification can be equipped with a different topology, and we do so as we proceed, thereby obtaining $\hat T(G)$, the \emph{augmented Teichm\"uller space}, which is a partial completion matching all our requirements of such an object. In particular, $\mcg(\Sigma)$ acts by homeomorphisms on $\hat T(\Sigma)$ and the quotient space of this action is homeomorphic to $\overline{\cal M}_{\mathrm{DM}}(\Sigma)$, the Deligne-Mumford compactification of the moduli space (see \cite{DFoRS}). \\

Moreover, the augmented Teichm\"uller space admits a natural extension of the analytic objects we considered over Teichm\"uller space in part two of this paper, namely the fibre space $\cal F(\Sigma)$ and the vector bundle $\cal V(\cal R)$. This will become crucial in Section \ref{ExtoATS}, where we construct extensions of the previously constructed sections of bundles over Teichm\"uller space to the completion.\\

Let us remark that there are several points of view of Teichm\"uller spaces not at all considered in this work, each of which lead to some different concept of completion, e.g., the \emph{Teichm\"uller} and the \emph{Thurston} completion, both of which are neither mutually homeomorphic nor homeomorphic to the completions considered here (see \cite{TAGoTS}, or for even more details see \cite{TBoTS}), nor any other metric completion of which there exist many
%%%%%%%%%%%%%%%%%%%%%%%%%%%%%%%%%%%%%%%%%%%%%%%%%%%%%%%%%%%%%%%%%%%%%%%%%%%%%%%%%%%%%%%%%%%%%%%%%%%%%%%%%%
\subsection{A Topological Stratification of Teichm\"uller Space and Horocyclic Coordinates} \label{Strat}
Guided by the idea that pinched curves are the source of non-completeness, let us first build a candidate topological space. Let us first give sets of curves which can be pinched simultaneously a name.
\begin{dfn}
    Let $\cal C$ be a set of oriented simple closed curves on $\Sigma$. $\cal C$ is called \underline{homotopically independent} iff the curves are disjoint, none of the curves is a loop around a puncture, and the free homotopy classes of all the curves in $\cal C$ are all non-trivial and distinct from one another. Such a set is called \underline{maximal} iff no curve can be adjoined to the set such that the resulting set is still homotopically independent.\end{dfn}
For example, any pair of pants decomposition of a Riemann surface of finite type is given by a maximal homotopically independent set of curves. The requirement on the curves not to bound a disc or a punctured disc implies that after pinching these curves, the resulting parts will still carry hyperbolic metrics, i.e., we don't obtain once or twice punctured spheres. The orientation of the curves required in the definition above is harmless and neglected most of the time; it is only needed when we describe the Earle-Marden coordinates. But since an orientation represents no difficulty whatsoever we invoke it in the definition in general.\\

From now on, let the symbol $\cal L$ denote a maximal homotopically independent set of curves on $\Sigma$. Define the space
$$ \overline {\mathrm{Teich}}(\Sigma) := \bigcup_{\cal C \in 2^{\cal L}} \mathrm{Teich}(\Sigma_{\cal C})\;,$$
where $\Sigma_C$ is the noded surface obtained from pinching the curves in $\cal C \subset \cal L$. The space is equipped with the \emph{finest} topology for which all the inclusions $\mathrm{Teich}(\Sigma_{\cal C}) \hookrightarrow \overline{ \mathrm{Teich}}(\Sigma)$ are continuous. Recall that $$\Teich(\Sigma_{\cal C}) \cong \Teich(\Sigma^1) \times \ldots \Teich(\Sigma^p)\;,$$
where $p$ is the number of parts $\Sigma^i$ resulting from the pinching as described in Section \ref{nodedRS}.\\

\begin{figure}[t]
        \includegraphics [scale=.8]{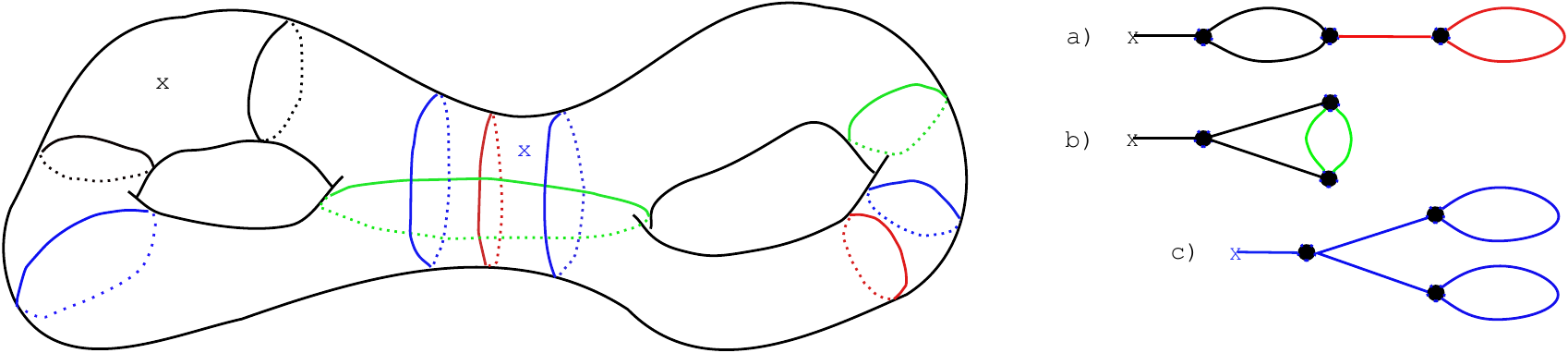}
    \caption{\label{CurvesGraph}Three different maximal sets of curves on a once-punctured genus two surface and the resulting graphs: a) black + red, b) black + green - both with the black puncture and c) blue with the blue puncture.}
\end{figure}

We will now explain the \emph{Earle-Marden} coordinate construction on the space $\overline{\Teich}(\Sigma)$. The space $\Teich(\Sigma_{\cal L})$ is a point, since by pinching a maximal set of curves, the parts $\Sigma_i$, $i= 1, \ldots, 2g-2 + n$, are all thrice punctured spheres, which have no moduli. Nevertheless, depending on the chosen system one gets possibly different graphs, as illustrated in Figure \ref{CurvesGraph}.
Let $z^i_1, \ldots, z^i_{a_i}$, $a_i \in \{1,\ldots,3\}$, denote local coordinate systems around the punctures of $\Sigma_i$ resulting from pinching a curve of $\Sigma$, provided $\Sigma_i$ is to the left of that curve, and equip the punctures on $\Sigma_i$ that are to the right of a pinched curve by the coordinate $w^j_b$, where the labels $j,b$ are taken to be the same as the ones of the coordinate $z^j_b$ to the left of that curve. Recall that a neighborhood $\tilde U_i$ of the noded surface $\Sigma_{\cal L}$ around a node is described by
$$ (z^i_j,w^i_j) \in U_i \subset \bb C^2\;, \qquad (\Sigma_{\cal S}, \cal O(\Sigma_{\cal S}))_{|\tilde U_i} \cong (U_i, \cal O(U_i))/\langle z^i_j w^i_j = 0\rangle \;.$$
Now apply the \emph{opening of nodes} construction to all pairs $(z^i_j, w^i_j)$, $i=1, \ldots, 2g-2+n$ and $j=1, \ldots, a_i$. More precisely, this means the following: Let $(z,w)$ denote such a pair and assume without loss of generality that they are defined on punctured discs. Replace the two punctured discs by
$$ D_t:= \left(\bb D^*_z \amalg \bb D^*_w) \right / \sim_t \qquad z \sim_t w \Leftrightarrow z \cdot w = t\;,$$
for any $t \in \bb D$. For $t\neq 0$ this is an annulus replacing the node with hyperbolic length of the closed curve $\phi: [0,1] \rightarrow D_t$, $\phi(s) = \sqrt{|t|} \mathrm{exp}(2 \pi i s)$, given\footnote{The punctured disc inherits a hyperbolic metric from the surface, which of course is \emph{not} complete. All hyperbolic metrics on a punctured disc are obtained by embedding the punctured disc into the complete punctured disc as punctured subdiscs $\bb D^*_{c^2}$, so repacing the discs $\bb D_{z}^*, \bb D_{w}^*$ by $\bb D_{c^2,z}^*, \bb D_{c^2,w}^*$ they become equipped with hyperbolic metrics from the inclusion $\bb D_{c^2} \hookrightarrow \bb D$. The length of the curve $\phi(s)$ in this metric is $2 \pi/\mathrm{log}\sqrt |t|$ and not the same as claimed above. This is because the opening construction is not so simple from the hyperbolic standpoint; the length must rather be computed in a different metric obtained by a procedure called \emph{grafting}. For a very detailed account on these matters we refer to \cite{THMatGotUC}.} by $2 \pi^2 / \mathrm{log}|t|$. \\

If we denote the surface obtained by opening the nodes with parameter $\vec t:= (t_1^1, \ldots, t^p_{a_p})$ by $\Sigma^{\vec t}$ and the total space by $\cal V(\Sigma_{\cal C})$, we get the following family of Riemann surfaces over the punctured polydisc,
$$\pi_{\vec t}: \cal V(\Sigma_{\cal C}) := \bigcup_{t^i_j \subset \bb D^*}  \Sigma^{\vec t} \rightarrow (\bb D^*)^{3g-3+n}\;,$$
for which, by the universal property of Teichm\"uller space, there exists a map $F: (\bb D^*)^{3g-3+n} \rightarrow \mathrm{Teich}(\Sigma)$, which identifies this family with the pull-back of the family $\cal V(\Sigma_{\cal C})$. 

Since the length of the closed geodesic is a monotone function of $|t|$ and the argument of $t$ of course corresponds to twisting around the geodesic, it seems reasonable to believe that all surfaces so obtained are not Teichm\"uller equivalent, because their Fenchel-Nielsen coordinates differ. There is a pitfall here: By pinching a curve one loses parts of the marking and so one obtains neighborhoods in Teichm\"uller spaces factored by subgroups of the mapping class group. Since this can switch markings, the map is not neccesarily a covering map, i.e., doesn't provide local coordinates. There has been some controversy in the literature, e.g., compare \cite{GCCfTS} with \cite{E-MCoTS:AC}. For a detailed account all the technicalities of this intriguing problem, see the very recent \cite{ATSaO}.\\

The difficulties described above can be dealt with by choosing \emph{special coordinates} around the punctures. Kra has written a detailed account on a certain choice of these, called \emph{horocyclic coordinates}, in \cite{HCfRSaMS}\footnote{This explicit and long paper deals with the construction and properties of these coordinates, which had been constructed and used by other authors, e.g., by Earle, Marden and Maskit, before, yet details were not available.}. The already mentioned paper (\cite{ATSaO}, p.~43-46) provides other coordinates and proves a more general theorem, namely the following.
\begin{thm}\label{EMCexist}
   Any point $[\Sigma, f, \Sigma']$ lies in the range of an injective map $F$ coming from the opening of $2g-2+n$ nodes for some choice of local coordinates $z_j^i$ and $w^i_j$.
\end{thm}
This means that we can think of any degenerating family of Riemann surfaces as given by the construction above as long as we don't care about the precise coordinates around the punctures. Moreover, the theorem is also true for any point $[\Sigma',f,\Sigma] \in \mathrm{Def}(\Sigma)$, but we have yet to explain the precise relation between $\overline{\mathrm{Teich}}(\Sigma)$ and $\mathrm{Def}(\Sigma)$, which will be clear at the end of Section \ref{AugmentedTS}.
%%%%%%%%%%%%%%%%%%%%%%%%%%%%%%%%%%%%%%%%%%%%%%%%%%%%%%%%%%%%%%%%%%%%%%%%%%%%%%%%%%%%%%%%%%%%%%%%%%%%%%%%%%%%%%%%%%%%%%%%%%%%%
\subsection{The Weil-Petersson Metric Completion} \label{WPComp}
The seminal paper describing the metric completion of Teichm\"uller space with respect to a natural metric, the Weil-Petersson metric, is \cite{EotWPMttBoTS} by H.~ Masur. We should start by introducing the Weil-Petersson (co)metric on Teichm\"uller space, and actually we have already done that in a different guise in Section \ref{BSoAF}. The tangent space to Teichm\"uller space at a point $[\Sigma,f,\Sigma']$ can naturally be identified with the \emph{harmonic Beltrami differentials}\footnote{Let us give some more detail here: We have two manifold models of Teichm\"uller space at hand, $\cal T_B(G)$ and $T(G)$. The latter one, since it is an open subset of a Banach space, suggests identifying $T_\phi T(G) \cong B_2(\bb D^c, G)$ for all points $\phi \in T(G) \subset B_2(\bb D^c, G)$, i.e., identifying the tangent space of Teichm\"uller space at \emph{any} point with the space of holomorphic quadratic differentials on the fixed conjugate Riemann surface of $\Sigma$, namely $\Sigma^*:= \bb D^c/G$. This seems a bit unnatural and indeed is, since the Bers embedding is base point dependent. So what is the tangent space to $\cal T_B(G)$? Let's start at the origin: Since the differential of the lifted Bers embedding is surjective we have a natural association,
$$ L^\infty_{(-1,1)} (\bb D, G) / \mathrm{Ker} D_0 \tilde \beta_3 \cong T_0\cal T_B(G)\;,$$
and there kernel turns out to be precisely the integrable holomorphic quadratic differentials,
$$ \mathrm{Ker} D_0 \tilde \beta_3 = A^1_2(\bb D,G)^\perp \subset L^\infty_{(-1,1)} (\bb D, G)\;. $$
On the other hand, if we start with an integrable quadratic differential $0 \neq \phi \in A^1_2(\bb D,G)$,  then by the transformation behaviour of the Poincar\'e density, $\lambda^{-2}_{\bb D} \bar \phi \in L^\infty_{(-1,1)}(\bb D, G)$, since $A^1_2(\bb D,G)\subset L^\infty_2(\bb D,G)$. Now this element is not in the kernel, because
$$ \langle \phi, \lambda^{-2}_{\bb D} \bar \phi \rangle^G_2 = \|\phi \|_{A_2^1(\bb D, G)}.$$
For dimensional reasons we then get $T_0\cal T_B(G) \cong \mathrm{HBelt}(\bb D, G) := \{\lambda_{\bb D}^{-2} \bar \phi, \: \phi \in A_2^1(\bb
D,G)\}$. Elements of the space $\mathrm{HBelt}(\bb D, G)$ are called the \emph{harmonic Beltrami differentials}. This identification is now more
intrinsic in the following way: The right-translation $R_\mu$ on $\cal T_B(G)$ induces an isomorphism of tangent spaces, $(R_\mu)_*T_0\cal T_B(G)
\cong T_\mu \cal T_B(G)$, and acts naturally on the identification, too, namely $(R_\mu)_* \mathrm{HBelt}(\bb D,G) = \mathrm{HBelt}(\bb D^\mu,
G^\mu)$. For more details, see, e.g., \cite{TCAToTS} or any other textbook on Teichm\"uller theory.} on $\Sigma$, and naturally dual to this space is
the space of holomorphic quadratic differntials on $\Sigma$. This space, described by automorphic forms on the covering, has a natural pairing, using which we define the Weil-Petersson metric on a tangent space,
$$ \langle t_1, t_2\rangle_{\mathrm{WP}} := \langle \phi_1, \phi_2 \rangle^{G^\mu}_2 \;, \qquad \mathrm{where}\quad t_i = \lambda_{D\mu}^{-2} \bar \phi_i \quad \mathrm{and} \quad \phi_i \in A_2^1(D^\mu, G^\mu)\;.$$
The induced metric on $\cal T_B(G)$, called the Weil-Petersson metric, is very remarkable: It is K\"ahler, has negative sectional curvature and though non-complete, it makes Teichm\"uller space geodesically convex. For more details we refer the reader to the review articles \cite{GotWPCoTS}, \cite{OtWPGotCoTS} and references therein. These accounts are quite beautiful, and we will not attempt to recreate them.\\

In \cite{NCotWPMfTS} Wolpert proved that Teichm\"uller space equipped with this metric is non-complete. Any non-complete metric space can be completed, and so one obtains the Weil-Petersson completion of Teichm\"uller space, $\overline{\mathrm{Teich}}(\Sigma)_{\mathrm{WP}}$. Masur, however, proved that this completion is homeomorphic to the one described in the previous section, i.e.,
$$ \overline{\mathrm{Teich}}_{\mathrm{WP}}(\Sigma) \cong \overline{\mathrm{Teich}}(\Sigma)\;.$$
He achieved this by providing qualitative information about the behaviour of the Weil-Petersson metric as one approaches the boundary.

In order to describe this, let $0<r<1$ be a fixed real number, $N$ the dimension of $\mathrm{Teich}(\Sigma)$ and let
$$(\vec t, \vec \tau)=(t_1,\ldots t_n, \tau_{n+1},\ldots, \tau_N) \in (\bb D^*)^n \times (\bb D - \bb D_r)^{N-n}\;,$$
be some choice of Earle-Marden coordinates on Teichm\"uller space such that $\Sigma$ is in the range (such coordinates exist by Theorem \ref{EMCexist}). We have split the coordinates into two groups, because we want to consider the pinching of the $t$'s and the interior variation of the $\tau$'s, i.e., the deformations of the latter are transversal to the degeneration. For the sake of exposition, let greek indices run form $1$ to $n$ and latin indices from $n+1$ to $N$. Expressed in these coordinates, Masur \cite{EotWPMttBoTS} proved that the Weil-Petersson metric tensor $\cal G$ satisfies
\begin{eqnarray} \label{WPasympt}
 \cal G_{\alpha i} \in O\left(|t_\alpha|^{-1}\log |t_\alpha|^{-3}\right) \;,\qquad  \cal G_{\alpha \beta} \in
  \begin{cases} O\left(|t_\alpha|^{-2} \log |t_\alpha|^{-3}\right), & \alpha = \beta\\
   O\left(|t_\alpha|^{-1} |t_\beta|^{-1} \log |t_\alpha|^{-3}\log |t_\beta|^{-3}\right), & \alpha \neq \beta\;, \end{cases}
\end{eqnarray}
for $\vec t \rightarrow 0$. Moreover, he not only proved that the components $\cal G_{ij}$ corresponding to non-pinched directions stay bounded but rather that they converge to the Weil-Petersson metric on the boundary Teichm\"uller space\footnote{The subset $\vec t =0$ of the metric completion corresponds to $\mathrm{Teich}(\Sigma_{\cal C})= \mathrm{Teich}(\Sigma_1) \times \ldots \times \mathrm{Teich(\Sigma_p)}$ where the set $\cal C$ is the set of curves pinched by setting $\vec t =0$ and so this also carries a Weil-Petersson metric.}, i.e., that the boundary strata are Teichm\"uller spaces together with their Weil-Petersson metrics. The asymptotic behaviour given in \eqref{WPasympt} has since then been strengthened considerably by Wolpert, see \cite{GotWPCoTS}.
\begin{thm} [\cite{EotWPMttBoTS}, Thm.~ 2 and Cor.~ 2]
   The Weil-Petersson metric on $\mathrm{Teich}(\Sigma)$ induces a metric $d_{\mathrm{WP}}$ on $\overline{\mathrm{Teich}}(\Sigma)$ compatible with the topology on $\overline{\mathrm{Teich}}(\Sigma)$. Moreover $d_{\mathrm{WP}}$ is mapping class group-invariant and induces a complete metric on $\overline {\cal M}(\Sigma)$ for which $\overline {\cal  M}(\Sigma)_{\mathrm{WP}} \cong \overline {\cal M}(\Sigma)_{\mathrm{DM}}$.
\end{thm}
The research field of metric completions of Teichm\"uller spaces is vast: There are many other natural metrics on Teichm\"uller space, e.g., the Teichm\"uller metric, the Bergman metric, the Arakelov metric, the McMullen metric, the Caratheodory metric, the Kobayashi metric, the Velling-Kirillov metric, the Takhtajan-Zograf metric - and the list probably goes on. A recent two paper series discussing various knowledge about some of these metrics on Teichm\"uller space is \cite{CMotMSoRS:1} and \cite{CMotMSoRS:2}, which also contain many references to previous work on this subject.
%%%%%%%%%%%%%%%%%%%%%%%%%%%%%%%%%%%%%%%%%%%%%%%%%%%%%%%%%%%%%%%%%%%%%%%%%%%%%%%%%%%%%%%%%%%%%%%%%%%%%%%%%
\subsection{The Bers Boundary of Teichm\"uller Space}\label{BersBound}
A rather straightforward idea of completing Teichm\"uller space is to take the closure of $T(G)$, the image of the Bers embedding. The purpose of this section, however, is to show that this is \emph{not} a good answer to the problem of completing Teichm\"uller space, since the space so obtained \emph{depends on the chosen base surface} $\Sigma$, and the change of basepoint maps do \emph{not all extend to homeomorphisms of the closure}. This was proved by Kerckhoff and Thurston in \cite{NotAotMGaBBoTS}, and is of course conceptually disappointing and rules out the possibility of obtaining a compactification of moduli space by $\overline T(G)/\cal{MCG}(G)$.\\

Another hint, that $\overline T(G)$ is not a good choice for a completion of Teichm\"uller space is that the dimension of the boundary has \emph{real codimension} one. This does not fit together with the idea of the last section, where we argued that the correct notion of degeneration produces a space of \emph{complex codimension} one\footnote{Recall that pinching oa curve $\gamma$ on $\Sigma$ produces two new punctures and either diminishes the genus by one or disconnects the surface into two surfaces, the genera of which add up to the genus of $\Sigma$. In either case
$$ \mathrm{dim}_{\bb C}\mathrm{Teich}(\Sigma_\gamma) = \begin{cases} 3(g-1)-3 + n + 2\\
 3g_1-3+n_1+1+3g_2-3+n_2+1  \end{cases} = 3g - 3 + n - 1 = \mathrm{dim}_{\bb C}\mathrm{Teich}(\Sigma) -1\;.$$}.\\

Nonetheless, the big advantage of the Bers compactification is that the points that are added to Teichm\"uller space \emph{have geometric meaning}: Recall that there was a 1-1 correspondence
$$c_G: \mathrm{Def}(G) \rightarrow \cal Q(G)\;,$$
and on the level of vector spaces $\cal Q(G) \cong B_2(\bb D^c,G)$, so any point of $B_2(\bb D^c,G)$ can be associated to an element of the deformation space of the Fuchsian group $G$, in particular the boundary points of $T(G)$. This section will deal with the study of the geometry of these boundary points.
\begin{dfn}
 The \underline{Bers boundary} of Teichm\"uller space, denoted $\partial T(G)$, is the topological boundary of the set
 $T(G)$ in $B_2(\bb D^c, G)$.
\end{dfn}
In other words $\partial T(G) = \mathrm{cl}\:T(G) - \mathrm{int}\:T(G)$. So we still have $\partial T(G) \subset \bb S(G) := \bb S \cap B_2(\bb D^c, G)$ since the set $\bb S$ was seen to be closed in $B_2(\bb D^c)$.
For a given $\phi \in \cal Q(G)$, let us introduce some short notation,
\begin{eqnarray}
    f^\phi:=(\mathrm{dev}\circ c_G^{-1})(\phi)\;,\qquad \chi^\phi:=(\mathrm{hom}\circ c_G^{-1})(\phi)\;,\qquad G^\phi := \chi^\phi(G)\;.
\end{eqnarray}
Call a group $G^\phi$ a \emph{boundary group} iff $\phi \in \partial T(G)$ for some other Fuchsian group $G$ of first kind. In \cite{OBoTSaoKG:1}, Bers conjectured the following, which was also the motivation for the name b-group.
\begin{conj} \label{b-conj}
    Every b-group is a boundary group.
\end{conj}
This conjecture is still open in full generality. It is, however, easy to see that boundary groups are b-groups, and we show this in the next lemma.
\begin{lemma} \label{isIso}
   Let $\phi \in \bb S$ and let $G$ be a Fuchsian group of first kind acting on the disc. Then $\chi^\phi$ is an isomorphism      and $G^\phi$ is a b-group with sc-hyp invariant component $f^\phi(\bb D^c)$.
\end{lemma}
\emph{Proof}. Suppose $\chi^\phi$ were not an isomorphism. Then there would be a $g \in G$ such that $\chi^\phi(g)= \bb 1$.
Then $$ f^\phi(w) = \chi^\phi(g) (f^\phi(w)) = (f^\phi \circ g) (w) \forall w \in \bb D^c\;,$$ which cannot be since any nontrivial $g$
fixes at most 3 points and f is injective. The injectivity allows us to write $$\chi^\phi(g) = f^\phi \circ g \circ (f^\phi)^{-1}$$ on
all of $f^\phi(\bb D^c)$, so we see imediately that $\chi^\phi(G)$ acts properly discontinuously on $f^\phi(\bb D^c)$ by
biholomorphisms. It remains to show that every point of $\partial f^\phi(\bb D^c)$ is a limit point of $\chi^\phi(G)$. Suppose
$p \in \partial f^\phi(\bb D^c)$ were not a limit point. Then no sequence $\chi^\phi(g_m)(w)$ for any $w \in f(\bb D^c)$ would
converge to $p$. But this is not true since we can find such a sequence for any point on $\partial \bb D^c$: 
$$ \chi^\phi(g_m) w = f(g_m(z)) \rightarrow p \in \partial \bb D^c\;. \Box$$
We now start on the geometry of $b$-groups, where we rely on the beautiful series of papers \cite{OBoTSaoKG:1}, \cite{OBoTSaoKG:2} and \cite{OBoTSaoKG:3}. The second one is especially important for us, since that paper explains the precise relationship between noded Riemann surfaces and certain groups on the Bers boundary. The relation is very beautiful, but quite tricky and involves many details that we won't be able to present here. We try to give the reader unfamiliar with these matters a guide, though, and will explain the geometric idea.\\

Maskit was able to prove the following classification of b-groups.
\begin{thm}[\cite{OBoTSaoKG:2}, Thm. 4] \label{classbGroups}
   Let $G'$ be a b-group with simply-connected invariant component $D$. Then there are three possibilities, namely
   \begin{enumerate}
   \item{$\Omega(G')$ has precisely two invariant simply-connected components,}
   \item{$\Omega(G') = D$, i.e., G is totally degenerate and}
   \item{$G'$ has \textrm{APT}'s.}
   \end{enumerate}
\end{thm}
This should be compared with another result from \cite{OBoTSaoKG:2} concerning the classification of b-groups. The classification is \emph{not disjoint}: Though a quasi-Fuchsian group is neither totally degenerate nor has APTs, a totally degenerate group \emph{can have} APTs (\cite{OBoTSaoKG:2}, Thm.~ 10).\\

Maskit proves further in (\cite {OBoTSaoKG:2}, Thm.~ 2) that groups satisfying (1) are quasi-Fuchsian, so from a deformation viewpoint they lie inside the Teichm\"uller space of a Fuchsian group $G$. The groups satisfying (2) are a little awkward in terms of their geometric interpretation, since they correspond to the \emph{vanishing} of the Riemann surface that was being deformed: Recall that we are thinking of deforming a Riemann surface $\Sigma = \bb D/G$ and we think of $G'$ as being $G^f$ of some $f \in \hat{\cal T}_S(G)$. The deformed surface is $\Sigma^f = (f(\bb D^c))^c/G^f$ while $f(\bb D^c)/G^f \cong \Sigma^*$ throughout the deformation. In case $\Omega(G')$ only has one component, this component must be $f(\bb D^c)$, so obviously $\Sigma^f=\emptyset$. The existence of such groups \emph{on the boundary} of Teichm\"uller space was proved by Bers in \cite{OBoTSaoKG:1}. \\

Hence we turn our attention to groups with APT's, so for the rest of this section let $G'$ be such a group. This is a good point to take a look back at Propositions \ref{2invcomp} thru \ref{finimp2comp} together with Ahlfors' Finiteness Theorem, Thm. \ref{AFT}. Because $\Omega(G')$ neither consists
of one nor two connected components, but more, the propositions imply that $\Omega(G')$ has \emph{infinitely many} connected components, all of which
are \emph{non-invariant} except for $D$ itself,
$$\Omega(G') = D \cup \bigcup_{j=1}^\infty D_j\;,\qquad \Omega(G')/G' = \Sigma^* \cup  \bigcup_{l=1}^p \Sigma_l\;.$$ By the Finiteness Theorem there are, however, only \emph{finitely many} quotient surfaces, which we denote by
\begin{eqnarray} \label{nRS} \Sigma^* := D/G'\;, \qquad \Sigma_j := \Delta_{j(l)}/G' = D_{j(l)} /G'_{D_{j(l)}}\;, \qquad \out{\mathrm{where} \;} \Delta_{j(l)} = \bigcup_{g\in G'} g(D_{j(l)})\;,\end{eqnarray}
where $j: \{1, \ldots ,p\} \rightarrow \bb N$ is an index function picking a component $D_{j(l)}$ covering the surface $\Sigma_l$. All of these are Riemann surfaces of finite analytic type. \\

Let us recall the definition of an accidental parabolic transformation, or APT for short. A b-group $G$ with invariant simply connected component $D$ admits a Fuchsian equivalent $G^F$, which is defined as
$$G^F:= \psi^{-1} \circ G \circ \psi\;, \qquad \mathrm{where} \; \psi: \bb D \rightarrow D \; \mathrm{is \:any \: Riemann \: mapping}\;.$$
We denote the conjugation map from $G$ to $G^F$ by $\chi^F$. An element $g \in G$ is an APT iff it is parabolic and $\chi^F(g)$ is loxodromic. \\

\begin{figure}[t]
   
    \includegraphics [scale=.8]{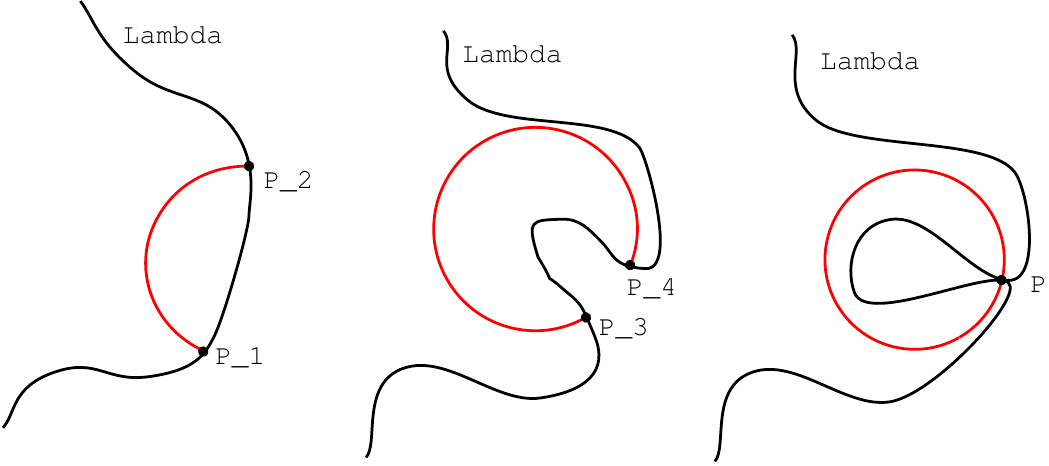}
    \caption{ \label{APTfig}This figure illustrates how a group with an APT can occur as the limit of quasi-Fuchsian groups. $\Lambda$ should be thought of as a piece of the limit set of the group (which, though drawn smooth for the sake of illustration, should rather be thought of as very irregular, since it is almost everywhere not differentiable and of Hausdorff dimension greater than one), and $P_1,P_2$ are the fixed points of the loxodromic element corresponding to the red axis. These two fixed points get mapped to the point $P$ in the third picture, the arc segment closes to a circle and now corresponds to a parabolic element whose Fuchsian equivalent will be loxodromic, i.e., an APT.}
\end{figure}
Let $g$ be an APT. Since a loxodromic element has two fixed points on the limit set, $\chi^F(g)$ has the same, say $p_1,p_2 \in \partial \bb D = \Lambda(G^F)$. Now the single fixed point $p$ of $g$ has to lie on $\Lambda(G)$. Because a homeomorphism from a Jordan domain extends to a continuous map of its closure, such an extension $\overline \psi$ of $\psi$ satisfies $\overline \psi(p_i) = p$. Bers shows in (\cite{OBoTSaoKG:1}, Prop.~ 5.1-5.5) that this is the only thing that can distinguish the Fuchsian equivalent from the group $G$ itself: $\chi^F$ preserves ellipticity and loxodromicity and maps parabolic elements to parabolic or loxodromic elements.\\

An illustration of this can be found in Figure \ref{APTfig}. Note, that this figure should be thought of as infinitesimal, because this process of two points $P_1,P_2$ touching in the limit must happen for infinitely many pairs of points, since we know that $\Omega(G')$ consists of infinitely many connected components, and hence it should appear infinitely often in any `finite size' illustration. \\

The appearance of infinitely many connected components can also be understood very easily with the help of the following lemma. The proof is an easy
consideration involving the trace and the fixed points of $g$.
\begin{lemma}\label{allAPT}
   Let $g \in G'$ be an $\mathrm{APT}$. Then for any $h \in G'$ and $n \in \bb Z$, $hgh^{-1}$ and $g^n$ are again $\mathrm{APT}$'s.
\end{lemma}
For a non-elementary group $G$ the cardinality of $[g]:=\{hgh^{-1}\mid h \in G'\}$ is infinite\footnote{Parabolic elements can only be the same if their fixed points agree. The fixed point of $hgh^{-1}$ is $h(p)$, where $p$ denotes the fixed point of $g$, and since $G'$ is non-elementary there are infinitely many possibilities for $h(p)$.} so each of these produce a new connected component of $\Omega(G')$. The lemma also allows us to introduce the notion of a \emph{basis} of the APT's. Call a loxodromic element $g \in G'$ \emph{primary} iff $g$ cannot be written as $h^n$ for some $h \in G'$ and $|n| \neq 1$.
\begin{dfn}
   Let $G'$ be a b-group. A set of elements of $G'$, $\{g_i\}$, is called a \underline{basis for the $\mathrm{APT}$'s} iff each $g_i$ is primary, $[g_i] \neq [g_j^{\pm 1}]$ and any $\mathrm{APT}$ of $G'$ is contained in some conjugacy class $[g_i^n]$.
\end{dfn}
The last requirement makes sense since by Lemma \ref{allAPT}, every element of $[g_i^n]$ is an APT.\\

Let us denote the map sending any closed curve on a Riemann surface $\Sigma$
passing through $p \in \Sigma$ to its class in the fundamental group $\pi_1(\Sigma,p)$ by $\tilde \tau$. This map induces a map $\tau$ sending free homotopy classes of loops to
conjugacy classes of group elements, in terms of which we can state the following lemma by Maskit.
\begin{lemma}[\cite{OBoTSaoKG:2}, Lemma 2]
        Let $G'$ be a b-group with invariant component $D$, and let $\{g_i\}$ be a basis for the $\mathrm{APT}$'s. Then there is a homotopically independent set of loops $\{\alpha_i\}$ on $\Sigma^*=D/G'$ such that $\tau(\alpha_i) = [g_i]$.
\end{lemma}
In other words, for any b-group with APT's, a basis for these corresponds to a set of loops on the non-deformed surface $\Sigma^*$ which are pinched in the complement, since the part of the axis $A'_{g} \subset D$ is a circle with one point missing and its open complement would have corresponded to the loop on the surface $(\Omega(G')\backslash D)/G'$. Hence, let $\Omega^+(G')$ denote the following set,
$$\Omega^+(G') := (\Omega(G') \backslash D) \cup \{p\in \Lambda(G') \mid \exists \: \mathrm{APT} \: g\in G' \: \mathrm{s.t.}\: g(p) = p\}\;.$$
Now, if we consider the quotient
$$ \Sigma^+ := \Omega^+(G')/G'\;,$$
we described how to understand this as a connected topological space containing all the $\Sigma_l$ from the decomposition \eqref{nRS}, which has been obtained from gluing together some punctures\footnote{More precisely: fill the punctures to obtain marked points and identify these.} of the surfaces $\Sigma_l$. This looks very much like a noded Riemann surface, and indeed is: The paper by Bers already establishes a connection with pinching curves by estimating widths of collars (i.e., traces of loxodromic elements); the complete dictionary of the precise relationship\footnote{This is not at all trivial: replacing curves by a point could a priori result in any kind of isolated singularity, not necessarily of the simple type $zw=0$. Moreover he deals with the case of surfaces of arbitrary signature and proves essentially that any combinatorical possibility can be realized.}, however, was given by Maskit in (\cite{OBoTSaoKG:2}, Theorem 5 and 6). The result is: \emph{Any b-group corresponds to a noded Riemann surface and vice versa}.\\

We have now seen how a b-group can produce noded Riemann surfaces. We have also stated the existence of totally degenerate b-groups, i.e., groups where the deformed surface has vanished. One can of course also think of both phenomena happening at the same time: Having APT's such that the deformed surface consists of several surfaces $\Sigma_j$ but \emph{not all} surfaces resulting from the pinching - some have vanished. To understand this we have to provide yet some more details.\\

How could one detect the vanishing of some parts of the surface? One possibility is by considering the Poincar\'e area of the surfaces involved. Recall that the area of a Riemann surface of finite type is determined by its genus and its number of punctures as follows,
$$ A(\Sigma) = 2\pi(2g-2+n)\;.$$
Hence, for a quasi-Fuchsian group with connected components $D_1$ and $D_2$,
$$ A(\Omega(G)/G) = A(D_1/G) + A(D_2/G)= 2 A(D_i/G)\;,$$
since the two factor surfaces of such a group have same data $(g_i,n_i)$. Let $\Sigma_\gamma$ be a surface obtained by pinching the simple closed geodesic $\gamma$. We have the two cases: Either $\gamma$ is a separating curve, in which case the area is given by
\begin{eqnarray*} A(\Sigma_\gamma) &=& A(\Sigma_1) + A(\Sigma_2) = 2\pi \left((2g_1-2+n_1+1) + (2g_2-2+n_2+1)\right)\\
      &=& 2\pi(2(g_1+g_2) -2 + n_1+n_2) = A(\Sigma)\;,\end{eqnarray*}
or $\gamma$ is non-separating, and then the area is given by
      $$A(\Sigma_\gamma) = 2 \pi (2(g-1)-2+n+2) = A(\Sigma)\;.$$
Hence we see inductively that pinching any set $\cal C$ of closed geodesics does not affect the total area of the surface $\Sigma_{\cal C}$, i.e., $A(\Sigma) = A(\Sigma_{\cal C})$. From a more general point of view, Bers \cite{IfFGKG} proves a quantitively sharpened version of the Ahlfors Finiteness Theorem, implying what are now known as the Bers Inequalities.
\begin{thm}
   Let $G$ be a function group with invariant component $D$. Denote the factor surfaces by $\Sigma^0 = D/G$ and $\Sigma_i$, $i=1\ldots N$, and let $q\geq 2$ be an integer. Then
   \begin{eqnarray} \label{BAI} \sum_{i=1\ldots N} A(\Sigma_i) \leq A(\Sigma^0)\;, \qquad \sum_{i=1\ldots N} \mathrm{dim} B_q(\Sigma_i) \leq \mathrm{dim}     B_q(\Sigma^0)\;.\end{eqnarray}If one of the two inequalities is an equality then so is the other.
\end{thm}
In the cases we described by pinching, the first of the above two inequalities was an equality. On the other hand, the above inequality is not always an equality, e.g., consider $G$ being totally degenerate.
\begin{dfn}
   Let $G'$ be a b-group and denote the factor surfaces as in \eqref{nRS}. $G'$ is called a \underline{regular b-group} iff equality holds in \eqref{BAI}. $G'$ is called \underline{partially degenerate} iff equality does not hold, yet $G'$ is not totally degenerate.
\end{dfn}
We have finally arrived at the desired group-theoretical notion corresponding to the topological operation of pinching curves: Regular b-groups are quasi-Fuchsian groups together with groups containing APT's such that no part of the quotient has disappeared. The paper of Maskit makes the
connection between noded Riemann surfaces and b-groups very precise. Let us summarize some of the facts discussed above, together with parts of
(\cite{OBoTSaoKG:2}, Thm.~ 7) applied to a regular b-group, for later reference in the following Theorem.
\begin{thm}
   Let $G'$ be a regular b-group without elliptic elements and denote the factor surfaces as in \eqref{nRS}. The genus and punctures of $\Sigma^0$ are     denoted by $(g^0,n^0)$ and the same data of $\Sigma_i$ by $(g_i,n_i)$. Let $k$ be the cardinality of a basis of APT's for $G'$. Then
   \begin{itemize}
     \item{each $\Sigma_i$ is hyperbolic,}
     \item{the number of parts is bounded by $p\leq k-1$,}
     \item{the number of punctures is given by $\sum n_i = n^0 + 2k$ and}
     \item{the dimension of bounded differentials is given by
           $$ \sum_{i=1\ldots p} \mathrm{dim} B_q(\Sigma_i) = \mathrm{dim} B_q(\Sigma^0) - k \;, \qquad 2\leq q \in \bb Z\;.$$}
           \end{itemize}
\end{thm}
Now that we understand the different possibilities for b-groups and their geometry, we have to look at their distribution on the boundary of Teichm\"uller space. The next section will then deal with convergence properties, i.e., in which sense these boundary groups are limits of the quasi-Fuchsian deformations of a given group. Let us call a boundary group with APT's a \emph{cusp}.
\begin{thm}[\cite{OBoTSaoKG:1}, Thm. 12 and 13]\label{BoundTeich}
   Let $G$ be Fuchsian of first kind and not a triangle group\footnote{A triangle group is a group, the factor surfaces of which are thrice punctured spheres, for which the Teichm\"uller space is a point.}. Then there are cusps on $\partial T(G)$. Moreover, the set of cusps has positive real codimension in $\partial T(G)$.
\end{thm}
From this it also follows that there not only are totally degenerate groups on $\partial T(G)$, but even that \emph{most} groups on the boundary are totally degenerate. This can be seen in several ways, one of which we describe here and another argument proving this will be provided in Section \ref{DontExtSect}. \\

The present argument will explain why the set of cusps has positive codimension on the boundary: An occurence of an APT in a group can be detected by considering the trace\footnote{Recall that the trace of a loxodromic element is in $\bb C \backslash [0,4]$ while the trace of a parabolic element is precisely 4.} of loxodromic elements during the deformation. The function
$$\mathrm{Tr}_g(\phi) := \mathrm{tr}^2\left((\mathrm{hom} \circ c_G^{-1})(\phi)(g)\right)\;, \qquad g \in G\;,$$
is holomorphic as a function $B_2(\bb D^c,G)\rightarrow \bb C$. The set of cusps (let us denote it with $CT(G)$ for the moment) can then be expressed as
$$ CT(G):= \bigcup_{\substack {g \in G \\ g \; \mathrm{loxodromic}}} \left(\mathrm{Tr}_g^{-1}(4) \cap \partial T(G) \right)\;.$$
This set can be seen as the countable union of subsets of proper analytic subvarieties of $B_2(\bb D^c, G)$. The fact about the positive real codimension follows.\\

The paper \cite{OBoTSaoKG:1} ends with some conjectures, many of which are solved by now. A very well-known one is the following, which in a sense gives the converse statement to the fact that most boundary groups are totally degenerate.
\begin{thm} [McMullen, \cite{CaD}]
   The set $CT(G)$ is dense in $\partial T(G)$.
\end{thm}
McMullen even proves that the set of \emph{maximal cusps}, i.e., regular b-groups on the boundary for which the factor surfaces are all thrice punctured spheres, is dense. A partial result to the conjecture by Bers already mentioned earlier, i.e., whether any b-group is a boundary group, is given by the following theorem.
\begin{thm}[Abikoff, \cite{OBoTSaoKG:3}]
   Every regular b-group is a boundary group.
\end{thm}
In full generality, however, the conjecture is still open and it has deep connections with other very difficult questions in Kleinian group theory and
the geometry of hyperbolic 3-manifolds, e.g., the ending lamination conjecture (see \cite{HMaKG}). In particular, it is not known if every totally
degenerate b-group is a boundary group.
%%%%%%%%%%%%%%%%%%%%%%%%%%%%%%%%%%%%%%%%%%%%%%%%%%%%%%%%%%%%%%%%%%%%%%%%%%%%%%%%%%%%%%%%%%5555
\subsection{Convergence of Kleinian Groups}
We have dealt with the set-theoretic aspects of the Bers boundary of Teichm\"uller space, but we have yet to understand the topological part, i.e., in which sense the objects (groups, domains, surfaces) converge. In this section we want to adress these questions a little bit from the point of view of Kleinian groups. The first part starts with some general notions and quotes some results, while in the second part we specialize these to sequences in Teichm\"uller space converging to the boundary. We won't provide any proofs in the first part, mainly because most proofs use \emph{3-dimensional} techniques\footnote{The action of Kleinian groups on $\cinf$ extends naturally to an action on the 3-dimensional hyperbolic ball $\bb B^3$, and any complete hyperbolic 3-manifold can be written as a quotient $\bb B^3 /\Gamma$ for some Kleinian group $\Gamma$. This perspective introduces powerful new techniques to the study of Kleinian groups, see, e.g., the beautiful book \cite{HMaKG}.}.

Let us begin by defining two notions of convergence for Kleinian groups.
\begin{dfn}
   Let $\chi_i: \Gamma \rightarrow \mathrm{PSL}(2,\bb C)$ be a sequence of homomorphisms of an abstract group $\Gamma$     into $\mathrm{PSL}(2,\bb C)$ and denote the image groups by $G_i$. The sequence is said to converge                     \begin{itemize}
      \item{\underline{algebraically} to $G$ iff $\lim \chi_i(\gamma)$ exists and is a M\"obius transformation for all             $\gamma \in \Gamma$, and any element $g \in G$ can be written as such a limit;}
      \item{\underline{geometrically} to $H$ iff for every subsequence $\{G_{n_i}\} \subset \{G_i\}$,
         $$H=\{g \in \mathrm{PSL}(2,\bb C) \mid g = \lim g_{n_i}\;, \; g_{n_i} \in G_{n_i}\}\;´.$$}
      \end{itemize}
\end{dfn}
The use of $\lim g_n$ is understood in the topology of $\mathrm{PSL}(2,\bb C)$. That the algebraic and geometric limits\footnote{A sidenote on the nomenclature: While algebraic convergence is quite self-explanatory, the reason for calling geometric convergence so stems from the theory of 3-dimensional hyperbolic manifolds. A sequence of Kleinian groups converges geometrically iff their fundamental polyhedra in $\bb H^3$ converge uniformly on compact sets, see \cite{AaGCoKG}, Prop.~ 3.10.} are again groups is easily seen. Moreover, if the sequence consists of finitely generated Kleinian groups, both limits are again Kleinian \cite{AaGCoKG}, but they are in general very different from one another\footnote{The
potential difference of the limits was the key ingredient in the proof by Thurston and Kerckhoff that the action of the modular group does not extend
continuously to the boundary of Teichm\"uller space.}. Still, both notions are somewhat related. In general we have the following.
\begin{prop}[\cite{AaGCoKG}, Prop.~ 3.8 and 3.10]
   Let $\chi_i: \Gamma \rightarrow G_i$ be an algebraically converging sequence converging to $\chi: \Gamma \rightarrow G$.
   Then there exists a geometrically converging subsequence $\{G_{n_i}\}$ converging to a Kleinian group $H$ and $H$ contains $G$ as a subgroup.
\end{prop}
Let us introduce a convergence notion for open sets in a topological space which is very much connected to the present considerations.
\begin{dfn} \label{CaratConv}
   Let $\{D_j\}$ be a sequence of open sets of a topological space $X$. Then the sequence converges \underline{in the sense of Carath\'eodory} to the open set $D \subset X$ iff
   \begin{enumerate}
      \item{for any compact $K \subset D$ there exists an $N \in \bb N$ such that $K \subset D_n$ for $n\geq N$ and}
      \item{any open $U$ with $U \subset D_j$ for infinitely many $j$ satisfies $U\subset D$. }
    \end{enumerate}
\end{dfn}
Of course we are thinking of the domains of discontinuity $\Omega(G_j)$ for a sequence of Kleinian groups. This fits into the picture by giving a
criterion for when the algeraic and geometric limit are the same.
\begin{thm} [\cite{AaGCoKG}, Prop.~ 4.2]
    Let $\rho_i: \Gamma \rightarrow G_i \subset \mathrm{PSL}(2,\bb C)$ be a sequence of isomorphisms converging algebraically to an isomorphism $\rho:\Gamma \rightarrow G$ with $\Omega(G)\neq \emptyset$. Then $\{G_i\}$ converges geometically to $G$ if the domains of discontinuity $\Omega(G_j)$ converge to $\Omega(G)$ in the sense of Carath\'eodory.
\end{thm}
Moreover, if the algebraic limit is \emph{geometrically finite}\footnote{This is a characteristic of Kleinian groups that only becomes visible in the 3-dimensional world: A group is called \emph{geometrically finite} iff any Dirichlet fundamental polyhedron for the action of $G$ on $\bb B^3$ is finite-sided. In contrast to $\cinf$, where a finitely generated group always has a finite-sided hyperbolic polygon as a fundamental domain, this is not true anymore in $\bb B^3$. Almost all technical difficulties in the theory of Kleinian groups are due to geometric infiniteness. Let us only mention that a b-group is geometrically finite iff it is quasi-Fuchsian or regular.} then the above theorem is an if and only if (\cite{AaGCoKG}, Thm.~ 4.8). A sequence converging algebraically and geometrically to the \emph{same limit} is called a \emph{strongly convergent} sequence.\\

After these general considerations, we now consider special sequences of groups, namely sequences induced from a norm-convergent sequence of quadratic differentials. Given a sequence $\{\phi_i\}\in B_2(\bb D^c,G)$, we have the induced sequence of representations $\chi_i:= \left(\mathrm{hom} \circ c_G^{-1}\right)(\phi_i) : G \rightarrow G_i \subset \mathrm{PSL}(2, \bb C)$. We denote the representation induced by the limit $\phi$ by $\chi^\phi$. Sequences contained in $T(G)$ for some $G$ converging to the boundary $\partial T(G)$ are of course our main interest. In particular, such sequences are contained in $\bb S$, which implies that the $\chi_i$ are all isomorphisms (see Lemma \ref{isIso}).
\begin{lemma}
        Let $G$ be a Fuchsian group and $\{\phi_i\} \subset \bb S(G)$ be a sequence converging in the $B_2$-norm to $\phi$. Then the sequence of                induced homomorphisms $\chi_i: G \rightarrow G_i$ converges algebraically to $\chi^\phi$.
\end{lemma}
\emph{Proof}. Because the developing map is schlicht,
$$f_i := \left(\mathrm{hom} \circ c_G^{-1}\right)(\phi_i)\;,$$
they can be represented as $\chi_i(g) = f_i \circ g \circ f_i^{-1}$ when acting on points of $\bb D^c$. $B_2$-convergence of the $\phi_i$ implies local uniform convergence of the $f_i$, and the limit function $f$ is schlicht by closedness of $\bb S(G)$. The corresponding $f$ determines a
deformation $(f,\chi)$, i.e., $\chi(g):= f \circ g \circ f^{-1}$ is a subgroup of $\mathrm{PSL}(2,\bb C)$, which coincides with $\chi^\phi$ and hence $\chi_i: G_i \rightarrow \mathrm{PSL}(2,\bb C)$, converges algebraically to $\chi$ by the local uniform convergence of the developing map. $\Box$\\

We are also able to understand convergence to totally degenerate boundary points, as we shall show next.
\begin{prop}
    Let $G$  be a Fuchsian group of first kind and $\{\phi_i\} \in T(G)$ a sequence converging in norm to $\phi \in \partial T(G)$ such that $G^\phi$ is totally degenerate. Then $\{G_i\}$ converges strongly to $G$.
\end{prop}
\emph{Proof.} The proof is based on showing Carath\'eodory convergence of the domains of discontinuity. First of all, it is a classical theorem in complex analysis that given a locally uniformly convergent sequence $\{h_j\} \in \cal S(\bb D)$ normalized by $h(0)=0$ and $h'(0) = 1$, the domains $h_j(\bb D)$ converge in the sense of Carath\'eodory to $h(\bb D)$. Let $\{f_i\}$ denote the sequence of developing maps of the deformations associated to the $\phi_i$. They converge locally uniformly and are 1-point normalized, hence the classical theorem applies, i.e., $f_j(\bb D)$ converges to $f(\bb D)$ in the sense of Caratheodory. What is left to show is that the complements converge, i.e., $f_i(\bb D)^c \rightarrow f(\bb D)^c = \emptyset$, since $\Omega(G_i) = f_i(\bb D) \cup f_i (\bb D)^c$ and $\Omega(G^\phi) = f^\phi(\bb D)$. The proof concludes with the following more general lemma.
\begin{lemma} \label{CarConvinLoccptSp}
   Let $X$ be a locally compact topological space, $U$ an open set such that $U^c = \emptyset$ and $U_i$ a sequence of open sets converging to $U$ in      the sense of Carath\'eodory. Then $U_i^c$ converge to the empty set in the sense of Carath\'eodory.
\end{lemma}
\emph{Proof.} Before we start with the proof, recall that we denote by $D^c$ the \emph{interior} of the complement of $D$. The complement of a set $D$ in a space is denoted by $D^*$.\\

We have to show that the sequence $U_i^c$ satisfies the two properties of Definition \ref{CaratConv}. The first is trivially satisfied, since the only compact set in the claimed limit set is the empty set itself. The second property states that any open set $V$ contained in infinitely many of the $U^c_i$ is contained in the limit. In our case this amounts to showing that $V = \emptyset$. Now 
$V \subset U_i^c$ implies $V^* \supset {U_i^c}^*$, and this again implies 
\begin{equation}\label{zui}  V^c \supset (U_i^c)^c\;. \end{equation}
Now, for any open set $D$, we have $D \subset(D^c)^c$: Let $p \in D$ be a point; then $p \not \in D^*$ and since for a closed set the closure of the interior is contained in the set itself,
$p \not \in \mathrm{cl}( \mathrm{int}(D^*))$, hence $p \in (\mathrm{cl}( \mathrm{int}(D^*)))^*$, which is an open set. On the other hand, $\mathrm{int}(D^*) \subset \mathrm{cl} (\mathrm{int}(D^*))$ so the reverse inclusion holds for the complements,
$$ (\mathrm{int}(D^*))^* \supset (\mathrm{cl} (\mathrm{int}(D^*)))^*\;.$$
Now since the interior is the largest open set contained in a closed set and the right hand side is open, this implies 
$$ (D^c)^c = \mathrm{int}(\mathrm{int}(D^*))^* \supset (\mathrm{cl} (\mathrm{int}(D^*)))^* \supset D\;.$$
If we use this fact, \eqref{zui} implies $V^c \supset U_i$. Now we use Carath\'eodory convergence of the $U_i$. Because $U \cap V \neq \emptyset$ by assumption and $X$ is locally compact, there exists a $K \subset U \cap V$. Since $K \subset U$, the first property of Caratheodory convergence implies that $K \subset U_i$ for $i$ big enough. But this contradicts the fact that $K \subset V$ because $V^c \subset U_i$. $\Box\: \Box$\\

The lemma we used in the proof is false in general: The Carath\'eodory convergence of a sequence of open sets to a limit set does not imply the Carath\'eodory convergence of the complements to the complement of the limit set. It was crucial that the limit set was $\emptyset$.\\

And in fact, this happens on boundaries of Teichm\"uller spaces: For any Fuchsian group $G$ such that $T(G)$ has dimension greater than one, there exist sequences $\phi_j \in T(G)$ converging to a boundary point $\phi \in \partial T(G)$ such that the algebraic limit and the geometric limit\footnote{Recall that an algebraically converging sequence always has a geometrically converging subsequence, so we can assume without loss of generality that it converges geometrically, too.} \emph{do not coincide}. A beautiful construction of such sequences can be found in \cite{NotAotMGaBBoTS}. \\

We note here that the proof that the modular group action on $T(G)$ \emph{never} extends to a continuous action on the Bers completion $\overline{T(G)}$ given in \cite{NotAotMGaBBoTS} uses the discrepancy between the algebraic and the geometric limit, which necessarily occurs somewhere.\\

For analytic considerations, it is often sufficient to be able to handle the convergence of specially chosen sequences. Bers introduced yet another convergence concept in \cite{TAotMGotCB}, which he called \emph{tame convergence}.
\begin{dfn}
   Let $\phi, \phi_j \in B_2(\bb D^c)$. The sequence $\{\phi_j\}$ is said to converge
   \begin{itemize}
      \item{\underline{weakly} to $\phi$ iff the norms $\|\phi_j\|_{B_2(\bb D^c)}$ are uniformly bounded and if the functions converge pointwise, i.e., $\lim \phi_j(z) =\phi(z)$ for all $z \in \bb D^c$.}
      \item{\underline{tamely} to $\phi$ iff $\phi, \phi_j \in \bb S\subset B_2(\bb D^c)$, the seqence converges weakly to $\phi$ and 
      for almost all $z \not \in f^\phi(\bb D^c)$, there exists an integer $J$ such that $z \not \in f^{\phi_j}(\bb D^c)$ for all $j\geq J$.} 
	\end{itemize}
\end{dfn}
It is easy to see that norm convergence implies weak convergence. Moreover, any bounded sequence contains a weakly convergent subsequence, because one can choose a countable dense set of points in $\bb D^c$, where the values of the functions are necessarily bounded. Choose a convergent subsequence for every point and use a diagonal argument to obtain a sequence converging at all the points of the dense set. By holomorphicity, the pointwise convergence then holds at all points. \\

Bers proved the following two facts about tame convergence.
\begin{prop}[\cite{TAotMGotCB}, Prop.~ 2.1]
   If $\phi \in \bb S$ and $\mathrm{mes}\:\partial f^\phi(\bb D^c) = 0$, then there exists a sequence $\{\phi_j\} \subset T(\bb 1)$ such that $\phi_j \rightarrow \phi$ tamely.
\end{prop}
This does not imply that any element in $\partial T(G)$ is the tame limit of a sequence in $T(G)$, of course. But an important special case does hold.
\begin{prop}[\cite{TAotMGotCB}, Prop.~ 2.2] \label{StrongConvPoss}
   Let $G$ be finitely generated of first kind and $\phi\in \partial T(G)$ a regular boundary point (i.e., $G^\phi$ is a regular b-group). Then there exists a sequence in $T(G)$ converging tamely to $\phi$.
\end{prop}
Kerckhoff and Thurston already observed that for a sequence in $T(G)$, tame convergence is equivalent to strong convergence of the corresponding groups. This means that for \emph{any regular boundary} point in $T(G)$ for $G$ Fuchsian of first kind, there exists a sequence in $T(G)$ such that the corresponding groups converge both algebraically and geometrically to the same limit, and hence, the domains of discontinuity also converge in the sense of Carath\'eodory. But, as argued above, this is not true for \emph{any} sequence in $T(G)$ converging (in norm) to that boundary point.

However, one can show that as long as one stays inside Teichm\"uller space, i.e., the limit $\phi$ of the converging sequence $\{\phi_j\}$ is again in $T(G)$, the associated sequence of groups converges strongly to the group $G^\phi$. More details on these subtle matters can be found in the book \cite{HMaKG}. 
%%%%%%%%%%%%%%%%%%%%%%%%%%%%%%%%%%%%%%%%%%%%%%%%%%%%%%%%%%%%%%%%%%%%%%%%%%%%%%%%%%%%%%%%%%%%%%%%%55555
\subsection{The Augmented Teichm\"uller Space}\label{AugmentedTS}
In this section we will finally introduce the relevant partial completion of Teichm\"uller space. From our thorough considerations in the previous sections, we have an idea of what the correct boundary points are - they should correspond to Riemann surfaces with nodes. Since the deformed Riemann surface is obtained as a quotient of $(f(\bb D^c))^c$ by the group $G^f$, it is sensible to only include \emph{regular b-groups} in the definition, thereby excluding the partially and totally degenerate ones, since for the latter some part of the deformed surface has vanished, as we have explained in Section \ref{BersBound}.
\begin{dfn}
   The \underline{Bers augmented Teichm\"uller space}, denoted by  $\hat T^B(G)$, is the following subspace of the Bers completion $\overline{T}(G)$,
   $$ \hat T^B(G) := \{ \phi \in T(G) \cup \partial T(G):\; G^\phi \mathrm{\;is\; a \; regular\; b-group}\}\;.$$
\end{dfn}
Obviously $T(G) \subset \hat T^B(G)$ since for all interior $\phi$, $G^\phi$ is quasi-Fuchsian. This space now handles the codimension problem of the boundary because of Theorem \ref{BoundTeich}. However, this does not yet circumvent the problem
described above regarding the non-continuous extension of base-point-change or even only of the $\cal{MCG}(G)$-action\footnote{In a sense, it makes the problem worse: The
mapping class group action did have a continuous extension to all totally degenerate points, i.e., the ones we got rid of.}. Abikoff resolved this
problem by defining a different topology on $\hat T^B(G)$. \\

We already explained in Section \ref{BersBound} that quotients of b-groups with APT's are precisely noded Riemann surfaces. Hence we see that
set-theoretically we have a bijection $$c: \mathrm{Def}(\Sigma_{\cal L}) \rightarrow \hat T^B(G)\;,$$ where $\Sigma$ is a Riemann surface of finite
type and is related to $G$ as always, $\Sigma = \bb D/G$, and $\Sigma_{\cal L}$ is a maximally noded Riemann surface obtained from $\Sigma$ by
pinching. If we restrict $c$ to $\mathrm{Teich}(\Sigma)\subset \mathrm{Def}(\Sigma_{\cal L})$, it coincides with the recipe we have described in
earlier sections on how to get different representations of Teichm\"uller space, and is therefore a homeomorphism. Thus it allows us to \emph{define}
a new topology on the set $\hat T^B(G)$ which is compatible with the subspace topology of $T(G)$ by declaring the map $c$ to be a homeomorphism. The topology so obtained is called the \emph{geodesic topology}.
\begin{dfn}
   Let $G$ be a Fuchsian group of first kind. The underlying set of $\hat T^B(G)$ equipped with the geodesic topology is called the \underline{augmented Teichm\"uller space} of $G$ and is denoted by $\hat T(G)$.
\end{dfn}
Note that while the Bers augmented Teichm\"uller space is defined for any Fuchsian group acting on the disc, we only define the augmented Teichm\"uller space for groups of first kind. This is convenient, since we only want to study  degeneration phenomena for surfaces of finite type.\\

Abikoff defined two additional topologies, which he called the \emph{conformal} and the \emph{isometric} topology. We will not give the formal definition here, but state the result Abikoff obtained in \cite{DFoRS}.
\begin{thm}
    The three topologies (geodesic, conformal, isometric) on $\hat T(G)$ are the same.
\end{thm}
He went on to prove the following result, which confirms that $\hat T(G)$ is an adequate object to study, since as a completion it satisfies the desired properties.
\begin{thm}
    The action of $\cal{MCG}(G)$ on $T(G)$ extends to a proper discontinuous action on $\hat T(G)$. Moreover, for the quotient one obtains
    $$\hat T(G) / \cal {MCG}(G) \cong \overline{\cal M}(\Sigma)_\mathrm{DM}\;.$$
\end{thm}
We will need one fact about the two topologies later on.
\begin{lemma}\label{RelationOfConvergence}
    Let $\{\phi_j\}$ be a sequence in $T(G)$ converging to $\phi \in \hat T(G)$. Then $\phi_j \rightarrow \phi$ in $\hat T^B(G)$, too. 
\end{lemma} 
\emph{Proof}. If $\phi \in T(G)$ there is nothing to prove, since we already explained above that the topologies agree on $T(G)$. Hence we assume $\phi \in \partial T(G)$. Hence in any neighborhood $U$ of $\phi$ there are infinitely many $\phi_j$, and since $\phi \in \overline{T}(G)$, the $\phi_j$ constitute a Cauchy sequence with respect to the $B_2$-norm. Let $\tilde \phi \in \overline{T}(G)$ be the limit of this sequence. It is now easy to see that $\tilde \phi = \phi$ because the topologies in the interior, i.e., on $U\cap T(G)$ agree. $\Box$\\

On the other hand, not every norm-converging sequence converges geodesically, of course.
%%%%%%%%%%%%%%%%%%%%%%%%%%%%%%%%%%%%%%%%%%%%%%%%%%%%%%%%%%%%%%%%%%%%%%%%%%%%%%%%%%%%%%%%%%%%%%%%%%%%%%%%%%%%%%%%%%%555
\addtocontents{toc}{\vspace{.3cm}}
\section{Extensions of Automorphic Forms over Augmented Teichm\"uller Space}\label{ExtoATS}
We want to extend our considerations about families of automorphic forms over Teichm\"uller space to families where the base is the augmented
Teichm\"uller space $\hat T(G)$ of a Fuchsian group $G$ of first kind. To make sense of automorphic forms over boundary points of $\hat T(G)$, we
first have to extend the fibre space $\cal F(G)$ to a fibre space over augmented Teichm\"uller space, which we do in the first section.\\

Our next goal is to extend the sections continuously to the boundary points.
\subsection{Fibre Spaces over Augmented Teichm\"uller Space}\label{FoATS}
Since the correspondence $c_G$ described in Section \ref{BersBound} relating the schlicht model of Teichm\"uller space to $T(G)$ is defined on all of
$\mathrm{Def}(G)$, we can define the set
$$ \hat{\cal T}_S(G):= \mathrm{dev}\left(c_G^{-1}(\hat T(G))\right) \subset \cal S^0(\bb D^c)\;,$$
which we will call the \emph{schlicht model of augmented Teichm\"uller space}\footnote{Hidalgo and Vasil'ev have, very recently, shown in \cite{NTS} how to generalize the notion of Beltrami differentials to incorporate noded Riemann surfaces as solutions to the Beltrami equation. Hence, one can also define the \emph{Beltrami model of augmented Teichm\"uller space}. This is conceptually very interesting and may lead to new insights on augmented Teichm\"uller spaces.}. With the help of the schlicht model, we can extend the fibre space
$\cal F(G)$ introduced in Section \ref{FibreoverTeich} to a fibre space over augmented Teichm\"uller space in the following way,
$$ \hat {\cal F}(G) \subset \hat{\cal T}_S(G) \times \bb C\;, \qquad \hat {\cal F}(G):= \left \{ (f,z): z \in     \Omega(G^f) \backslash f(\bb D^c) \right\}\;.$$
This is a fibre space, since for points $f$ of augmented Teichm\"uller space the groups $G^f$ are not totally degenerate, i.e., $\Omega(G^f) \neq f(\bb D^c)$. The natural projection to the base will be denoted by
$$ \pi_{\hat {\cal F}}: \hat{\cal F}(G) \rightarrow  \hat{\cal T}_S(G)\;.$$
By the symbol $\hat F(G)$, we will denote the result applying the Bers embedding to the first factor of the product $\hat{\cal F}(G) \subset
\hat{\cal T}_S(G) \times \bb C$, which of course again is a fibre space $\pi_{\hat F}:\hat F(G) \rightarrow \hat T(G)$.\\

Observe that these fibre spaces contain the former fibre spaces over Teichm\"uller space and that the projections agree\footnote{This is a good point
to make a comment on the notation: We have introduced several symbols for the different representations of Teichm\"uller space, but we will not do
this for the fibre space $\cal F$ or its augmentation. Of course from a purist point of view, \eqref{same} is not quite true, since $\cal F(G)$ had
$\cal T_B(G)$, and not $\cal T_S(G)$, as its base. In any case, the fibre above a point is the quasidisc (or the union of such) corresponding to the point
in Teichm\"uller space and is independent of how we think of Teichm\"uller space itself. Since the different representations of Teichm\"uller space are
canonically related, no confusion should arise.} on this subset,
\begin{eqnarray} \label{same}
   \hat{\cal F}(G) \cap \pi_{\hat {\cal F}}^{-1}(\cal T_S(G)) = \cal F(G)\;, \qquad \pi_{\hat{\cal F}}\big|_{\cal F(G)}= \pi_{\cal F}\:,
\end{eqnarray}
since $\Omega(G^f) \backslash f(\bb D^c) = w^f(\bb D)$ if $f \in \cal T_S(G)$ and $w^f$ denotes any quasiconformal extension of $f$ to $\cinf$.\\

We know that fibres of the space $\cal F(G)$ are quasi-discs, so in particular they are Jordan domains. We need a topological description of the new fibres,
which we give in the following lemma for convenience, though it contains no new facts.
\begin{lemma}\label{NewFibres}
   $\hat F(G)$ is a connected subset of $\bb C^{3g-2+n}$, where $(g,n)$ is the type of the surface $\bb D/G$. Any fibre $D_p$ over a point $p \in \hat T(G)$ is a bounded open set contained in $\bb D_4$ consisting of either one simply-connected component, in which case $p\in T(G)$, or infinitely many simply connected, components when $p \in \hat T(G) \backslash T(G)$.
\end{lemma}
\emph{Proof.} 
$\hat F(G)\subset \bb C^{3g-2+n}$ is obvious. The connectedness of $\hat F(G)$ follows from Proposition \ref{StrongConvPoss}, which implies that any boundary point can be approached by a sequence such that the domains of discontinuity of the associated groups converge in the sense of Carath\'eodory. Since the fibres are locally compact, any point $p$ in a boundary fibre has a compact neighborhood $K$ in the fibre, and by the first property of Carath\'eodory convergence, $K$ is already contained in the fibres over interior points nearby. Therefore, there is a path from $p$ to the interior $\cal F(G)$, which obviously is connected.\\   

The boundedness of $D_p$ follows from the Koebe-1/4 Theorem, since $\hat T(G) \subset \bb S$. The other statements were already shown in Section \ref{NodedRS}. $\Box$
%%%%%%%%%%%%%%%%%%%%%%%%%%%%%%%%%%%%%%%%%%%%%%%%%%%%%%%%%%%%%%%%%%%%%%%%%%5555
\subsection{Continuous Extension of Sections}\label{DontExtSect}
We want to tackle the problem of extending the vector bundle $\cal V(\cal R)$ and its sections to the augmented setting. For this, let $t$ be a tuple of Bers coordinates, i.e., given a base $\{\varphi_i\}$ of $B_2(\bb D^c, G)$, the $t_i$ are the coefficients of a given $\phi \in B_2(\bb D^c, G)$ with respect to this basis. They obvioulsy provide global holomorphic coordinates for $T(G)$.\\

Recall that the vector bundle $\cal V(\cal R)$ was given\footnote{Here again and
in what follows, we switch freely between the different representations of Teichm\"uller spaces as the base of the fibre spaces.} as the space
$$\pi_{\cal R}: \cal V(\cal R):=\bigcup_{t \in T(G)} A^1_{\rho_s^t}(D^t, G^t)
\rightarrow T(G) \;,$$ where $\cal R$ was the \emph{unique holomorphic extension} of a given s-factor $\rho_s:\bb D \times G \rightarrow \bb C^*$ on
the 0-fibre of $\cal
F(G)$. The restriction  of $\cal R$ to an arbitrary fibre $D^t$ over $t$ was denoted by $\rho_s^t$. \\

The expressions $g^t \in G^t$ extend continuously\footnote{In fact, the correspondence between the
deformation space of groups $\mathrm{Def}(G)$ and $B_2(\bb D^c,G)$ via the Schwarzian derivative does not change the regularity of any auxiliary
parameter. Hence $g^t$ is holomorphic in $t$ for all $t \in \bb C^N$, where $N$ is the dimension of $B_2(\bb D^c,G)$.} to all $t \in \hat T(G)$.
Therefore there exists a continuous extension of $\cal R$ to a function
$$ \hat {\cal R} : \hat F(G) \times G \rightarrow \bb C^*\:,$$
and the restriction to any fibre $\Omega(G^t)\backslash f^t(\bb D^c)$ yields an s-factor of automorphy $\rho_s^t$. This holds because the only allowed $t$-dependence of
$\cal R$ is the implicit dependence via the holomorphic $t$-dependence of the group elements $g^t$, and the equation for a factor of automorphy
follows from the chain rule of differentiation. Observe that although the new fibres are not simply connected  anymore, no ambiguity arises for non-integer $s$ from more choices of logarithms.
The choice of logarithm in the expression $g'^s$ on the 0-fibre extends in a unique way to all fibres by continuity. \\

Hence the extension of the vector bundle can be defined easily as a fibre space,
$$ \pi_{\hat{\cal R}}: \hat{\cal V}(\cal R):=\bigcup_{t \in \hat T(G)} A^1_{\rho_s^t}(\Omega(G^t)\backslash f^t(\bb D^c), G^t) \rightarrow \hat T(G)\;.$$
This is not quite a vector bundle anymore, because the dimension of the fibres over points on the boundary $\partial T(G)$ does not agree with the dimension of the spaces over interior points. This can be seen as follows: The dimension of $A^1_{\rho_s^t}(\Omega(G^t)\backslash f^t(\bb D^c), G^t)$ for $t \in T(G)$ is given by the Riemann-Roch formula (Thm. \ref{R-R}),
\begin{eqnarray} \label{LocalDim} \mathrm{dim}_{\bb C} A^1_{\rho_s^t}(\Omega(G^t)\backslash f^t(\bb D^c), G^t) = (2s-1)(g-1) + n[s]=:N_s(g,n)\;.
\end{eqnarray}
Now let $t$ be a boundary point. Then 
$$\Omega(G^t)/G^t = \Sigma^0 \cup\left(\Sigma_1 \cup \ldots\cup \Sigma_p\right)\;,$$
and hence 
$$\left(\Omega(G^t) \backslash f^t(\bb D^c)\right)/G^t =\Sigma_1 \cup \ldots\cup \Sigma_p\;.$$
By factoring the group action into orbits (for the notation, see Section \ref{KGandU}) we see that
$$ A^1_{\rho_s^t}(\Omega(G^t)\backslash f^t(\bb D^c), G^t) \cong \bigoplus_{k=1\ldots p} A^1_{\rho_s^t}(D_k, G_{D_k})\;,$$
i.e., we obtain the sum of the induced automorphic forms on the parts of the noded surface from which the nodes have been removed.
\begin{lemma}
  Let $N_s(g,n)$ be as in \eqref{LocalDim} and $s\geq 2$ an integer. For a boundary point $t\in\hat T(G)$, let $\Sigma^0, \Sigma_i$ be the parts   of the quotient	surface and let $(g^0,n^0)$ resp.~ $(g_i,n_i)$ their types.
  Then 
  \begin{equation}\label{TrueAlways} \mathrm{dim}_{\bb C} A^1_{\rho_s^t}(\Omega(G^t)\backslash f^t(\bb D^c), G^t) = N_s(g,n) - \frac12\left(\sum_{i=1\ldots p} n_i - n\right)\;. \end{equation}
\end{lemma}
\emph{Proof}.
We prove this by induction. Observe that the additional term
$$\frac12\left(\sum_{i=1\ldots p} n_i - n\right)$$
is precisely the number of curves on $\Sigma = \bb D/G$ which have been pinched to nodes and then removed to obtain $\Sigma_1 \cup\ldots\cup\Sigma_p$, and hence the statement is true for no pinched curves. As we know, there are two different types of pinching. First, if we pinch a non-separating curve, we produce a surface of genus $g^0-1$ with $n^0 + 2$ punctures. Hence
\begin{equation*}
   \begin{aligned} N_s(g^0-1,n^0+2) &= (2s-1)(g^0-2) + (n_0 + 2)[s]\\  &= (2s-1)(g^0-1) + n^0[s] + 2([s] - s) +1 \\&= N_s(g^0,n^0)-1\;,
   \end{aligned}
\end{equation*}  
because for integer $s$, $[s] = s-1$. On the other hand, if we pinch a separating curve, we produce two surfaces with $g_1 + g_2 = g^0$ and $n_1+n_2 = n^0 + 2$. So
\begin{equation*}
   \begin{aligned}
   N_s(g_1,n_1) +  N_s(g_1,n_1) &= (2s-1)(g_1 + g_2 -2) + (n_1 + n_2)[s]\\  &= (2s-1)(g^0-1) + n^0[s] + 2([s] - s) +1 \\&= N_s(g^0,n^0)-1\;.
   \end{aligned}
\end{equation*}
Hence the claimed formula is true for pinching one curve and it did not depend on the specific starting data $(g^0,n^0)$. By induction, \eqref{TrueAlways} is true for any set $\cal C \subset \cal L$ of pinched curves. $\Box$\\

The analogous statement for non-integer $s$ is a little more sublte to obtain. Let us explain this in the case of pinching a non-separating curve. Observe that the computation above is not valid for $s \not \in \bb Z$, because if $N_s(g^0-1,n^0+2)$ is an integer, $N_s(g^0,n^0)-1$ isn't, since they differ by the term $2([s] - s)$ which is not an integer. In case $s \in \bb Z[(2(g^0 - 1)-2)]$, the space of automorphic forms for this factor still is in a natural correspondence to the space of sections of some holomorphic line bundle on the compactification of the surface of type $(g^0-1,n^0+2)$, but still $N_s(g^0,n^0)$ is not an integer. In this case we have to rewrite Riemann-Roch in terms of the degrees of the bundles first\footnote{Recall that if the surface is non-compact, the value of $s$ of an $s$-factor does not have a geometric meaning and so there are equivalent factors with different values of $s$. The induced $s$-factor for the same value of $s$ does not correspond to the canonical bundle on hte compactification of the surface of type $(g^0,n^0)$, but this is implicitly used in the Riemann-Roch formula.} and then compute the result. In this case, the dimension of the spaces will again differ by 1. We leave the details to the reader. 
%%%%%%%%%%%%%%%%%%%%%%%%%%%%%%%%%%%%%%%%%%%%%%%%%%%%%%%%%%%%%%%%%%%%%%%%%%%%%%%%%%%%%%%%%%%
\begin{thm} \label{ExttoBound}
    The sections of $\cal V(\cal R)$ constructed in \ref{ExtendSections} extend continuously to sections of         $\hat{\cal V}(\cal R)$.
\end{thm} 
\emph{Proof}. Let us recall how we the sections of $\cal V(\cal R)$ were consturcted: Let $f \in \cal O(\bb D_{4+\epsilon})$ be any function. Apply the $\Theta$-operator fibrewise with respect to the induced factor of automorphy on the fibre: This is well-defined, since $f \in \cal
O(\bb D_{4+\epsilon})$ implies that $f$ is bounded on all fibres $D^t$ and since $s\geq 2$, Lemma \ref{TwoMinusS} then implies that $f \in A_{\rho_s^t}^1(D^t,G^t)$.\\
This extends word by word,  
   $$ \Theta_{\hat{\cal R}}[f]: \hat T(G) \rightarrow \hat{\cal V}(\cal R)\:, \qquad \Theta_{\hat{\cal R}}[f](t,z) := \sum_{g \in G^t} \rho_{g^t}^{-1}(z)(f\circ g^t)(z)\;,$$
Continuity follows, because convergence in the augmented topology of course implies convergence in the Bers topology (i.e., convergence in $\hat T(G)$ implies convergence in $\hat T^B(G)$), which we showed in Lemma \ref{RelationOfConvergence} and the construction above is continuous with respect to $t$ with $t$ viewed as a function in $\bb C^{3g-3+n}$, which induces the Bers topology. $\Box$\\

A very interesting thing happens now: For a tuple of functions we defined the Wronskian in Section \ref{PolAB}. Now let $\cal R$ be the family of factors of automorphy corresponding to the factor $\rho_s$ on the 0-fibre with $s\geq 2$. Let $\hat \Psi_1, \ldots, \hat \Psi_N$ be the continuous extensions of $N$ sections of $\cal V(\cal R)$ such that they induce a basis of the fibres in a neighborhood of the 0-fibre as constructed in Theorem (Thm. \ref{BasisAE}). Then they induce a basis of the fibres in all fibres over $T(G) \backslash N(\Psi)$ where $N(\Psi)$ was the vanishing locus of the Wronskian of the tuple $\Psi_i$.\\

The Wronskian of the tuple $\hat \Psi_i$ is a continuous extension of the Wronskian of the $\Psi_i$. Now by dimensional reasons it \emph{must vanish} on all of $\hat T(G) \backslash T(G)$, because the vector spaces over boundary points are of dimension $<N$. On the other hand, let $p \in \bb P^1(\bb C^{3g-2+n})$ be an arbitrary complex line, and define $\hat T_p(G):= \hat T(G) \cap p$. The Wronskian restricted to $T_p(G)$ is a holomorphic function of one variable. Hence if $(\hat T(G) \backslash T(G))\cap p$ were dense, the Wronskian would be identically zero on $T_p(G)$. But generically, this is \emph{not the case}, as we know. The same argument works of course for any 1-dimensional analytic subset, and this explains the analytic subsets $N(\Psi)$ we defined in Section \ref{PolAB}: They are precisely those analytic subsets for which the intersection with $\partial T(G)$ contains a \emph{dense subset of regular boundary points}.\\

It should also be noted that the extended sections $\Theta_{\hat {\cal R}}[f]$ lie in a fibre subspace $\cal R \hat{\cal V}\subset \hat{\cal V}$. The fibres only differ over boundary points. More precisely, $\pi_{\hat{\cal R}}^{-1}(t) \cap \cal R \hat{\cal V}\subset \hat{\cal V}$ is given by tuples of differentials on the parts of the surface obtained by removing the nodes which have \emph{related residue} at two punctures which formerly belonged to the same node. Let us be more precise about this point: The residue is a local quantity, so we restrict ourselves to the local picture of degeneration, i.e., a cylinder degenerating to a pair of punctured discs as described in Section \ref{Strat}. The residue is then obtained as an integral along a curve around the puncture, which we choose without loss of generality to be a circle of small radius, of the local representative of the section. Now if we open the node, the value of the integral changes continuously, since the section is continuous over $\hat T(G)$. The integral of course only depends on the free homotopy class of the curve, and after opening the node the curve is freely homotopic to a circle in the other disc. By pinching again, we see that the residues are related. And what do we mean precisely by related? Well, the local functions in this picture are solutions to an $s$-automorphy equation, and so the transition function contains via the choice of a $(2g-2)$nd root a $(2g-2)$nd root of unity as ambiguity, which depends of the data of the bundle. Hence related means $\mathrm{res}_{p_1} \phi_1(0)= c \cdot \mathrm{res}_{p_2} \phi_2(0)$ where $c \in \exp (2 \pi i \bb Z[(2g-2)^{-1}]) = \exp (\pi i \bb Z[(g-1)^{-1}])$. For example, for integer $s$, the value of $c$ is 1 for even $s$ and $-1$ for odd $s$, and elements of a fibre of $\cal R \hat{\cal V}\subset \hat{\cal V}$ are in this case called \emph{regular s-differentials} (e.g., in \cite{FDTSaG} or \cite{SoDRS}). \\ 

Over a boundary point, the fibre of $\cal R \hat{\cal V}\subset \hat{\cal V}$ is of codimension $k$ in the fibre of $\hat {\cal V}$ because for each node, having related residue as discussed above is a linear condition.  
%%%%%%%%%%%%%%%%%%%%%%%%%%%%%%%%%%%%%%%%%%%%%%%%%%%%%%%%%%%%%%%%%%%%%%%%%%%%%%%%%%%%%%555
\subsection{Asymptotic Behaviour of the Constructed Sections}\label{AsofSections}
Let us look at the asymptotic behaviour of products of sections $\Psi_i$ constructed with the help of the varying $\Theta$-operator like in Proposition \ref{ReqFamTheta}, i.e., $\Psi_i$ is given by
\begin{equation}
   \Psi_i = \Theta_{\cal R}[h_i]\;, \qquad h_i \in \cal O(\bb D_{4+\epsilon})\;,
\end{equation}
where $\cal R$ is the family of $s$-factors associated to the factor $\rho_s$ on the distinguished fibre with real $s\geq 2$ and let $t$ denote some coordinates in Teichm\"uller space, as usual. We consider the $t$-dependence of the hermitean product $\langle \Psi_1, \Psi_2\rangle^{G^t}_s$,
\begin{equation*}
\begin{aligned}
   \left|\langle \Psi_1, \Psi_2\rangle^{G^t}_s\right|&= \left|\int_{\cal F^t} \Theta_{\rho^t_s}[h_1](z) \overline{\Theta_{\rho^t_s}[h_2]} \lambda_{D^t}^{2-2s} d^2z\right|
   =\left|\int_{D^t} \Theta_{\rho^t_s}[h_1] \overline{h_2(z)} \lambda^{2-2s}_{D^t} d^2z\right|\\
   &\leq \| \Theta_{\rho^t_s}[h_1] \|_{B_s(\bb D^t)} \int_{D^t} |h_2(z)| \lambda^{2-s}_{D^t} d^2z
   = \| \Theta_{\rho^t_s}[h_1] \|_{B_s(\bb D^t)} \|h_2\|_{A^1_{s}(D^t)}\;,
\end{aligned}
\end{equation*}
where the first step in the computation follows from Lemma \ref{scalar}. By the equivalence of norms given in Theorem \ref{incl} and the fact that the norm of the $\Theta$-operator is at most one, we obtain
$$ \| \Theta_{\rho^t_s}[h_1] \|_{B_s(\bb D^t)} \leq M(t) \|\Theta_{\rho^t_s}[h_1] \|_{A^1_s(\bb D^t)}\leq M \|h_1\|_{A_s^1(\bb D^t)}\;.$$
Altogether we arrive at
\begin{equation}
       \left|\langle \Psi_1, \Psi_2\rangle^{G^t}_s\right| \leq M(t)  \|h_1\|_{A_s^1( D^t)}  \|h_2\|_{A^1_{s}(D^t)}\;,
\end{equation}
and by assumption we also know that the $A^1_s$-norms of the functios $h_i$ are uniformly bounded. Hence we arrive at the conclusion that all the asymptotic behaviout lies in the behaviour the $t$-dependent constant $M(t)$.\\

Previously, in Equation \eqref{EmbeddingConstant} we obtained in the following bound for this constant,
\begin{eqnarray} M \leq \sup_{z \in \bb D} \left(\lambda_{\bb D}^{-2s}(z)\alpha_s(z,z)\right)\;,\end{eqnarray}
where the function $\alpha_s(z,w)$ is given by applying the $\Theta$-operator associated to the canonical $s$-factor for the Fuchsian equivalent $G^F_t$ of $G^t$ to the automorphic Bergman kernel $K_{\bb D,s}(z,w)$ in the first variable. The automorphic Bergman kernel is given by 
$$ K_{\bb D,s}(z,w) = \frac{(2s-1)}{\pi} (1-z\bar w)^{-s}\;,$$ 
and hence the 
$$ M(t) \leq \sup_{z \in \bb D} \sum_{g_t \in G^F_t} |g_t'|^{-s} \frac{(1-|z|^2)^{2s}}{(1-g_t(z) \bar z)^{s}}\;.$$
The behaviour of this expression under the deformation of a loxodromic transformation $g \in G$ to an accidental parabolic transformation in the Fuchsian equivalent of the boundary group contains all the asymptotic behaviour. This, however, is not so easy to understand in a quantitative manner. We hope to come back to this matter in the future.
\addtocontents{toc}{\vspace{.5cm}}
%%%%%%%%%%%%%%%%%%%%%%%%%%%%%%%%%%%%%%%%%%%%%%%%%%%%%%%%%%%%%%%%%5
%%%%%%%%%%%%%%%%%%%%%%%%%%%%%%%%%%%%%%%%%%%%%%%%%%%%%%%%%%%%%%%%%5
%%%%%%%%%%%%%%%%%%%%%%%%%%%%%%%%%%%%%%%%%%%%%%%%%%%%%%%%%%%%%%%%%5
\newpage
%%%%%%%%%%%%%%%%%%%%%%%%%%%%%%%%%%%%%%%%%%%%%%%%%%%%%%%%%%%%%%%%%5
%%%%%%%%%%%%%%%%%%%%%%%%%%%%%%%%%%%%%%%%%%%%%%%%%%%%%%%%%%%%%%%%%5

\end{document}